\documentclass[aop,preprint]{imsart}
\usepackage{bm}
\usepackage{mathrsfs}
\usepackage{amsfonts,amsmath, amssymb,amsthm,amscd, mathtools}
\usepackage[colorlinks,citecolor=blue,urlcolor=blue]{hyperref}
\usepackage[numbers]{natbib}
\usepackage{verbatim}
\usepackage{graphicx}
\usepackage[latin1]{inputenc}
\usepackage{latexsym}
\usepackage{lscape}
\usepackage{epsfig}
\usepackage{xcolor}

\mathtoolsset{showonlyrefs}

\startlocaldefs

\newcommand{\e}[0]{\epsilon}

\newcommand{\EE}{\ensuremath{\mathbb{E}}}
\newcommand{\PP}{\ensuremath{\mathbb{P}}}

\newcommand{\N}{\ensuremath{\mathbb{N}}}
\newcommand{\R}{\ensuremath{\mathbb{R}}}

\theoremstyle{definition}

\theoremstyle{definition}

\theoremstyle{definition}

\theoremstyle{definition}

\newtheorem{thm}{Theorem}[section]
\newtheorem{cor}[thm]{Corollary}
\newtheorem{lemma}[thm]{Lemma}

\newtheorem{prop}[thm]{Proposition}
\newtheorem{assumption}[thm]{Assumption}
\theoremstyle{definition}
\newtheorem{defn}[thm]{Definition}
\newtheorem{rem}[thm]{Remark}
\newtheorem{eg}[thm]{Example}

\newtheorem{observe}[thm]{Observation}

\newcommand{\mc}{\mathcal}

\newcommand{\bt}{\mathbf{t}}

\newcommand{\eqdist}{\stackrel{(d)}{=}}

\def\beq{\begin{equation}}
\def\eeq{\end{equation}}
\def\ba{\begin{enumerate}[(a)]}
\def\bei{\begin{enumerate}[(i)]}
\def\be{\begin{enumerate}[(1)]}
\def\ee{\end{enumerate}}
\def\bi{\begin{itemize}}
\def\ei{\end{itemize}}
\def\beg{\begin{eg}}
\def\eeg{\end{eg}}
\def\bd{\begin{defn}}
\def\ed{\end{defn}}
\def\bt{\begin{thm}}
\def\et{\end{thm}}
\def\bl{\begin{lemma}}
\def\el{\end{lemma}}

\def\bc{\begin{cor}}
\def\ec{\end{cor}}
\def\bp{\begin{prop}}
\def\ep{\end{prop}}
\def\bo{\begin{observe}}
\def\eo{\end{observe}}
\def\bas{\begin{assumption}}
\def\eas{\end{assumption}}
\def\RR{\mathbb{R}}

\def\EE{\mathbb{E}}
\def\ZZ{\mathbb{Z}}
\def\NN{\mathbb{N}}

\def\PP{\mathbb{P}}

\newcommand{\boundgood}{\mathsf{G}}

\newcommand{\down}{\mathsf{Fall}}
\newcommand{\up}{\mathsf{Rise}}

\endlocaldefs

\begin{document}
\begin{frontmatter}
\title{KPZ equation correlations in time}
\runtitle{KPZ correlations in time}

\begin{aug}
\author[A]{\fnms{Ivan} \snm{Corwin}\ead[label=e1]{ic2354@columbia.edu}},
\author[B]{\fnms{Promit} \snm{Ghosal}\ead[label=e2]{promit@mit.edu}}
\and
\author[C]{\fnms{Alan} \snm{Hammond}\ead[label=e3]{alanmh@berkeley.edu}}

\address[A]{Department of Mathematics, Columbia University,
 \printead{e1}}

\address[B]{Department of Mathematics, Massachusetts Institute of Technology, \printead{e2}}
  
\address[C]{Departments of Mathematics and Statistics, U.C. Berkeley,
 \printead{e3}}  
\end{aug}

\begin{abstract}
We consider the narrow wedge solution to the Kardar-Parisi-Zhang stochastic PDE under the characteristic $3:2:1$ scaling of time, space and fluctuations. We study the correlation of fluctuations at two different times. We show that when the times are close to each other, the correlation approaches one at a power-law rate with exponent $2/3$, while when the two times are remote from each other, the correlation tends to zero at a power-law rate with exponent $-1/3$. We also prove exponential-type tail bounds for differences of the solution at two space-time points.

Three main tools are pivotal to proving these results: 1) a representation for the two-time distribution in terms of two independent narrow wedge solutions; 2) the Brownian Gibbs property of the KPZ line ensemble; and 3) recently proved one-point tail bounds on the narrow wedge solution.
\end{abstract}



\end{frontmatter}


%
%
%
%






\section{Introduction}

The Kardar-Parisi-Zhang stochastic PDE,
 written formally as
\begin{align}\label{eq:KPZ}
\partial_t \mathcal{H}(t, x) = \frac{1}{2}\partial^2_x \mathcal{H}(t,x) + \frac{1}{2} (\partial_x \mathcal{H}(t,x))^2 + \xi(t,x) \, ,
\end{align}
governs the random evolution in time $t \geq 0$ of an interface $\mathcal{H}(t,\cdot): \R \to \R$. The interface grows in the normal direction, while being smoothed by surface tension and roughened by a field of randomness that is independent between distinct space-time points -- a field represented in the equation by  the space-time Gaussian white noise $\xi$.
In the physically relevant {\it Cole-Hopf} solution, it is specified that $\mathcal{H}(t,x) = \log\mathcal{Z}(t,x)$, where $\mathcal{Z}$ is the solution to the stochastic heat equation (SHE)
\begin{align}\label{eq:SHE}
\partial_t \mathcal{Z}(t,x) &= \frac{1}{2}\partial^2_x \mathcal{Z}(t,x) + \xi(t,x) \mathcal{Z}(t,x) \, .
\end{align}
The fundamental solution to the SHE --  to be denoted here by $\mathcal{Z}^{\mathbf{nw}}(t,x)$ -- has initial data $\mathcal{Z}(0,x)$ set equal to a Dirac delta function $\delta(x)$ at $x=0$.
In this article, we study the resulting {\em narrow wedge} solution to the KPZ equation, namely
$\mathcal{H}^{\mathbf{nw}}(t,x)=\log \mathcal{Z}^{\mathbf{nw}}(t,x)$.
The initial condition $\mathcal{H}^{\mathbf{nw}}(0,x)$ is not well-defined; however, for small $t$, $\mathcal{H}^{\mathbf{nw}}(t,x)$ resembles the curve $-\tfrac{x^2}{2t}$,
so that the name `narrow wedge' may be seen as a substitute for the informative but prosaic term `narrow parabola'.

 The solution theory for the SHE is standard and is discussed in~\cite{Walsh,Corwin12}. The mathematical analysis of the KPZ equation has, however, offered many outstanding challenges: notably, for showing that the equation accurately approximates the evolution of physical growth processes; for developing a solution theory that is robust enough to permit such approximations; and for studying properties of the solution's probability distribution and various asymptotics.  Indeed, the Cole-Hopf KPZ solution has been shown in~\cite{Bertini1997} to coincide with limits of certain discrete growth processes; while the development of solution theory has been the object of intense recent activity, including the theory of regularity structures \cite{Hairer13}; energy solutions \cite{GJ2014a, GP2015a}; paracontrolled distributions \cite{GIP15,GP17}; and the renormalization group \cite{Kupiainen16}. The reader may consult~\cite{CorwinShen19} for a discussion of recent advances in the theory of singular stochastic PDE. Our present investigation focusses on the third area of challenge -- the properties of the KPZ equation's probability distribution.

Since the work of Kardar, Parisi and Zhang \cite{KPZ86}, the KPZ equation has been predicted to behave non-trivially in its long-time limit under a $3:2:1$ scaling of time, space and fluctuations. The $1/3$ fluctuation exponent was first verified in \cite{BQS11} for stationary initial data; in \cite{Amir11}, this exponent, and the corresponding limiting one-point distribution, was identified for narrow wedge initial data. The $2/3$ spatial exponent was first verified in \cite{CorHam16} for the narrow wedge.

The purpose of the present paper is to verify for the first time the KPZ equation's $3/3$ temporal exponent in the case of narrow wedge initial data. We will prove that the correlation of centred and scaled fluctuations at a pair of distinct moments in the time scale $t$ transitions between one as the times approach each other, and zero as they separate -- and we will quantify this transition with precise power-law exponents for the speed.

First we specify notation that represents a $3:2:1$ scaled version of $\mathcal{H}^{\mathbf{nw}}(t,x)$ in a manner suitable for the presentation and proof of our main results. The parameter $t > 0$ specifies a time scale; the value of  $\alpha > 0$ specifies  time judged on this scale; and $x \in \R$ specifies location judged on the suitable spatial scale $t^{2/3}$.
Indeed, for $t,\alpha>0$ and $x\in \RR$, we set
\begin{align}\label{eq:Upsilon}
\mathfrak{h}_{t}(\alpha,x) := t^{-1/3} \Big( \mathcal{H}^{\mathbf{nw}}\big(\alpha t, t^{2/3}x\big)+ \alpha t/24 \Big)\, .
\end{align}
The use of a pair $(\alpha,t)$ of temporal parameters may seem to introduce a touch of redundancy, but various scalings are in our view simplified by the use of these parameters; and we hope that conceptual clarity is offered by the division of roles between $t$, whose value  fixes a time scale, and~$\alpha$, which varies on the given scale.

\subsection{Main results}
Our first main result gives bounds on the correlation between $\mathfrak{h}_{t}(1,0)$ and $\mathfrak{h}_{t}(\alpha,0)$ for $\alpha>2$. In this case of two remote times, correlations decay as $\alpha^{-1/3}$ in the limit of high~$\alpha$. Recall that the correlation  between two random variables $X$ and $Y$ is defined to be
$$\mathrm{Corr}\big(X,Y\big)= \frac{\mathrm{Cov}\big(X,Y\big)}{\sqrt{\mathrm{Var}(X)}\sqrt{\mathrm{Var}(Y)}},\qquad \textrm{where}\qquad \mathrm{Cov}\big(X,Y\big) = \mathbb{E}[XY]-\mathbb{E}[X]\,\mathbb{E}[Y].$$

\bt[Remote correlations]\label{thm:Main1} There exist $t_0, c_1, c_2 >0$ such that, for all $t>t_0$ and $\alpha>2$,
\begin{align}\label{eq:Corr1}
c_1 \alpha^{-1/3} \leq \mathrm{Corr}\big(\mathfrak{h}_{t}(1,0), \mathfrak{h}_{t}(\alpha , 0)\big) \leq c_2 \alpha^{-1/3} \, .
\end{align}
\et
\smallskip

The second main result  bounds the correlation between $\mathfrak{h}_{t}(1,0)$ and $\mathfrak{h}_{t}((\beta+1) , 0)$ when $\beta$ is smaller than $1/2$. As $\beta$ approaches zero, the correlation between $\mathfrak{h}_{t}(1,0)$ and $\mathfrak{h}_{t}(\beta+1, 0)$ approaches one with a discrepancy of order $\beta^{2/3}$.

\bt[Adjacent correlations]\label{thm:Main2}
For any $t_0>0$, there exist $c_1=c_1(t_0)>0$ and $c_2=c_2(t_0)>0$ such that, for all $t>t_0$ and $\beta\in (0,1/2)$ satisfying $\beta t>t_0$,
\begin{align}\label{eq:Corr2}
c_1 \beta^{2/3}\leq 1-\mathrm{Corr}\big( \mathfrak{h}_{t}(1,0), \mathfrak{h}_{t}(1+\beta, 0)\big)\leq c_2\beta^{2/3} \, .
\end{align}
\et
\smallskip

These theorems invite a number of  interesting questions. Notice that, in Theorem \ref{thm:Main1}, the minimal time $t_0$ must be large enough for the bounds to hold. First, we may ask: what happens to the two bounds for small $t$? In fact, in our proof of this theorem, we show a stronger upper bound: for any $t_0>0$, there exists $c_2(t_0)>0$ so that, for all $t>t_0$ and $\alpha>2$, $\mathrm{Corr}\Big(\mathfrak{h}_{t}(1,0), \mathfrak{h}_{t}(\alpha , 0)\Big)$ is at most $c_2 \alpha^{-1/3}$. Our proof does not, however, produce a matching lower bound. It is unclear to us if a new phenomenon occurs at short time that renders such a bound invalid or if instead our technique of proof is unsuitable for obtaining a bound of this form.

A second question also touches on short-time behaviour. In Theorem \ref{thm:Main2} we require $\beta t >t_0$ for some arbitrary yet fixed $t_0$. For $t$ fixed, this limits us to values of $\beta>t_0/t$. (As $t\nearrow \infty$, this lower bound tends to zero.) What happens for $\beta$ smaller than $t_0/t$? The source of the restriction $\beta t >t_0$ comes from the application of the one-point tail behaviour from Propositions \ref{ppn:onepointlowertail} and \ref{ppn:onepointuppertail}. These tail bounds are valid provided that time exceeds $t_0$. When the time gap $\beta t$ is less than $t_0$, we need short-time tail bounds which are presently not available in the literature. As we remark at the end of Section \ref{sec:broader}, this condition that $\beta>t_0/t$ (as well as both of our main theorems) has an interpretation in terms of \emph{aging}. 

 Finally, a third question concerns the specification of the constants $c_1$ and $c_2$ in Theorem \ref{thm:Main1} and \ref{thm:Main2}. In the physics literature, de Nardis and Le Doussal \cite{de_Nardis_Le_Doussal_2017,de_Nardis_2018} provide various predictions concerning the two time distribution, in particular the constants. For instance, they claim that as $\alpha,t\to \infty$ the correlation in  Theorem \ref{thm:Main1} behaves like $\alpha^{-1/3}$ times a limiting constant written in terms of the  Airy$_2$ process $\mathcal{A}_2(\cdot)$ and a two-sided Brownian motion $B(\cdot)$ as
\begin{align}
\lim_{\alpha\to \infty}\lim_{t\to \infty} \alpha^{\frac{1}{3}}\mathrm{Corr}(\mathfrak{h}_t(\alpha,0), \mathfrak{h}_{t}(1,0)) = \frac{\mathrm{Corr}\big(\mathcal{A}_2(0),\max_{y}\{\mathcal{A}_2(y)-y^2+\sqrt{2}B(y)\}\big)}{\mathrm{Var}\big(\mathcal{A}_2(0)\big)}
\end{align}
They also predict that the conditional density of $\mathfrak{h}_{t}(1+\beta,0)- \mathfrak{h}_{t}(1,0)$ given $\mathfrak{h}_{t}(1,0)$ converge to the Baik-Rains distribution as $t$ goes to $\infty$ followed by $\beta $ tends to $0$. From this prediction one sees that in this limit, the correlation in Theorem \ref{thm:Main2} should behave like $\beta^{2/3}$ times a constant given by the ratio of the variances of Baik-Rains distribution and the Tracy-Widom distribution. We have not tried to optimize our methods to access these precise values of the constants, though it is certainly worth more study.

We further mention some other natural themes related to our first two main results. Theorems~\ref{thm:Main1} and~\ref{thm:Main2} both probe the two-time correlation when the spatial coordinate is zero. In this context, two space-time points lie along a space-time line of velocity zero emanating from the origin. In general, space-time lines of fixed velocity emanating from the origin are called \emph{characteristics}; and, by trivial affine shifts -- see Proposition \ref{ppn:Stationarity} -- our results imply the same correlation behaviour along such characteristics. On the other hand, we have not addressed what happens when, say,  $\mathrm{Corr}\big(\mathfrak{h}_{t}(1,0), \mathfrak{h}_{t}(\alpha , x)\big)$ is considered for $x\neq 0$. This is a problem that may plausibly be addressed by the methods of the present article. That said, if the aim is to prove results concerning more general initial data than the narrow wedge for the KPZ equation, our methods will not be useful without significant further input. Indeed, as the proof sketch offered in Section \ref{sec:idea} will explain, we rely upon a special property of narrow wedge initial data, namely its connection to the KPZ line ensemble and this ensemble's Gibbs property; and this essential aspect is lacking when general initial data is considered.

We state two further theorems  concerning local fluctuations of the KPZ equation in {\em space} and {\em time}.

\bt\label{t.locreg}
For  any $t_0 >0$ there exist $s_0 = s_0(t_0)>0$ and $c=c(t_0)>0$  such that for all $x \in \R$, $t\geq t_0$, $s\geq s_0$, and $\e \in (0,1]$,
\begin{align}\label{eq:SuphBd}
 \PP \, \bigg( \, \sup_{z \in [x,x+\e]} \Big\vert \mathfrak{h}_{t}(1,z)+ \tfrac{z^2}{2} - \mathfrak{h}_{t}(1,x)- \tfrac{x^2}{2} \Big\vert \, \geq \,  s\e^{1/2}  \, \bigg) \leq   \exp \big( - c s^{3/2} \big) \, .
 \end{align}
\et

This theorem gauges fluctuation of the parabolically shifted process  $\mathfrak{h}_{t}(1,x)+\tfrac{x^2}{2}$ since, as we will recall in  Proposition \ref{ppn:Stationarity}, this shift renders the process stationary in the $x$ variable.
Before turning to our theorem concerning fluctuations in time, we state a corollary of the preceding result, which gives an estimate on the spatial modulus of continuity for $\mathfrak{h}_{t}$. The corollary holds uniformly in the limit of high $t$, so that, when the result is allied with the one-point convergence that will be the subject of Proposition~\ref{ppn:Stationarity}, we learn that the spatial process is tight and that any subsequential limit shares the  H\"{o}lder continuity known for the finite $t$ processes.

\bc\label{thm:SpaceHolder}
For $t_0>0$ and any interval $[a,b]\subset \RR$, define
\begin{align}\label{eq:sholder}
\mathcal{C} :=\sup_{x_1\neq x_2\in [a,b]} |x_1-x_2|^{-\frac{1}{2}}\Big( \log\frac{|b-a|}{|x_1-x_2|}\Big)^{-2/3}\Big|\mathfrak{h}_{t}(1,x_1)+\tfrac{x_1^2}{2} - \mathfrak{h}_{t}(1,x_2) - \tfrac{x_2^2}{2} \Big| \, .
\end{align}
Then there exist $s_0=s_0(t_0,|b-a|)>0$ and $c=c(t_0,|b-a|)>0$  such that, for $s\geq s_0$ and $t>t_0$,
$$
\mathbb{P}\big(\mathcal{C}>s \big)\leq \exp\big(-cs^{3/2}\big) \, .
$$
\ec

Our last result concerns the upper and lower tails for the difference in fluctuations at two times. These tail bounds will be used in the proofs of Theorems \ref{thm:Main1} and \ref{thm:Main2}.

\bt\label{ppn:HDivTail}
For any $t_0>0$, there exist $s_0=s_0(t_0)>0$, $c_1=c_1(t_0)>0$, $c_2=c_2(t_0)>0$ and $c_3= c_3(t_0)>0$  such that, for  $s>s_0$, $t>t_0$ and $\beta \in (0,1/2)$ for which $t \beta>t_0$,
 \begin{align}
\exp\big(-c_1 s^{3/2}\big) \leq \mathbb{P}\Big(\mathfrak{h}_{t}(1+\beta,0)-\mathfrak{h}_{t}(1,0) \geq \beta^{1/3}s\Big)&\leq \exp\big(-c_2 s^{3/4}\big) \, ,\label{eq:TailBd1}\\
 \mathbb{P}\Big(\mathfrak{h}_{t}(1+\beta,0)-\mathfrak{h}_{t}(1,0)\leq -\beta^{1/3}s\Big)&\leq \exp\big(-c_3s^{3/2}\big) \, .  \label{eq:TailBd2}
 \end{align}
 \et

Just as in Theorem \ref{thm:Main2}, we suspect different behaviour will arise in the short-time limit as $\beta\to 0$ for $t$ fixed. Because of this lack of a short-time result, we do not present a temporal modulus of continuity result arising from this theorem in the style of the spatial Corollary \ref{thm:SpaceHolder}.  Notice, also, that we do not record a lower bound on the probability in~\eqref{eq:TailBd2}; we have not pursued this since we presently have no application for it. Finally, although we focus here on $x=0$ at both times, it should also be possible to address two different spatial locations.

\subsection{Idea of proof}\label{sec:idea}
Before describing our methods, we explain why existing approaches to studying the KPZ equation do not yield our main results. The  KPZ equation has attracted attention of specialists in stochastic PDE -- via regularity structures, paracontrolled distributions or the renormalization group, for example -- and of probabilists using integrable methods involving, for instance, Macdonald processes and the Bethe ansatz. It is natural enough to ask why our results do not follow from techniques in either of these areas. Stochastic PDE methods are well suited to the study of  local problems but they have little to say about the distribution of KPZ equation  under the characteristic $3:2:1$ scaling. For example, that the one-point distribution has fluctuations of order~$t^{1/3}$ is inaccessible via these techniques. Integrable probability has been able to identity both the $t^{1/3}$ scaling as well as the limiting one-point fluctuations. However, as of yet, there has been no rigorous progress in deriving  explicit formulas for the joint distribution of the KPZ equation at several space-time points. Even if such formulas existed, it might not be easy to extract our results from them. For instance, for zero temperature models (as discussed in Section \ref{sec:broader}), explicit multi-point formulas exist, but they have yet to prove valuable for  extracting correlation decay results in the style of Theorems~\ref{thm:Main1} and~\ref{thm:Main2}.

How do we proceed? Our study is principally probabilistic, and a vital aspect is the use of Gibbsian line ensembles -- objects that are born integrably but whose lives are in large part lived probabilistically. Indeed, the analysis of a Gibbsian line ensemble
is one of three pivotal tools concerning the KPZ equation that will enable our proofs -- tools whose use is here briefly described, and in more detail in Section \ref{sec:Tools}.

The first  tool, the {\em composition law} in Proposition \ref{ppn:TimeEvl}, realizes the two-time distribution in terms of an exponentiated convolution of two independent narrow wedge KPZ equation solutions.

The second  tool is the pertinent Gibbsian line ensemble. This is the {\em KPZ line ensemble}, which is recalled in Proposition \ref{NWtoLineEnsemble}. The narrow wedge solution to the KPZ equation for any fixed time is embedded as the lowest indexed curve of an infinite ensemble of curves which jointly enjoy a {\em Brownian Gibbs property}. This property says that, fixing an index and an interval, the law of that indexed curve on that interval only depends on the boundary data (the curve indexed one above and one below on that interval, and the starting and ending points) and is comparable to the law of Brownian bridge with this endpoint pair, the comparison made succinctly via a Radon-Nikodym reweighting depending on said boundary data. Moreover, the Brownian Gibbs property implies {\em stochastic monotonicity} (Lemma \ref{Coupling1}) which shows that shifting the boundary data in a given direction likewise shifts the law of the curve conditioned on that data.

The third tool consists of {\em tail bounds} (Propositions \ref{ppn:onepointlowertail} and \ref{ppn:onepointuppertail}) for the distribution of $\mathfrak{h}_t(1,0)$.

How do these tools combine to produce our main results? To compare $\mathfrak{h}_t(1,0)$ with $\mathfrak{h}_t(\alpha,0)$, we first use the composition law to realize $\mathfrak{h}_t(\alpha,0)$ in terms of the process $\mathfrak{h}_t(1,\cdot)$ and an independent and scaled (based on the value of $\alpha$) narrow wedge KPZ equation solution which we denote later by $\mathfrak{h}_{\alpha t \downarrow t}(\cdot)$. The composition law is a {\em softening} of a variational problem  in which one would instead maximize the sum of the two narrow wedge solutions over their spatial argument. Our composition law only matches this variational problem in the $t\nearrow\infty$ limit, but the limiting case offers a convenient venue for a brief description of our approach.


Indeed, for the limiting maximization problem -- with $t$ formally infinite --  controlling the difference or correlation between  $\mathfrak{h}_t(1,0)$ with $\mathfrak{h}_t(\alpha,0)$ boils down to three main steps. First we must show that the maximization likely occurs for values of $x$ near zero. Second we must show that the value of the functions in consideration at the maximizing location are close to their values at zero. And finally, we must show that one of the two random variables $\mathfrak{h}_t(1,0)$ or $\mathfrak{h}_{\alpha t \downarrow t}(0)$ is small compared to the other one. Which process is small depends on whether we are working with $\alpha$ large -- the case of remote correlations -- in which case $\mathfrak{h}_t(1,0)$ should be correspondingly small; or if we are considering $\alpha=1+\beta$ for $\beta$ small -- the case of adjacent correlations -- in which case  it is $\mathfrak{h}_{\alpha t \downarrow t}(0)$ that should be small.

These steps can be realized by using the Gibbsian line ensemble tool along with the tail bounds. Using the tail bounds we may control (with exponentially small tail probability) the boundary data for the lowest labeled curve of the KPZ line ensemble -- i.e., the narrow wedge solution at fixed time. Allying this with stochastic monotonicity enables us to transform our problems into questions involving the fluctuations of Brownian bridges, which can be treated quite classically. This general argument motif  of using the Gibbs property to transfer one-point information into spatial regularity originated in \cite{CH14} where the Airy line ensemble was introduced and studied, and has been developed in various directions in many subsequent works on Gibbsian line ensembles such as \cite{CorHam16,Corwin2018,Hammond18, Hammond18a, Hammond18b, Hammond18c, HammondSarkar, HammondSarkarBaus, Basicproperties, Landscape, caputo2019, CIW2019}. An example in which we closely follow a proof in the literature is Proposition \ref{p.locreg}, which mimics \cite[Proposition~$3.5$]{Hammond18b}. In most other cases, while we follow the general motif, we are forced to develop new variations since our present work is the first instance of applying Gibbsian line ensemble methods to study the temporal regularity of a KPZ class model.

Our proofs must operate when $t$ is finite, rather than in the limiting case of $t \nearrow \infty$. We have mentioned that   the finite-$t$ KPZ equation composition law involves not a variational problem, but its softened convolution formula. We thus need to vary the proposed approach, arguing that the principal contribution to the integrals in the concerned composition law occurs near the integrand's maximizer.
A second complication arises from our seeking in  Theorems \ref{thm:Main1} and \ref{thm:Main2} to control correlations, while the  tools we have indicated merely control tail bounds on fluctuations. The needed transition is achieved in the proofs of these theorems via some rather general arguments contained in Appendix \ref{sec:Appendix}.

While the variational problem version of the composition law is a helpful, albeit only heuristic, way of thinking in the context of the KPZ equation, it is precisely the composition law for the {\em KPZ fixed point} where the narrow wedge solutions are replaced by suitably scaled Airy processes. In fact, we first developed our arguments in this simpler context. However, several recent works described next in Section \ref{sec:broader} probe temporal correlation for the KPZ fixed point and other {\em zero temperature} models such as exponential, geometric and Brownian last passage percolation. While the methods used in these investigations do not seem amenable (at least with present tools) to application to {\em positive temperature} models like the KPZ equation, the Gibbsian line ensemble technique is readily lifted to this level, so that we may employ it here. In Section \ref{sec:KPZfixedpoint}, we describe the KPZ fixed point analogue of our work in greater detail.

Our Gibbsian line ensemble approach to studying temporal correlations has further potential. On the zero temperature side, it should be applicable to all of the just mentioned last passage percolation models. On the positive temperature side, there is a growing body of models which can be embedded into a Gibbsian line ensemble including the semi-discrete polymer  \cite{ON12}; log-gamma / strict-weak polymers  \cite{COSZ}; and stochastic six vertex model / ASEP \cite{BBW16,Corwin2018}. In each of these positive temperature cases (besides the KPZ equation), there are some challenges to implementing the methods from the present work. For example, suitably strong tail bounds have yet to be demonstrated for the polymer models; and  the composition law is considerably more complicated for the stochastic six vertex model and ASEP.
That these models embed into Gibbsian line ensemble is a facet of their underlying {\em integrability}. Indeed, the study of Gibbsian line ensembles and their marriage of integrable and non-integrable probabilistic ideas have been quite fruitful recently.

\subsection{Broader context}\label{sec:broader}


The main results, Theorems \ref{thm:Main1} and \ref{thm:Main2}, fit into a broader effort in the last few years to understand the temporal correlation structure for the KPZ equation and other models in its eponymous universality class. The experimental observations of \cite{Takeuchi2012} about temporal correlations brought this question into focus. Soon after, due to the applicability of the replica trick, there were a few informative non-rigorous investigations in the physics community of the two-time distribution for the KPZ equation. Using certain combinatorial approximations, \cite{Dotsenko_2013,Dotsenko_2015,Dotsenko_2016} derived a formula for the two-time distribution. Later work of \cite{de_Nardis_Le_Doussal_2017}  argued against the validity of this approximation and the resulting formula, and derived a formula for the two-time distribution when one of the fluctuations is taken deep in its tail (which simplifies some combinatorics in the Bethe ansatz used there). That work also provided an argument based on this approximate formula for the type of decay of correlation bounds proved here (along with predictions for the limiting values of the constants in those bounds). Further details are contained in \cite{de_Nardis_2018} and an independent calculation leading to the same conclusion is provided in \cite{D17}. In~\cite{De_Nardis_Le_Doussal_Takeuchi}, time correlations are further investigated using numerical simulations.

Ours are the first rigorous results regarding the temporal correlations of the KPZ equation or any other \emph{positive temperature} models. As we will review momentarily, there has been considerable recent mathematical activity in the analysis of temporal correlations of \emph{zero-temperature} models in the KPZ universality class which enjoy \emph{determinantal structure}. While these results have no direct implications for the KPZ equation, they do inform our expectation for its behaviour.
On the other hand, with the exception of slow decorrelation, the methods used in the works that we now mention do not seem to generalize to positive temperature.

\emph{Slow decorrelation}, shown to hold true for many KPZ class models in \cite{Ferrari_2008} and  \cite{corwin2012}, implies that, for any $\eta<1$ and any pair of times $t_1= t$ and $t_2 = t+t^{\eta}$, as $t\nearrow \infty$ the fluctuations (up to centring and scaling) converge to the same limiting random variable. This choice corresponds to \emph{very adjacent} times in out setup -- i.e., taking $\beta\to 0$ as an inverse power of $t$ in Theorem \ref{thm:Main2}.

Studying exponential last passage percolation as well as limiting Airy process variational formulas for the two-time distribution, \cite{FS16} presented two non-rigorous approaches to studying remote and adjacent correlation decay. Such results are proved in~\cite{FO18} and~\cite{BG18} for exponential and Poissonian last passage percolation. Neither of the latter works rely on explicit formulas for the two-time distribution, but rather on more probabilistic characterizations. The work of \cite{FO18} is also able to address more general initial data than just narrow wedge.  The exponents for the adjacent correlations under flat or stationary initial data are predicted to be $2/3$. For flat data, \cite{FS16} predicted that the remote correlation exponent in exponential  last passage percolation is $-1$, instead of $-1/3$;   recently, \cite{BGZ19} proved this prediction. For the KPZ equation with flat  data, we expect to see the same exponent as in \cite{BGZ19}. It is worth investigating whether our Gibbsian line ensemble methods are useful in addressing correlations for more general initial data for the KPZ equation than narrow wedge.

There are, in fact, some proven explicit formulas for the two-time distribution of certain zero-temperature determinantal KPZ class models. The first formula was derived in \cite{Johansson2017} and concerns Brownian last passage percolation -- and, through a limit transition, the KPZ fixed point. The formula is complicated enough that extracting remote and adjacent correlation seems arduous. Recently, \cite{J18} has derived a new and much more manageable formula for geometric last passage percolation which has permitted the study, in \cite{J19}, of limits of the two-time distribution for adjacent and remote times. So far, explicit formulas have not been used to prove correlation results. In \cite{BZ17}, multi-time formulas for periodic last passage percolation (or TASEP) have been derived in the time scale on which the system relaxes to equilibrium.  However, these formulas are again rather complicated and it is unclear how tractable they may render the zero-temperature counterparts to the questions that we consider.  Quite recently, \cite{J19B}  and \cite{Liu_2019} proved multi-time formulas and asymptotics in the non-periodic case.  As of yet, there are no multi-time exact formulas known for the KPZ equation or any other positive temperature KPZ class model.

Excepting \cite{FO18,BGZ19}, which  work more generally, the articles that we have mentioned treat specific types of initial data. Recently, there have been significant advances~\cite{KPZfixed},~\cite{Hammond18},~\cite{Landscape} regarding the KPZ fixed point with general initial data. In fact, all these works contain H\"older continuity results for the fixed point; see also \cite{P17} and \cite{HammondSarkar} for related results. Our spatial and temporal fluctuation results (and consequential modulus of continuity estimates) agree in the limit of high~$t$ with the H\"older continuity of the KPZ fixed point.

We close this discussion by recalling our initial motivation to study this problem. In 2010, Amir Dembo and Jean-Dominique Deuschel asked one of us whether the KPZ equation \emph{ages}. 
 As explained in \cite{Dembo2007} (a work providing many helpful references on the subject), aging is a phenomenon in glassy materials in which ``older systems relax in a slower manner than younger ones''. (We refer to \cite[Section 2.3]{BenArous} for a detailed discussion on aging.) The nature of this relaxation can be studied via the two-time correlation; aging corresponds to correlation crossing over from one to zero as the times move from being adjacent to remote. In this sense, our main theorems prove that the KPZ equation does, indeed, age. In fact,  a referee has pointed out to us that the assumption that $t>t_0$ in Theorems \ref{thm:Main1} and \ref{thm:Main2} is natural from the perspective of aging. For aging to occur, typically a system has to evolve (or relax) for some time from its initial configuration. This analogy with aging requires further investigation and could be complemented by the detailed study of the correlation exponents for KPZ started from its stationary initial data.

 The authors may attest that the years elapsed since they learnt of this question have furnished them with vivid indications of the power of the phenomenon of aging -- and not merely through the study of the KPZ equation.

\subsection{Outline}
Section \ref{sec:Tools} reviews several important and known properties of the KPZ equation, including its composition law, its relation to the KPZ line ensemble, and its one-point tail bounds. Since our analysis of the two-time distribution involves a cousin of a classical variational problem at zero temperature, and since zero temperature counterparts to our study have recently been undertaken, it is profitable -- though not necessary for understanding our proofs -- to view our results through the lens of zero temperature; and in
Section \ref{sec:KPZfixedpoint} we do so.
Our main technical contribution begins in Section \ref{sec:spatialtails} where we demonstrate how to extend one-point tail bounds to bounds on spatial fluctuation tails.
While Propositions~\ref{ppn:LowerTail} and \ref{ppn:UpperTail} contain global spatial fluctuation results which measure the size of spatial fluctuations on all of space, Propositions \ref{p.locreg}  and \ref{cor:HitParabola} contain local fluctuation results. This section also contains the proof of the local spatial fluctuation Theorem~\ref{t.locreg} which follows quite readily by combining Propositions \ref{ppn:LowerTail}, \ref{ppn:UpperTail} and \ref{p.locreg}.
The remote correlation decay  Theorem \ref{thm:Main1}, adjacent correlation Theorem \ref{thm:Main2} and spatial modulus of continuity Corollary~\ref{thm:SpaceHolder}
are successively proved in Sections~\ref{sec:remote},~\ref{sec:adjacent} and~\ref{sec:Holderspace}.
Section \ref{sec:Appendix}, the appendix, contains several general probabilistic results which we develop to relate tail probabilities (to which we generally have access) to covariances and correlations.

\subsubsection*{Acknowledgments}
Aspects of this project were prompted by the question that we have mentioned, asked by Amir Dembo and Jean-Dominique Deuschel during a conference at MSRI in fall 2010.
In an earlier version of this work, Xuan Wu and Yier Lin pointed out a missing detail in the proof of  Proposition~\ref{p.locreg} regarding controlling the probability of the $\up$ event therein. The authors thank Xuan and Yier for their close reading and comments. The authors also wish to acknowledge useful comments from Patrik Ferrari, Pierre Le Doussal and Kazumasa Takeuchi.
Ivan Corwin was partially supported by the Packard Fellowship for Science and Engineering, and by the NSF through grants DMS-$1811143$ and DMS-$1664650$. Alan Hammond was partially supported by the NSF through grants DMS-$1512908$ and DMS-$1855550$ and by a Miller Professorship at U.C. Berkeley.

\section{Tools}\label{sec:Tools}

In this section, we recall significant results that we will need. In Section \ref{sec:KPZLE}, the KPZ line ensemble is introduced (in Proposition \ref{NWtoLineEnsemble}). This ensemble  enjoys two key properties: (1) its lowest indexed curve coincides in law with the fixed time narrow wedge Cole-Hopf solution to the KPZ equation; and (2) it enjoys a certain Brownian Gibbs property. The general notion of a line ensemble and the Brownian Gibbs property are recorded in Definition~\ref{LineEns}, while Lemma~\ref{Coupling1} records certain monotonicity results associated with this Gibbs property. Section \ref{sec:KPZeqinputs} contains results about the solution to the KPZ equation such as its stationarity (Proposition \ref{ppn:Stationarity}); positive association (Propositions~\ref{FKGProp} and \ref{ppn:CondEIn});  composition law (Proposition \ref{ppn:TimeEvl}); and tail bounds (Propositions~\ref{ppn:onepointlowertail} and~\ref{ppn:onepointuppertail}).

Before commencing, let us record a few pieces of {\em notation} which will be used throughout. We write  $\N = \{1,2,\ldots\}$. We will often discuss probabilities without specifying the probability space. When we use symbols for an event in such a probability space, we will often use the styles $\mathsf{A}$ or $\mathcal{A}$ instead of the standard $A$. For two events $\mathsf{A}$ and $\mathsf{B}$ we will sometimes write $\mathbb{P}(\mathsf{A}, \mathsf{B})$ instead of $\mathbb{P}(\mathsf{A}\cap \mathsf{B})$. In the statements and proofs of many of our results, we will use the hopefully obvious notation $c=c(\cdot)>0$ to represent a positive constant $c$ that depends on  the variable in place of $\cdot$.  In some of our proofs, we will allow constants such as this to vary line to line and within lines to simplify the exposition and avoid introducing too many constants. Finally, we will sometimes use the shorthand
$$\mathfrak{h}_{t}(x):=\mathfrak{h}_{t}(1,x)$$
when the time parameter is fixed and our interest is in the spatial process.
\subsection{KPZ line ensemble}\label{sec:KPZLE}

It was shown in \cite{CorHam16} that the narrow wedge solution of the KPZ equation may be  embedded as the lowest indexed curve of an infinite ensemble of curves which enjoys a \emph{Brownian Gibbs} property. Describing this requires a little notation. The term \emph{line} in `line ensemble' alludes to the piecewise constant curves in counterpart ensembles associated to models such as Poissonian last passage percolation. The term is a misnomer in the present context because the ensembles in question have a locally Brownian, and hence continuous, structure.

\bd[Brownian Gibbs line ensembles; Definitions 2.1 and 2.2 of \cite{CorHam16}]\label{LineEns}
Let $\Sigma\subset \ZZ$ and $\Lambda\subset\RR$ be  intervals. Let $X$ be the set of continuous functions $f:\Sigma \times \Lambda\to \RR$, equipped with the topology of uniform convergence on compact subsets;  and let $\mathcal{C}$ be the $\sigma$-field generated by $X$.

A $(\Sigma\times \Lambda)$-\emph{indexed line ensemble} $\mathcal{L}$ is a random variable $\mathcal{L}$ defined on a probability space $(\Omega,\mathcal{B},\mathbb{P})$ taking values in $X$ such that $\mathcal{L}$ is a $(\mathcal{B}, \mathcal{C})$ measurable function. In other words, $\mathcal{L}$ is a set of random continuous curves indexed by $\Sigma$ where each of those curves maps $\Lambda$ to $\RR$. An element of $\Sigma$ is a curve index, and we will write $\mathcal{L}^k(x)$ instead of $\mathcal{L}(k,x)$ for $k\in \Sigma$ and $x\in \Lambda$; we will write $\mathcal{L}^k$ for the entire index $k$ curve. In a  standard notational abuse regarding random variables, we have  suppressed the dependence on $\omega\in \Omega$. We will generally replace $\Lambda$ by the interval $(a,b)$ for $a<b$.  We will also often consider $\Sigma$ either to equal $\NN$ or $\{k_1,\ldots, k_2\}$ for some pair of integers $k_1<k_2$.

In order to formulate the Brownian Gibbs property, we need a background measure on line ensembles -- the free Brownian bridge line ensemble measure. For any two integers  $ k_1 <k_2$, two vectors of reals $\vec{x}, \vec{y} \in \RR^{k_1-k_2+1}$, and an interval $(a,b)$, we say that a $\{k_1, \ldots , k_2\}\times (a,b)$-indexed line ensemble is a \emph{Brownian bridge line ensemble with entrance data $\vec{x}$ and  exit data $\vec{y}$} if its law, which we denote by $\mathbb{P}^{k_1, k_2, (a,b), \vec{x}, \vec{y}}_{\mathrm{free}}$, is equal to that of $k_2- k_1+1$ independent Brownian bridges starting at values~$\vec{x}$ at $a$  and ending at values $\vec{y}$ at $b$. We use the notation $\mathbb{E}^{k_1, k_2, (a,b), \vec{x}, \vec{y}}_{\mathrm{free}}$ to denote the expectation with respect to the probability measure $\mathbb{P}^{k_1, k_2, (a,b), \vec{x}, \vec{y}}_{\mathrm{free}}$. When $k_1=k_2=1$, we write $\mathbb{P}^{(a,b), \vec{x}, \vec{y}}_{\mathrm{free}}$. It is natural to think of $a$ and $b$ as times and $\vec{x}$ and $\vec{y}$ as starting and ending locations for the Brownian bridges. However, these notions of time and space are {\em not} the same as in the KPZ equation, so we will avoid this usage.

Suppose given a continuous function  $\mathbf{H}:\RR\to [0,\infty)$ which we will call a \emph{Hamiltonian}. Our attention will be fixed almost exclusively on a choice of the form
\begin{align}
\mathbf{H}_t(x)= e^{t^{1/3}x} \, \, \, \, \textrm{for given $t > 0$} \, .\label{eq:Ht}
\end{align}

For two measurable functions $f,g:(a,b)\to \RR$,  a $\{k_1, \ldots , k_2\}\times (a,b)$-{\em indexed $\mathbf{H}$-Brownian bridge line ensemble with entrance data $\vec{x}$, exit data $\vec{y}$ and boundary data $(f,g)$} is a collection of random curves  $\mathcal{L}^{k_1}, \ldots , \mathcal{L}^{k_2}: (a,b)\to \RR$ whose law we denote by $\mathbb{P}^{k_1, k_2, (a,b), \vec{x}, \vec{y}, f,g}_{\mathbf{H}}$. This law is specified  by the Radon-Nikodym derivative
\begin{align}
\frac{d\mathbb{P}^{k_1, k_2,(a,b), \vec{x}, \vec{y}, f, g}_{\mathbf{H}}}{d\mathbb{P}^{k_1, k_2, (a,b), \vec{x}, \vec{y}}_{\mathrm{free}}}(\mathcal{L}^{k_1}, \ldots , \mathcal{L}^{k_2}) = \frac{W^{k_1, k_2, (a,b), \vec{x}, \vec{y}, f,g}_{\mathbf{H}}(\mathcal{L}^{k_1}, \ldots , \mathcal{L}^{k_2})}{Z^{k_1, k_2, (a,b), \vec{x}, \vec{y}, f,g}_{\mathbf{H}}} \, ,
\end{align}
where  $Z^{k_1, k_2, (a,b), \vec{x}, \vec{y}, f, g}_{\mathbf{H}}$ is the normalizing constant which produces a probability measure, and
\begin{align}\label{eq:W}
W^{k_1, k_2, (a, b), \vec{x}, \vec{y}, f,g}_{\mathbf{H}}(\mathcal{L}^{k_1}, \ldots , \mathcal{L}^{k_2})  = \exp\bigg\{- \sum_{i=k_1}^{k_2+1} \int \mathbf{H}\big(\mathcal{L}^{i}(x)- \mathcal{L}^{i-1}(x)\big) dx \bigg\} \, .
\end{align}
In the right-hand side of the preceding display, the boundary conditions are in force via the setting of $\mathcal{L}^{k_1-1}$ equal to $f$, or to   $+\infty$ if $k_1-1\notin \Sigma$; and of $\mathcal{L}^{k_2+1}$ equal to $g$, or to $-\infty$ if $k_2+1\notin \Sigma$.

We will say that a $(\Sigma\times\Lambda)$-indexed line ensemble $\mathcal{L}$ enjoys the \emph{$\mathbf{H}$-Brownian Gibbs property} if, for all $K = \{k_1,\ldots, k_2\}\subset \Sigma$ and $(a,b)\subset \Lambda$, the following distributional equality holds:
\begin{align}
\mathrm{Law}\Big(\mathcal{L}_{K\times (a,b)} \text{ conditioned on }\mathcal{L}_{\Sigma\times \Lambda \backslash K \times (a,b)}\Big) = \mathbb{P}^{k_1,k_2, \vec{x}, \vec{y}, f, g}_{\mathbf{H}} \, ,
\end{align}
where $\vec{x}= (\mathcal{L}^{k_1}(a), \ldots , \mathcal{L}^{k_2}(a))$, $\vec{y} =(\mathcal{L}^{k_1}(b),\ldots , \mathcal{L}^{k_2}(b))$, and where again $f = \mathcal{L}^{k_1-1}$ (or $+\infty$ if $k_1-1\notin \Sigma$) and $g = \mathcal{L}^{k_2+1}$ (or $-\infty$ if $k_2+1\notin \Sigma$). That is, the ensemble's restriction to $K\times (a,b)$ is influenced by the complementary information only via only boundary data (the starting and ending points and the neighbouring curves); and, given this pertinent data, the law of the restriction is a $\mathbf{H}$-Brownian Gibbs line ensemble with the boundary parameters induced by the data.

The $\mathbf{H}$-Brownian Gibbs property may be viewed as a spatial Markov property. Just as for Markov processes, it is useful to have a \emph{strong} version of the Gibbs property which is valid with respect to \emph{stopping domains}. This we now describe. For a line ensemble $\mathcal{L}$ as above, let $\mathfrak{F}_{\textrm{ext}}\big(K\times (a,b)\big)$ denote the $\sigma$-field generated by the curves $K\times (a,b)$. A pair $(\mathfrak{a},\mathfrak{b})$ of random variables is called a $K$-stopping domain if
$\big\{\mathfrak{a} \leq a, \mathfrak{b}\geq b\big\} \in \mathfrak{F}_{\textrm{ext}}\big(K\times (a,b)\big)$.
Let $C^{K}(a,b)$ denote the set of continuous $K$-indexed functions $(f^{k_1},\ldots, f^{k_2})$, each from $(a,b)$ to $\RR$; and let
$$
C^K := \Big\{ (a,b,f^{k_1},\ldots, f^{k_2}):a<b \textrm{ and } (f^{k_1},\ldots, f^{k_2})\in C^K(a,b)\Big\}\, .
$$
Write $\mathcal{B}(C^K)$ for the set of all Borel measurable functions from $C^K$ to $\RR$. A $K$-stopping domain $(\mathfrak{a},\mathfrak{b})$ satisfies the \emph{strong $\mathbf{H}$-Brownian Gibbs property} if, for all $F\in \mathcal{B}(C^K)$, $\mathbb{P}$-almost surely
$$
\mathbb{E} \Big[ F\big(\mathfrak{a},\mathfrak{b},\mathcal{L}\big\vert_{K\times (\mathfrak{a},\mathfrak{b})}\big) \Big\vert \mathfrak{F}_{\textrm{ext}}\big(K\times (a,b)\big)\Big] = \mathbb{E}^{k_1,k_2, (\ell,r),\vec{x}, \vec{y}, f, g}_{\mathbf{H}} \Big[F\big(\ell,r,\mathcal{\tilde L}^{k_1},\ldots, \mathcal{L}^{k_2}\big)\Big] \, ,
$$
where, on the right-hand side, $\ell = \mathfrak{a}$, $r= \mathfrak{b}$, $\vec{x} =(\mathcal{L}^i(\mathfrak{a}))_{i\in K}$, $\vec{y} =(\mathcal{L}^i(\mathfrak{b}))_{i\in K}$, $f = \mathcal{L}^{k_1-1}$ (or $+\infty$ if $k_2+1\notin \Sigma$), $g=\mathcal{L}^{k_2+1}$ (or $-\infty$ if $k_2+1\notin \Sigma$), and the curves $\mathcal{\tilde L}^{k_1},\ldots, \mathcal{L}^{k_2}$ have law $\mathbb{P}^{k_1,k_2, (\ell,r), \vec{x}, \vec{y}, f, g}_{\mathbf{H}}$.
\ed

\bl[Lemma 2.5 of \cite{CorHam16}]\label{lem:strongBGP}
Any line ensemble which enjoys the $\mathbf{H}$-Brownian Gibbs property also enjoys the strong $\mathbf{H}$-Brownian Gibbs property.
\el

Line ensembles with the $\mathbf{H}$-Brownian Gibbs property  benefit from certain  stochastic monotonicities. 

\bd[Domination of measure]\label{DomM}
Let $\mathcal{L}_1$ and $\mathcal{L}_2$ be two  $(\Sigma\times \Lambda)$-indexed line ensembles with respective laws $\mathbb{P}_1$ and $\mathbb{P}_2$. We say that $\mathbb{P}_{1}$ dominates $\mathbb{P}_{2}$ if there exists a coupling of  $\mathcal{L}_1$ and~$\mathcal{L}_2$ such that $\mathcal{L}^{j}_{1}(x)\geq \mathcal{L}^{j}_{2}(x)$ for all $j\in \Sigma$ and $x\in \Lambda$.
\ed


\bl[Stochastic monotonicity: Lemmas~$2.6$ and $2.7$ of \cite{CorHam16}]\label{Coupling1}
Fix finite intervals $K\subset \Sigma $ and $(a,b)\subset \Lambda$; and, for $i\in \{1,2\}$, vectors $\vec{x}_{i} = \big( x_i^{(k)}: k \in K \big)$ and  $\vec{y}_{i} = \big( y_i^{(k)}: k \in K \big)$ in $\R^K$ that satisfy $x^{(k)}_{2}\leq x^{(k)}_{1}$ and $y^{(k)}_{2}\leq y^{(k)}_{1}$ for $k\in K$; as well as measurable functions $f_i: (a,b)\to\RR \cup \{+\infty\} $ and $g_i: (a,b) \to \RR \cup \{-\infty\}$ such that $f_2(s)\leq f_1(s)$ and  $g_2(s)\leq  g_1(s)$ for $s\in (a,b)$. For $i \in \{1,2\}$, let
$\PP_i$ denote the law $\mathbb{P}^{k_1, k_2, (a,b), \vec{x}_i, \vec{y}_i, f_i,g_i}_{\mathbf{H}}$, so that a $\PP_i$-distributed random variable
 $\mathcal{L}_i= \{\mathcal{L}_i^{k}(s)\}_{k\in K, s\in(a,b)}$ is a $K\times (a, b)$-indexed line ensemble.
 If $\mathbf{H}:[0,\infty)\to \RR$ is convex, then $\mathbb{P}_{1}$ dominates $\mathbb{P}_{2}$ -- that is, a common probability space $(\Omega,\mathcal{B},\mathbb{P})$ may be constructed on which the two measures are supported such that, almost surely, $\mathcal{L}_{1}^{k}(s)\geq \mathcal{L}_{2}^{k}(s)$ for $k\in K$ and $s\in (a,b)$.
 \el

Since $\mathbf{H}_t(x)$ in \eqref{eq:Ht} is convex, Lemma~\ref{Coupling1} applies to any $\mathbf{H}_t(x)$-Brownian Gibbs  line ensemble.

Our next result recalls the  unscaled and scaled KPZ line ensemble constructed in~\cite{CorHam16}.

\bp[Theorem~$2.15$ of~\cite{CorHam16}]\label{NWtoLineEnsemble} Let $t>0$. There exists an $\NN\times \RR$-indexed line ensemble $\mathcal{H}_t =\{\mathcal{H}^{(n)}_{t}(x)\}_{n\in \NN, x\in \RR}$ such that:
\begin{enumerate}
\item the lowest indexed curve $\mathcal{H}^{(1)}_{t}(x)$ is equal in distribution (as a process in $x$) to the Cole-Hopf solution $\mathcal{H}^{\mathbf{nw}}(t,x)$ of KPZ started from the narrow wedge initial data;
 \item $\mathcal{H}_t$ satisfies the $\mathbf{H}_{1}$-Brownian Gibbs property (see Definition~\ref{LineEns});
 \item and the scaled KPZ line ensemble $\{\mathfrak{h}^{(n)}_t(x)\}_{n\in \NN,x\in \RR}$, defined by 
 \begin{align}\label{eq:UsilonNDef}
 \mathfrak{h}^{(n)}_t(x) \, = \,  t^{-1/3} \Big( \mathcal{H}^{(n)}_{t}\big(t^{2/3} x\big)+ t/24 \Big) \, ,
\end{align}
satisfies the $\mathbf{H}_{t}$-Brownian Gibbs property.
\end{enumerate}
\ep

This result shows that the lowest indexed curve $\mathfrak{h}_t^{(1)}$ in the scaled KPZ line ensemble has the law of the centred and scaled narrow wedge solution $\mathfrak{h}_t(x):=\mathfrak{h}_t(1,x)$ of the KPZ equation defined in~\eqref{eq:Upsilon}. This property is vital because it permits the Brownian Gibbs property of the ensemble to be brought to bear as a tool for analysing $\mathfrak{h}_t(x)$: one-point KPZ equation inputs thus lead to spatial regularity results for this process in~$x$.

We mention in passing that \cite[Theorem~$2.15$]{CorHam16} further asserts a similarity in compact regions between the scaled ensemble's curves and Brownian bridges, and it does so uniformly in $t \geq 1$. We will, however, make no use of this assertion.

\subsection{Input results for the KPZ equation}\label{sec:KPZeqinputs}

Recall the notation  $\mathfrak{h}_t(x) =\mathfrak{h}_t(1,x)$.

\bp[Stationarity]\label{ppn:Stationarity}
The one-point distribution of $\mathfrak{h}_{t}(x)+\frac{x^2}{2}$ is independent of $x$ and converges weakly to a limit as $t\nearrow \infty$. Furthermore, the processes $x\mapsto \mathfrak{h}_{t}(x)+\frac{x^2}{2}$ and $x\mapsto \mathfrak{h}_{t}(-x)+\frac{x^2}{2}$ are equal in law.
\ep
\begin{proof}
The first sentence follows immediately from  \cite[Corollary~1.3 and Proposition~1.4]{Amir11}. The second is a straightforward consequence of the reflection invariance of space-time white noise.
\end{proof}

The next two results are variants  for the KPZ equation  of the FKG inequality. They can be proved by appealing to the standard FKG inequality for prelimiting models such as ASEP.

\bp[Positive association and the FKG inequality]\label{FKGProp} For any $k\in \N$, $t_1,\ldots ,t_k\geq 0$, $x_1,\ldots , x_k\in \RR$ and $s_1, \ldots ,s_k\in \RR$,
 \begin{align}\label{eq:FKGStatement}
 \mathbb{P}\Big(\bigcap_{\ell=1}^{k}\big\{\mathfrak{h}_{t_\ell}(x_\ell)\leq s_\ell\big\}\Big)\geq \prod_{\ell =1}^{k}\mathbb{P}\Big(\mathfrak{h}_{t_\ell}(x_\ell)\leq s_{\ell}\Big).
 \end{align}
In particular, for $t_1, t_2\in \RR_{>0}$, $x_1, x_2 \in \RR$ and $s_1, s_2\in \RR$,
\begin{align}\label{eq:RevFKG}
\mathbb{P}\Big(\mathfrak{h}_{t_1}(x_1)>s_1, \mathfrak{h}_{t_2}(x_2)>s_2\Big)\geq \mathbb{P}\Big(\mathfrak{h}_{t_1}(x_1)>s_1\Big)\mathbb{P}\Big(\mathfrak{h}_{t_2}(x_2)>s_2\Big).
\end{align}

 \ep
 \begin{proof}
 This result follows from \cite[Proposition~1]{CQ11} after centring and scaling. 
 \end{proof}

The second FKG result asserts that conditioning the KPZ equation solution at time $t$ on a larger (or smaller) value increases (or decreases) the conditional expectation at a later time $\alpha t$.

\bp[Monotonicity under conditioning]\label{ppn:CondEIn}
For $t>0$, $\alpha>1$, $x_1, x_2, r\in\RR$ and $u~>~v~\in~\RR$,
\begin{align}
\mathbb{P}\big(\mathfrak{h}_{t}(1,x_1)>v\big)&\times\mathbb{P}\big(\mathfrak{h}_{t}(\alpha,x_2)>r, \mathfrak{h}_{t}(1,x_1)>u\big)\\&\geq \mathbb{P}\big(\mathfrak{h}_{t}(1,x_1)>u\big)\times \mathbb{P}\big(\mathfrak{h}_{t}(\alpha,x_2)>r, \mathfrak{h}_{t}(1,x_1)>v\big),\\
\mathbb{P}\big(\mathfrak{h}_{t}(1,x_1)\leq u\big)&\times \mathbb{P}\big(\mathfrak{h}_{t}(\alpha,x_2)>r, \mathfrak{h}_{t}(1,x_1\big)\leq v)\\&\geq \mathbb{P}\big(\mathfrak{h}_{t}(1,x_1)\leq v\big)\times \mathbb{P}\big(\mathfrak{h}_{t}(\alpha,x_2) > r , \mathfrak{h}_{t}(1,x_1)\leq u\big).
\end{align}

%
\ep

\begin{proof}
The two bounds are proved similarly hence we only treat the first. Consider the SHE with Dirac delta initial data (the proof works for general initial data too). We {\em claim} that,  for $s<t$,
\begin{align}
  \mathbb{P}(\mathcal{Z}(s,x_1)>e^v)&\times \mathbb{P}\big(\mathcal{Z}(t,x_2)>e^r, \mathcal{Z}(s,x_1)>e^u\big)
  \\&\geq \mathbb{P}(\mathcal{Z}(s,x_1)>e^u)\times \mathbb{P}\big(\mathcal{Z}(t,x_2)>e^r, \mathcal{Z}(s,x_1)>e^v\big) \, .
\end{align}
The proposition's first bound follows from this claim by taking logarithms, centring and scaling.

The proof of the claim is almost the same as that of \cite[Proposition~1]{CQ11}. It relies on (1) the results of \cite{Amir11}, which approximate $\mathcal{Z}(t,x)$ by  the microscopic Cole-Hopf (or G\"{a}rtner) transform $\mathcal{Z}^{\epsilon}(t,x)$  of  ASEP under weak asymmetry scaling; and (2) the FKG inequality for ASEP, a bound due to this model's graphical construction  (see the proof of \cite[Proposition~1]{CQ11} or \cite{Liggett85,Liggett99} for details). The FKG inequality implies that the desired identity claimed for $\mathcal{Z}$ holds for $\mathcal{Z}^{\epsilon}$; by convergence, this bound holds in the limit as well.
\end{proof}


The next result is the composition law for the KPZ equation. By its use, aspects of the two-time distribution will be inferred from  the fixed-time  narrow wedge KPZ solution. The mainstays of the result's proof are the random semi-group property and the time-reversal symmetry enjoyed by the SHE. Our presentation of the derivation is brief in view of its similarities to the proof of \cite[Lemma 1.18]{CorHam16}. It is reasonable to wonder whether the composition law can be applied for study more than two disjoint times. For such a purpose, there is a composition law but it cannot be formulated purely in terms of the narrow wedge solution: it would be necessary to understand the joint distribution of the KPZ equation started from various shifted narrow wedges. In the $t\nearrow \infty$ limit, this data is expected to be measurable with respect to the space-time Airy sheet.

For $t>0$, define a $t$-indexed composition map $I_t(f,g)$ between two functions $f(\cdot)$ and $g(\cdot)$:
\begin{align}\label{eq:Icomp}
I_t(f,g):=t^{-1/3}\log \int^{\infty}_{-\infty} e^{t^{1/3}\big(f(t^{-2/3}y) +g(-t^{-2/3} y)\big)} dy \, .
\end{align}

\bp[Composition law]\label{ppn:TimeEvl}
For any fixed $t>0$ and $\alpha>1$, there exists a spatial process $\mathfrak{h}_{\alpha t\downarrow t}(\cdot)$ supported on the same probability space as the KPZ equation solution such that:
\begin{enumerate}
\item $\mathfrak{h}_{\alpha t\downarrow t}(\cdot)$ is distributed according to the law of the process $\mathfrak{h}_{t}(\alpha-1, \cdot)$;
\item $\mathfrak{h}_{\alpha t\downarrow t}(\cdot)$ is independent of $\mathfrak{h}_t(\cdot) :=\mathfrak{h}_t(1,\cdot)$; and
\item $\mathfrak{h}_t(\alpha,0) = I_{t}\big(\mathfrak{h}_{t}, \mathfrak{h}_{\alpha t\downarrow t}\big)$.
\end{enumerate}
\begin{rem}
The notation $\mathfrak{h}_{\alpha t\downarrow t}(\cdot)$  will be, in Sections \ref{sec:remote} and \ref{sec:adjacent}, complemented by similar notation $\mathfrak{h}_{0\uparrow t}(\cdot)$ in place of $\mathfrak{h}_t(1,\cdot)$ and $\mathfrak{h}_{0\uparrow \alpha t}(\cdot)$ in place of $\mathfrak{h}_t(\alpha,\cdot)$. The idea prompting this usage is that the solution from time $0$ to $\alpha t$ can be constructed by combining the solution from time zero to time $t$ with the solution from $\alpha t$ to $t$ -- in a sense which will be clear in the proof that we now give.
\end{rem}

%
%
\ep


\begin{proof}[Proof of Proposition~\ref{ppn:TimeEvl}]
For $s<t$ and $x,y\in \RR$, let $\mathcal{Z}^{\mathbf{nw}}_{s,x}(t,y)$ be the solution at time $t$ and position $y$ of the SHE started at time $s$ with Dirac delta initial data at position $x$. As these four parameters vary, we assume that all solutions are coupled on a probability space upon which their common space-time white noise is defined. It is due to this and the linearity of the SHE that
\begin{align}\label{eq:Conv}
\mathcal{Z}^{\mathbf{nw}}(t,y) := \mathcal{Z}^{\mathbf{nw}}_{0,0}(t,y) = \int_{\RR} \mathcal{Z}^{\mathbf{nw}}_{0,0}(s,x)\mathcal{Z}^{\mathbf{nw}}_{s,x}(t,y) dx \, ,
\end{align}
 where $\mathcal{Z}^{\mathbf{nw}}_{0,0}(s,x)$ and $\mathcal{Z}^{\mathbf{nw}}_{s,x}(t,y)$ are independent.  Though this convolution formula is known, we did not find a careful proof in the literature. Thus, for completeness we will provide such a proof below. Assume for the moment its validity. We also need an {\em interchange} property of the SHE: namely that, for $s<t$ and $y\in \RR$ fixed,  $\mathcal{Z}^{\mathbf{nw}}_{s,x}(t,y)$ is equal in law as a process in $x$ to $\mathcal{Z}^{\mathbf{nw}}_{s,y}(t,x)$ -- the change between the two expressions is in the interchange of $x$ and $y$. This is readily shown as a consequence of the chaos series formula \eqref{eq:ZEq} and the invariance in law of $\xi$ under reflections.
Expressing the convolution and interchange properties in terms of $\mathfrak{h}_t$ immediately yields the proposition.

 We now return to show \eqref{eq:Conv}. There are many (now standard) ways to prove this (e.g., from the Feynman-Kac representation of the solution of the SHE). We proceed via the chaos series for the SHE (see \cite{C18} or \cite{alberts2014} for background). For any $0\leq s<t$ and $x,y\in \mathbb{R}$, $\mathcal{Z}^{\mathbf{nw}}_{s,x}(t,y)$ is given as a special case of the following chaos series expansion (see Theorem~2.2 of \cite{C18}):
\begin{align}
\mathcal{Z}^{\mathbf{nw}}_{s,x}(t,y) = \sum_{k=0}^{\infty} \int_{\Delta_{k}(s,t)} \int_{\mathbb{R}^{k}} P_{k;s,x;t,y}(\vec{s},\vec{x})d\xi^{\otimes_k}(\vec{s},\vec{x})\label{eq:ZEq}.
\end{align}
Here we have used the following conventions. We write $\vec{s}=(s_1,\ldots, s_k)\in \RR^{k}_{\geq 0}$,  $\vec{x}= (x_1,\ldots, x_k)\in \RR^{k}$ and define the set of ordered times
$$\Delta_{k}(s,t) = \{\vec{s}: s\leq s_1\leq s_2\leq \ldots \leq s_k\leq t\}.$$
The integration in \eqref{eq:ZEq} is a multiple It\^o stochastic integral against the white noise $\xi$ and the term $P_{k;s,x;t,y}(\vec{s},\vec{x})$ is the density function for a one-dimensional Brownian motion starting from $(s,x)$ to go through the time-space points $(s_1,x_1),\ldots, (s_k,x_k)$ and end up at $(t,y)$. This transition density has the following product formula using the Gaussian heat kernel $p(s,y) := (2\pi s)^{-1/2}\exp(-y^2/2s)$ and the conventions $s_0=s$,  $s_{k+1}=t$, $x_0=x$ and $x_{k+1}=y$:
\begin{align}
P_{k;s,x;t,y}(\vec{s},\vec{x}) = \prod_{i=0}^{k} p(s_{i+1}-s_i,x_{i+1}-x_i).
\end{align}
The heat kernel $p(\cdot,\cdot)$ satisfies the simple convolution identity
\begin{align}\label{eq:pconv}
p(t,x) = \int p(s,y) p(t-s,x-y)  dy
\end{align}
for any $0\leq s<t$. It is from this identity that \eqref{eq:Conv} will follow, as we now show.

Fix $r\in (s,t)$.  Without changing the expressed value, we may replace $\int_{\Delta_k(s,t)}$  in \eqref{eq:ZEq} by the quantity $\sum_{i=0}^k \int_{\Delta_k(s,t)} \mathbf{1}_{s_i\leq r<s_{i+1}}$ (the sum of indicator functions gives the value one). This implies that
\begin{align}\label{eq:indinsert}
\mathcal{Z}^{\mathbf{nw}}_{s,x}(t,y) =  \sum_{k=0}^{\infty}\sum_{i=0}^k \int_{\Delta_{k}(s,t)} \int_{\mathbb{R}^{k}} \mathbf{1}_{s_i\leq r<s_{i+1}}  P_{k;s,x;t,y}(\vec{s},\vec{x})d\xi^{\otimes_k}(\vec{s},\vec{x}).
\end{align}
Using \eqref{eq:pconv}, we may rewrite
$$
 \mathbf{1}_{s_i\leq r<s_{i+1}}  P_{k;s,x;t,y}(\vec{s},\vec{x}) =  \mathbf{1}_{s_i\leq r<s_{i+1}}  \int_{\RR} P_{i;s,x;r,z}(\vec{s}_{[1,i]},\vec{x}_{[1,i]})P_{k-i;r,z;t,y}(\vec{s}_{[i+1,k]},\vec{x}_{[i+1,k]})dz
$$
where, for $1\leq a\leq b\leq k$,  $\vec{s}_{[a,b]}$ denotes $(s_{a},\ldots, s_b)$ and likewise for $\vec{x}$.
Inserting this into \eqref{eq:indinsert}, replacing $\int_{\Delta_k(s,t)} \mathbf{1}_{s_i\leq r<s_{i+1}}$ by $\int_{\Delta_i(s,r)}\int_{\Delta_{k-i}(r,t)}$ and relabeling $\vec{s}_{[1,i]}=\vec{u}$, $\vec{s}_{[i+1,k]}=\vec{v}$, $\vec{x}_{[1,i]}=\vec{a}$, $\vec{x}_{[i+1,k]}=\vec{b}$,
we find that
\begin{align}
\mathcal{Z}^{\mathbf{nw}}_{s,x}(t,y) =  \sum_{k=0}^{\infty}\sum_{i=0}^k \int_{\Delta_i(s,r)}\int_{\Delta_{k-i}(r,t)} \int_{\mathbb{R}^{i}}\int_{\mathbb{R}^{k-i}} \int_{\RR}& P_{i;s,x;r,z}(\vec{u},\vec{a})P_{k-i;r,z;t,y}(\vec{v},\vec{b})dz\\& \times  d\xi^{\otimes_i}(\vec{u},\vec{a})d\xi^{\otimes_{k-i}}(\vec{v},\vec{b}).
\end{align}
Here we also used that the white noise integration can be split since the times range over disjoint intervals.
Making the change of variables $j=k-i$, the double sum $\sum_{k=0}^{\infty}\sum_{i=0}^k$ can be replaced by $\sum_{i=0}^{\infty}\sum_{j=0}^\infty$. Finally, bringing the integral in $z$ to the outside, we recognize that what remains is the product of two chaos series of the form \eqref{eq:ZEq}, which is exactly the result claimed in  \eqref{eq:Conv}. The resumming and reordering of integrals is readily justified since all sums are convergent in $L^{2}$ (with respect to the probability space on which $\xi$ is defined -- see, for example, \cite[Theorem 2.2]{C18} for details). Independence of $ \mathcal{Z}^{\mathbf{nw}}_{0,0}(s,x)$ and $\mathcal{Z}^{\mathbf{nw}}_{s,x}(t,y)$ in \eqref{eq:Conv} follows immediately from their being defined with respect to disjoint portions of the space-time white noise.

\end{proof}

Our final inputs are one-point tail bounds for the KPZ equation, recently proved in \cite{CG18a,CG18b}.

\bp[Uniform lower-tail bound]\label{ppn:onepointlowertail}
For any  $t_0>0$, there exist $s_0 = s_0(t_0)>0$ and $c =c(t_0)>0$ such that, for  $t>t_0$, $s>s_0$ and $x\in \RR$,
\begin{align}
\mathbb{P}\Big(\mathfrak{h}_{t}(x)+\frac{x^2}{2} \leq -s \Big) &\leq \exp\big(-c s^{5/2}\big) \, . \label{eq:LowTailUp}
\end{align}
\ep
\begin{proof}
By Proposition \ref{ppn:Stationarity}, we may take $x=0$. The result then follows from~\cite[Theorem 1.1]{CG18a} because $\Upsilon_t$ in the quoted result is the same as $\mathfrak{h}_t(0)$  up to a constant change of scale.
\end{proof}

\bp[Uniform upper-tail bound]\label{ppn:onepointuppertail}
For any $t_0>0$, there exist $s_0 = s_0(t_0)>0$ and $c_1=c_1(t_0)>c_2=c_2(t_0)>0$ such that, for $t\geq t_0$, $s>s_0$ and $x\in \RR$,
\begin{align}
\exp\big(-c_1 s^{3/2}\big) \leq \mathbb{P}\Big(\mathfrak{h}_t(x)+\frac{x^2}{2}\geq s\Big) \leq \exp\big(-c_2 s^{3/2}\big) \, .
\end{align}
\ep
\begin{proof}
By Proposition \ref{ppn:Stationarity}, we may take $x=0$. The result follows from \cite[Theorem~1.10]{CG18b}. 
\end{proof}

\section{Analogous results for the KPZ fixed point}\label{sec:KPZfixedpoint}

It is believed that $\mathfrak{h}_t(\alpha, x)$ converges in the limit of high~$t$ and as a time-space process to the narrow-wedge initial data solution of the \emph{KPZ fixed point}. This important universal object is a Markov process on random functions whose existence was conjectured in \cite{Corwin2015}; it has  recently been constructed in \cite{KPZfixed} for any fixed initial data via its transition probabilities, and in \cite{Landscape} simultaneously for all initial data via the \emph{Airy sheet}. A special case of the putative universality of the KPZ fixed point is the  conjecture -- made, for example, in \cite[Conjecture 1.5]{Amir11} --  that the process  $x\mapsto 2^{1/3}\big(\mathfrak{h}_t(x)+\tfrac{x^2}{2}\big)$ converges to $x\mapsto \mathcal{A}(x)$, where $\mathcal{A}(\cdot)$ is the Airy$_2$ process introduced in~\cite{Prahofer2002}. A similar though stronger assertion has been expressed in~\cite[Conjecture 2.17]{CorHam16} for the KPZ line ensemble: namely that, after adding in the parabolic shift $\tfrac{x^2}{2}$ and scaling by $2^{1/3}$, this ensemble converges in the limit of high $t$ to the Airy line ensemble constructed in \cite{CH14}.

In this section, we review our results and methods through the lens offered by zero temperature -- that is, we discuss the counterpart problems and solutions in the limit where the KPZ time parameter~$t$ becomes high. Positive temperature structures have zero temperature counterparts with simple and vivid interpretations, and there has been much recent effort to understand counterpart problems in the limiting $t \nearrow \infty$ case. Thus it is that, while the upcoming discussion is logically needless for comprehension of this paper's results and proofs, we hope that this summary may aid the reader's perspective on our results, their derivations and their relation to recent advances.

We start by noting the $t \nearrow \infty$ counterpart to the composition law Proposition \ref{ppn:TimeEvl}.  For functions $f$ and $g$, the high $t$ limit of $I_t(f,g)$ defined in \eqref{eq:Icomp} is the variational problem
\begin{align}\label{eq:Ilimits}
I_t(f,g):=t^{-1/3}\log\Big(t^{2/3} \int^{\infty}_{-\infty} e^{t^{1/3}\big(f(y) +g(-y)\big)} dy\Big) \, \xrightarrow[t\nearrow \infty]{} \, \sup_{y\in \mathbb{R}} \Big\{f(y)+g(-y)\Big\}=:I_{\infty}(f,g) \, .
\end{align}
This Laplace method type of transition from the logarithm of the integral of an exponential in $I_t$ to the supremum in $I_{\infty}$ is the hallmark of passing from positive to zero temperature.

The form of the limiting composition law permits counterparts to our principal results to be formatted in terms of  two independent Airy$_2$ processes, $\mathcal{A}$ and $\mathcal{\widetilde{A}}$. Define $\mathcal{B}(x):= 2^{-1/3} \mathcal{A}(x) -\tfrac{x^2}{2}$ and $\mathcal{\widetilde{B}}(x):= 2^{-1/3} \mathcal{\widetilde{A}}(x) -\tfrac{x^2}{2}$. Set $\mathcal{B}_1 = \mathcal{B}(0)$ and, for $\beta>0$, define
\begin{align}\label{eq:mathcalB}
\mathcal{B}_{1+\beta} \, = \, \sup  \Big\{\mathcal{B}(x)+ \beta^{1/3} \mathcal{\widetilde{B}}\big(-x \beta^{-2/3}\big) : x \in \R \Big\} \, .
\end{align}
The scaling here applied to the term $\mathcal{\widetilde{B}}$
results in a process in $x$ suitable for the description of scaled last passage percolation values over a scaled duration equal to $\beta$; the resulting term may be viewed as
a time $\beta$ version of the Airy$_2$ process.

The pair $\big(\mathcal{B}_1, \mathcal{B}_{1+\beta}\big)$ is counterpart to $\big(\mathfrak{h}_t(1,0), \mathfrak{h}_{t}(1+\beta,0) \big)$. In law, this pair has the joint distribution of the narrow-wedge initial data KPZ fixed point at
the space-time point pair $(0,1)$ and $(0,1+\beta)$. This distributional equality holds for  given $\beta$, and no such assertion is being made regarding processes in $\beta$.

Besides the composition law, the other key tools described in Section \ref{sec:Tools} have direct analogues for the KPZ fixed point. The KPZ line ensemble is replaced by the Airy line ensemble \cite{CH14}. After a parabolic shift, the latter enjoys the version of the Brownian Gibbs property given formally by $\mathbf{H}(x) = \infty \mathbf{1}_{x\geq 0}$,
in which intersection of adjacent curves is forbidden.
Stochastic monotonicity is unaffected by the limit $t \nearrow \infty$. The KPZ fixed point also satisfies the same stationarity properties and FKG inequalities, while the tail bounds in Propositions \ref{ppn:onepointlowertail} and \ref{ppn:onepointuppertail} are replaced by such bounds for the GUE Tracy-Widom distribution \cite{TracyWidom94}.

Combining these alterations with the noted transition from $I_t$ to $I_{\infty}$ composition laws, we now state zero-temperature counterparts to our main theorems. We do not gives proofs of the statements in this article, although our arguments offer templates for such proofs. Indeed, the zero-temperature context is an alternative and, in certain regards, simpler mode for interpreting statements and proofs -- the presentation of the following statements is intended to aid the reader who wishes to view this article through the zero-temperature prism.




First, two assertions concerning exponents for remote and adjacent two-time correlation. The recent works~\cite{FO18} and~\cite{BG18} offer corresponding theorems for certain last passage percolation models.

\textbf{Theorem \ref{thm:Main1} analogue:} There exist $c_1, c_2 >0$ such that, for $\alpha>2$,
\begin{align}
c_1 \alpha^{-1/3} \leq \mathrm{Corr}\big(\mathcal{B}_1, \mathcal{B}_{\alpha}\big) \leq c_2 \alpha^{-1/3}.
\end{align}

\textbf{Theorem \ref{thm:Main2} analogue:} There exist $c_1,c_2>0$ such that, for $\beta\in (0,\frac{1}{2})$,
\begin{align}
c_1 \beta^{2/3}\leq 1-\mathrm{Corr}\big(\mathcal{B}_1, \mathcal{B}_{1+\beta}\big) \leq c_2\beta^{2/3}.
\end{align}

The next spatial regularity result has been proved for scaled Brownian last passage percolation -- see \cite[Theorem $1.1$]{Hammond18b}.

\textbf{Theorem \ref{t.locreg} analogue:} There exist $s_0>0$ and $c>0$ such that, for $x \in \R$, $s\geq s_0$ and $\e \in (0,1]$,
\begin{align}
 \PP \, \bigg( \, \sup_{z \in [x,x+\e]} \big\vert \mathcal{B}(z) - \mathcal{B}(x)\big\vert \, \geq \,  s\e^{1/2}  \, \bigg) \leq   \exp \big( - c s^{3/2} \big) \, .
 \end{align}

The next modulus of continuity inference follows -- see \cite[Theorem $1.3$]{Hammond18b} for such a result which holds uniformly over choices of $\mc{B}$ arising from a large class of initial data, in place of the narrow wedge considered here.

\textbf{Corollary \ref{thm:SpaceHolder} analogue:}
For any interval $[a,b]\subset \RR$, define
$$
\mathcal{C} :=\sup_{x_1\neq x_2\in [a,b]} \, |x_1-x_2|^{-1/2}\bigg( \log\frac{|b-a|}{|x_1-x_2|}\bigg)^{-2/3}\big\vert \mathcal{B}(x_1)- \mathcal{B}(x_2)\big\vert \, .
$$
Then there exist $s_0=s_0(|b-a|)>0$ and $c=c(|b-a|)>0$   such that, for $s\geq s_0$,
$$
\mathbb{P}\big(\mathcal{C}>s \big)\leq \exp\big(-cs^{3/2}\big) \, .
$$

The spatial-temporal modulus of continuity for the KPZ fixed point is also probed in~\cite[Proposition 1.6]{Landscape}, a result which implies the next stated tail on the law of fluctuation between nearby times at a given location at zero temperature. A similar result is obtained for Poissonian last passage percolation in~\cite{HammondSarkar}.

\textbf{Theorem \ref{ppn:HDivTail} analogue:} There exist $s_0>0$ and $c_1, c_2, c_3>0$   such that, for $s>s_0$ and $\beta \in (0,1/2)$,
 \begin{align}
\exp\big(-c_1s^{3/2}\big) \leq \mathbb{P}\Big(\mathcal{B}_{1+\beta}-\mathcal{B}_1 \geq \beta^{1/3}s\Big)&\leq \exp\big(-c_2 s^{3/4}\big) \, ,\\
 \mathbb{P}\Big(\mathcal{B}_{1+\beta}-\mathcal{B}_1\leq -\beta^{1/3}s\Big)&\leq \exp\big(-c_3s^{3/2}\big) \, .
 \end{align}

\section{Spatial process tail bounds}\label{sec:spatialtails}
In this section, we prove bounds on the tails of various functionals of the spatial process $\mathfrak{h}_t(\cdot)$, such as its infimum, supremum and increment on a fixed interval.
Four propositions will be stated and proved, the first two concerning global properties of this process, the later two addressing local ones; in the latter vein, we will also prove the local spatial regularity Theorem~\ref{t.locreg} here. The proofs in this section rely upon: (1) the KPZ line ensemble Brownian Gibbs property, Proposition~\ref{NWtoLineEnsemble}; (2) the monotone coupling Lemma \ref{Coupling1}; (3) Brownian bridge calculations; and (4) the one-point tail bounds Propositions \ref{ppn:onepointlowertail} and \ref{ppn:onepointuppertail}. Note that we do not use the composition law in this section.

In the proofs of the first two propositions,
the constant $c=c(t_0,\nu)>0$ may change value from line to line and even between consecutive inequalities. Moreover, a bound involving $c$ and $s$ implicitly asserts that there exist $s_0=s_0(t_0,\nu)$ and $c(t_0,\nu)$ such that, for all $s\geq s_0$ and $t\geq t_0$, the recorded bound holds. The explicit form is used in the propositions'  statements, and the abbreviating device is employed in their proofs. Recall also that $\neg \mathcal{E}$ denotes the complement of the event $\mathcal{E}$.

\bp\label{ppn:LowerTail}
For any $t_0>0$ and $\nu\in (0,1)$, there exist $s_0=s_0(t_0,\nu)>0$ and  $c=c(t_0,\nu)>0$ such that, for $t\geq t_0$ and $s>s_0$,
\begin{align}
\mathbb{P}(\mathcal{A})\leq \exp\big(-c s^{5/2}\big)\qquad \textrm{where}\qquad \mathcal{A}:=\Big\{\inf_{x\in \RR}\Big( \mathfrak{h}_t(x)+\frac{(1+\nu) x^2}{2}\Big) \leq -s \Big\} \, .\label{eq:InfU}
\end{align}
\ep

\begin{proof}

 For $n\in \ZZ$, define $\zeta_n := n/s$ and
\begin{align}
\mathcal{E}_n &:= \Big\{\mathfrak{h}_{t}(\zeta_n)\leq -\frac{(1+\tfrac{\nu}{2})\zeta^2_n}{2}- (1-\epsilon)s\Big\} \, ,\\
\mathcal{F}_n &:= \Big\{\mathfrak{h}_{t}(y)\leq -\frac{(1+\nu)y^2}{2}- (1-\epsilon/2)s, \quad \forall y \, \in [\zeta_n, \zeta_{n+1}]\Big\} \, .
\end{align}
By the union bound,
 \begin{align}
 \mathbb{P}(\mathcal{A})\leq \sum_{n\in \ZZ}\mathbb{P}(\mathcal{E}_n)+ \sum_{n\in \ZZ}\mathbb{P}\Big(\mathcal{A}\cap\Big\{\bigcap_{m\in \ZZ} \neg\mathcal{E}_n \Big\}\cap \mathcal{F}_n\Big)+ \mathbb{P}\Big(\mathcal{A}\cap \Big\{\bigcap_{n\in \ZZ}  \neg \mathcal{E}_n \Big\}\cap \Big\{\bigcap_{n\in \ZZ} \neg \mathcal{F}_n \Big\}\Big) \, . \qquad\label{eq:SplitProb4}
 \end{align}
 We bound each of the three right-hand summands. Note that $\mathcal{A}\cap \neg \mathcal{F}_n=\emptyset$ for all $n\in \ZZ$. Hence,
 \begin{align}\label{eq:CapProbBd}
 \mathbb{P}\Big(\mathcal{A}\cap \Big\{\bigcap_{n\in \ZZ}  \neg \mathcal{E}_n \Big\}\cap \Big\{\bigcap_{n\in \ZZ} \neg \mathcal{F}_n \Big\}\Big) = 0 \, .
 \end{align}
  The first summand in  \eqref{eq:SplitProb4} is bounded by
 \begin{align}
 \sum_{n\in \ZZ}\mathbb{P}(\mathcal{E}_n)\leq \sum_{n\in \ZZ}\exp\Big(-c\big(\tfrac{\nu}{2}\zeta^2_n+s\big)^{5/2}\Big) \leq \sum_{n\in \ZZ}\exp\Big(-c\big(\tfrac{\nu}{2}\zeta^2_n\big)^{5/2} -cs^{5/2}\Big) \, . \label{eq:SeqProbBd}
\end{align}
The first inequality here is due to the bound on $\mathbb{P}(\mathcal{E}_n)$ that results from Proposition \ref{ppn:onepointlowertail}, while  the second is obtained by applying the reverse Minkowski inequality -- this being the bound that, for any $a,b>0$, $(a+b)^{5/2}\geq a^{5/2}+b^{5/2}$. The right-hand sum in \eqref{eq:SeqProbBd} is now bounded above by a suitable multiple of the corresponding integral:
 \begin{align}\label{eq:SumIntConv}
 \sum_{n\in \ZZ}\exp\Big(-c\big(\tfrac{\nu}{2}\zeta^2_n)^{5/2}\Big)\leq s^{5}\int^{\infty}_{-\infty} \exp(- c|x|^{5}) dx \leq c s^5. 
 \end{align}
The right-hand side of \eqref{eq:SeqProbBd}  is thus seen to be at most $c s^{5} e^{-c's^{5/2}}$, which is in turn bounded above by $e^{-c''s^{5/2}}$, where in this instance,  we distinguish between constants in an attempt at avoiding confusion.
We thus find that
  \begin{align}\label{eq:ESum}
  \sum_{n\in \ZZ}\mathbb{P}(\mathcal{E}_n)\leq \exp(-c s^{5/2}) \, .
  \end{align}

The second right-hand term in~\eqref{eq:SplitProb4} remains to be addressed.  Indeed, since $\mathcal{A}\cap\big\{\bigcap_{m\in \ZZ} \neg\mathcal{E}_n \big\}\cap \mathcal{F}_n$ is contained in $\neg \mathcal{E}_{n} \cap \neg\mathcal{E}_{n+1} \cap \mathcal{F}_n$ for all $n\in \ZZ$, the next bound suffices in light of  \eqref{eq:CapProbBd} and \eqref{eq:ESum} to prove \eqref{eq:InfU} and hence the proposition: 
\begin{align}\label{eq:EESum}
\sum_{n\in \ZZ}\mathbb{P}\Big( \neg \mathcal{E}_{n} \cap \neg\mathcal{E}_{n+1} \cap \mathcal{F}_n\Big)\leq \exp(- cs^{3}).
\end{align}
To bound $\mathbb{P}\big( \neg \mathcal{E}_{n} \cap \neg\mathcal{E}_{n+1} \cap \mathcal{F}_n\big)$ for all $n\in \ZZ$, we will make use of Proposition~\ref{NWtoLineEnsemble}, which shows that the lowest indexed curve $\mathfrak{h}^{(1)}_t(\cdot)$ of the KPZ line ensemble has the same distribution as $\mathfrak{h}_t(\cdot)$. Owing to this, we may replace $\mathfrak{h}_{t}(\cdot)$ in the definitions of $\mathcal{E}_n$ and $\mathcal{F}_n$ by $\mathfrak{h}^{(1)}_t(\cdot)$. We will work with these modified definitions for the rest of the proof (which is to say, the derivation~\eqref{eq:EESum}); we will thus be able to use the $\mathbf{H}_t$-Brownian Gibbs property associated with the KPZ line ensemble.

Let $\mathfrak{F}_n=\mathfrak{F}_{\mathrm{ext}}\big(\{1\}\times (\zeta_n, \zeta_{n+1})\big)$ be the $\sigma$-algebra generated by $\{\mathfrak{h}^{(n)}_t(x)\}_{n\in \NN,x\in \RR }$ outside $\{\mathfrak{h}^{(1)}_t(x): x\in (\zeta_n,\zeta_{n+1})\}$.  By the strong $\mathbf{H}_t$-Brownian  Gibbs property  Lemma \ref{lem:strongBGP} for $\{\mathfrak{h}^{(n)}_t\}_{n\in \NN, x\in \RR}$,
 \begin{align}
&&&&&&&&\mathbb{P}\big( \neg \mathcal{E}_{n} \cap \neg\mathcal{E}_{n+1} \cap \mathcal{F}_n \big) = \mathbb{E}\Big[\mathbf{1}_{\neg \mathcal{E}_{n} \cap \neg\mathcal{E}_{n+1}}\cdot \mathbb{E}[\mathcal{F}_n|\mathfrak{F}_n]\Big]
= \mathbb{E}\Big[ \mathbf{1}_{\neg \mathcal{E}_{n} \cap \neg\mathcal{E}_{n+1}}\cdot \mathbb{P}_{\mathbf{H}_t}(\mathcal{F}_n)\Big] \qquad\label{eq:PCE}
\end{align}
where $\mathbb{P}_{\mathbf{H}_t}:=\mathbb{P}^{1,1,(\zeta_n, \zeta_{n+1}), \mathfrak{h}^{(1)}_t(\zeta_n), \mathfrak{h}^{(1)}_t(\zeta_{n+1}), +\infty, \mathfrak{h}^{(2)}_t}_{\mathbf{H}_t}$.
By Lemma~\ref{Coupling1}, there exists a monotone coupling between $\mathbb{P}_{\mathbf{H}_t}$  and $\tilde{\mathbb{P}}_{\mathbf{H}_t}:=\mathbb{P}^{1,1,(\zeta_n, \zeta_{n+1}), \mathfrak{h}^{(1)}_t(\zeta_n), \mathfrak{h}^{(1)}_t(\zeta_{n+1}), +\infty, -\infty}_{\mathbf{H}_t}$ such that
 \begin{align}\label{eq:Monotone11}
 \mathbb{P}_{\mathbf{H}_t}(\mathcal{F}_n)\leq \tilde{\mathbb{P}}_{\mathbf{H}_t}(\mathcal{F}_n).
\end{align}
note that $\tilde{\mathbb{P}}_{\mathbf{H}_t}$ is same as $\mathbb{P}^{(\zeta_n, \zeta_{n+1}), \mathfrak{h}^{(1)}_t(\zeta_n), \mathfrak{h}^{(1)}_t(\zeta_{n+1})}_{\mathrm{free}}$ which is the law of a Brownian bridge. For $n\in \ZZ$, define $\theta_n:= (1-\epsilon)s + \frac{(1+2^{-1}\nu)}{2}\zeta^2_n$. Then $\mathbf{1}_{\neg \mathcal{E}_{n} \cap \neg\mathcal{E}_{n+1}} = \mathbf{1}_{\mathfrak{h}^{(1)}_{t}(\zeta_n)> -\theta}\cdot \mathbf{1}_{\mathfrak{h}^{(1)}_{t}(\zeta_{n+1})> -\theta_{n+1}}$.
Under Brownian bridge law  on the interval $(\zeta_n, \zeta_{n+1})$, the probability of $\mathcal{F}_n$ increases with pointwise decrease of the sample paths at the endpoints. Thus,
\begin{align}\label{eq:BB}
\mathbf{1}_{\neg \mathcal{E}_{n} \cap \neg\mathcal{E}_{n+1}}\cdot \tilde{\mathbb{P}}_{\mathbf{H}_t}(\mathcal{F}_n) \, \leq \, \mathbb{P}^{(\zeta_n, \zeta_{n+1}), -\theta_n, -\theta_{n+1}}_{\mathrm{free}}(\mathcal{F}_n) \, .
\end{align}
For a Brownian bridge $B(\cdot)$ with $B(\zeta_n)= -\theta_n$ and $B(\zeta_{n+1}) = -\theta_{n+1}$,
\begin{align}
\mathbb{P}^{(\zeta_n, \zeta_{n+1}), -\theta_n, -\theta_{n+1}}_{\mathrm{free}}(\mathcal{F}_n)\leq \mathbb{P}\Big(\min_{t\in [\zeta_n,\zeta_{n+1}]}B(t)\leq -\{\theta_n\vee \theta_{n+1}\}-\frac{\epsilon}{2}s - \frac{\nu}{4}\zeta^2_n\Big).
\end{align}
Combining this with \eqref{eq:BB}, \eqref{eq:Monotone11} and \eqref{eq:PCE} yields
\begin{align}
\mathbb{P}\big((\neg \mathcal{E}_{n})\cap(\neg\mathcal{E}_{n+1})\cap \mathcal{F}_n\big) \leq \mathbb{P}\Big(\min_{t\in [\zeta_n,\zeta_{n+1}]}B(t)\leq -\{\theta_n\vee \theta_{n+1}\}-\frac{\epsilon}{2}s - \frac{\nu}{4}\zeta^2_n\Big) \, . \label{eq:BB&Mon}
\end{align}
Employing an elementary Brownian bridge estimate \cite[Lemma 2.5]{CG18b} shows that
\begin{align}
\mathbb{P}\big((\neg \mathcal{E}_{n})\cap(\neg\mathcal{E}_{n+1})\cap \mathcal{F}_n\big) \leq  \exp\Big(-s\big(\frac{\nu}{4}\zeta^2_n + \frac{\epsilon s}{2}\big)^2\Big)\leq \exp\Big(-\frac{s\nu^2}{16}\zeta^4_n- \frac{\epsilon^2s^3}{2}\Big) \, ,
\end{align}
where the second inequality is due to $(a+b)^2\geq a^2+b^2$ for any $a,b>0$. Summing both sides of the above inequality over $n\in \ZZ$ and bounding the right-hand sum in the manner of \eqref{eq:SumIntConv}, we arrive at~\eqref{eq:EESum} and thus complete the proof of Proposition~\ref{ppn:LowerTail}.
\end{proof}


\bp\label{ppn:UpperTail}
For any $t_0>0$ and $\nu\in (0,1)$, there exist $s_0= s_0(t_0,\nu)>0$ and $c_1=c_1(t_0,\nu)>c_2=c_2(t_0,\nu)>0$ such that, for $t\geq t_0$ and $s>s_0$,
\begin{align}
\exp\big(- c_1s^{3/2}\big) \leq \mathbb{P}(\mathsf{A}) \leq \exp\big(- c_2 s^{3/2}\big)\qquad \textrm{where}\quad \mathsf{A}= \Big\{\sup_{x\in \RR}\Big( \mathfrak{h}_t(x)+\frac{(1-\nu)x^2}{2}\Big) \geq s\Big\}.\qquad\label{eq:SupU}
\end{align}
\ep

\begin{proof}

The first inequality of \eqref{eq:SupU} follows from Proposition \ref{ppn:onepointuppertail} since  $\{\mathfrak{h}_t(0)\geq s\}\subseteq \big\{\sup\limits_{x\in \RR}\mathfrak{h}_t(x)\geq s \big\}$.

Turning to prove the second inequality  of \eqref{eq:SupU}, for $n\in \ZZ$, set $\zeta_n = n/s$ and
\begin{align}
\mathsf{E}_n &:= \Big\{\mathfrak{h}_t(\zeta_n)\geq  -\frac{(1-\tfrac{\nu}{2})}{2}\zeta^2_n + \frac{s}{2}\Big\}\\
\mathsf{F}_n &:= \Big\{\mathfrak{h}_t(x)\geq -\frac{(1-\nu)}{2}y^2+ s, \text{ for some }y \in (\zeta_n, \zeta_{n+1})  \Big\}.
\end{align}
(These events are similar to those in the derivation of Proposition \ref{ppn:LowerTail}.  Since the present events concern the upper tail, and their earlier cousins the lower tail,
we replace calligraphic by sans-serif font to denote them.)
Seeking the second inequality  of \eqref{eq:SupU}, note that
\begin{align}
\mathbb{P}(\mathsf{A})\leq \sum_{n\in \ZZ}\mathbb{P}(\mathsf{E}_n) + \sum_{n\in \ZZ}\mathbb{P}\big( \neg \mathsf{E}_n \cap \neg \mathsf{E}_{n+1} \cap \mathsf{F}_n\big) + \mathbb{P}\Big(\mathsf{A}\cap \Big\{\bigcap_{n\in \ZZ} \neg \mathsf{E}_n \Big\}\cap \Big\{\bigcap_{n\in \ZZ} \neg \mathsf{F}_n \Big\}\Big) \, .\qquad \label{eq:ProbSplit12}
\end{align}
Note that $\mathsf{A}\cap \neg \mathsf{F}_n = \emptyset$, so that  the last term on the right-hand side of \eqref{eq:ProbSplit12} equals zero.

For the first term on the right-hand side of \eqref{eq:ProbSplit12}, the second inequality in Proposition~\ref{ppn:UpperTail} shows that
\begin{align}
\sum_{n\in \ZZ}\mathbb{P}(\mathsf{E}_n) &\leq \sum_{n\in \ZZ} \exp\Big(-c\big(\tfrac{1-2^{-1}\nu}{2}\zeta^2_n+s\big)^{3/2}\Big)\leq \exp\big(-c s^{3/2}\big) \, .\label{eq:E_nBd}
\end{align}
The second inequality here follows from the argument used in \eqref{eq:SeqProbBd} and \eqref{eq:SumIntConv}.%
%

Thus, the proof of Proposition~\ref{ppn:UpperTail} will be finished if we can  show that
\begin{align}\label{eq:EEFseries}
 \sum_{n\in \ZZ}\mathbb{P}\big( \neg \mathsf{E}_n \cap \neg \mathsf{E}_{n+1} \cap \mathsf{F}_n\big) \leq \exp\big(-c s^{3/2}\big) \, .
\end{align}
For $\widetilde{\mathsf{E}}_n := \Big\{ \mathfrak{h}_t(\zeta_n) \geq -\frac{(1+2^{-1}\nu)}{2}\zeta^2_n- s^{2/3}\Big\}$, we have that
\begin{align*}
\mathbb{P}\big( \neg \mathsf{E}_{n-1} \cap \neg \mathsf{E}_{n+1} \cap \mathsf{F}_n\big)&\leq \mathbb{P}(\mathsf{B}_n) + \mathbb{P}( \neg\widetilde{\mathsf{E}}_{n-1})+ \mathbb{P}(\neg \widetilde{\mathsf{E}}_{n+1}) \, ,
\end{align*}
where we have defined the event
$\mathsf{B}_n = \neg \mathsf{E}_{n-1} \cap \widetilde{\mathsf{E}}_{n-1} \cap  \neg \widetilde{\mathsf{E}}_{n+1} \cap \widetilde{\mathsf{E}}_{n+1} \cap \widetilde{\mathsf{F}}_n$.
We will bound each term on the last displayed right-hand side.
Proposition~\ref{ppn:onepointlowertail} implies that 
   \[\mathbb{P}(\neg\widetilde{\mathsf{E}}_n)\leq \exp\Big(- c\big(s^{2/3} +\tfrac{\nu}{4} \zeta^2_n\big)^{5/2}\Big) \, . \]
   Summing over $n\in \ZZ$ in the same way as in \eqref{eq:E_nBd} yields
    \begin{align}\label{eq:SumRes}
    \sum_{n\in \ZZ}  \big(\mathbb{P}\big( \neg\widetilde{\mathsf{E}}_{n-1})+ \mathbb{P}(\neg\widetilde{\mathsf{E}}_{n+1})\big) \leq \exp\big(- cs^{5/2}\big) \, .
    \end{align}
It remains to show that $\sum_{n\in \ZZ}\mathbb{P} (\mathsf{B}_n) \leq \exp\big(-c s^{3/2}\big)$.
This readily follows once we show
    \begin{align}\label{eq:CplexBd}
   \mathbb{P} (\mathsf{B}_n)\leq 2\mathbb{P}\,\Big(\mathfrak{h}_t(0)\geq \tfrac{\nu}{16}\zeta^2_n+\tfrac{1}{8}s\Big) \, .
    \end{align}
Indeed, the latter right-hand side can be bounded above by appealing to the right-hand inequality in Proposition \ref{ppn:onepointuppertail}; the bound on the sum of $\mathbb{P} (\mathsf{B}_n)$ follows then by the logic that governs  \eqref{eq:E_nBd}.

\begin{figure}[t]
\includegraphics[width=.5\linewidth]{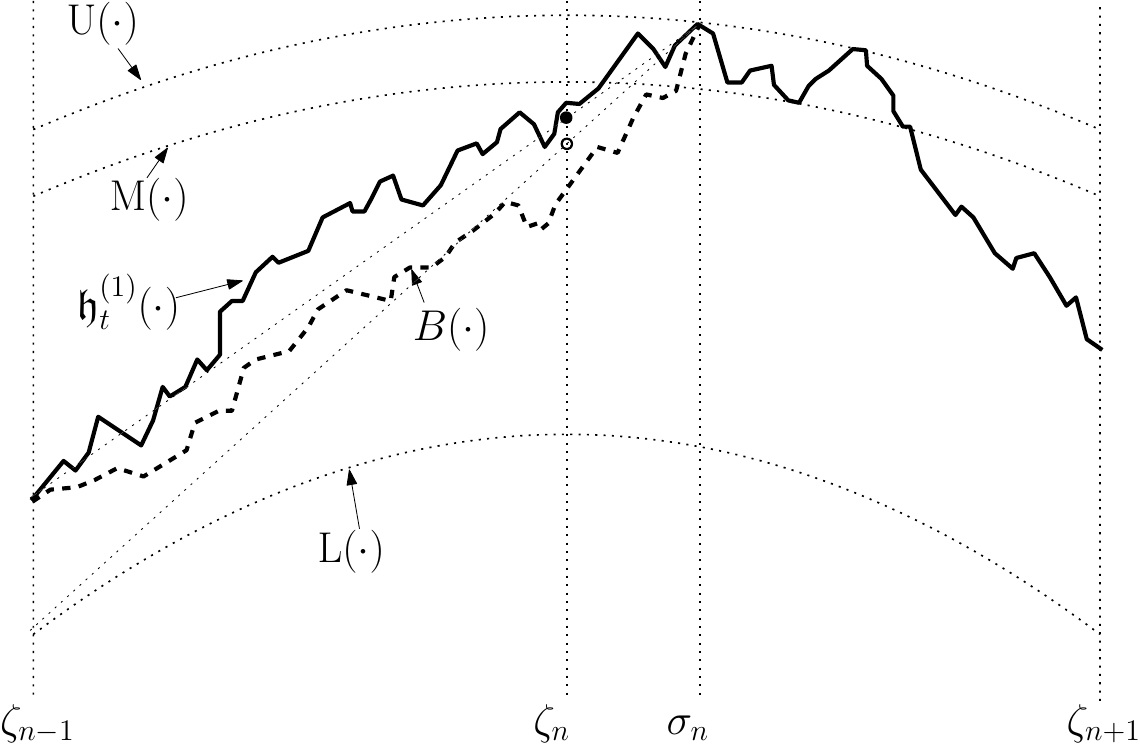}
\caption{Illustration for the proof of (\ref{eq:CplexBd}).}
\label{fig:Figure2}
\end{figure}

The remainder of this proof is devoted to showing \eqref{eq:CplexBd}. We will rely upon the equality in distribution  between  the narrow wedge solution $\mathfrak{h}_t(\cdot)$ and  the lowest labelled curve $\mathfrak{h}^{(1)}_t(\cdot)$ of the KPZ line ensemble that is offered by Proposition \ref{NWtoLineEnsemble}. Consequently, in the definitions of the events  $\mathsf{E}_n$, $\mathsf{F}_n$, $\widetilde{\mathsf{E}}_n$ and $\mathsf{B}_n$, we may substitute $\mathfrak{h}_t$ with $\mathfrak{h}^{(1)}_t$. Our proof of \eqref{eq:CplexBd} parallels the proof of \cite[Proposition~4.4]{CH14}; see also \cite[Lemma~4.1]{CorHam16}.

Define three curves (and consult Figure~\ref{fig:Figure2} for an illustration of the main objects in this proof):
   \begin{align}
   U(y):= -\frac{(1-\nu)}{2}y^2 +s, \quad L(y):= -\frac{(1+\tfrac{\nu}{2})}{2}y^2 -s^{2/3}, \quad M(y):=  -\frac{(1-\tfrac{\nu}{2})}{2}y^2+\frac{s}{2} \, .
\end{align}
If $\neg \mathsf{E}_{n-1} \cap \widetilde{\mathsf{E}}_{n-1}$ and $\neg \mathsf{E}_{n-1} \cap \widetilde{\mathsf{E}}_{n-1}$ occurs, then $\mathfrak{h}^{(1)}_t(\cdot)$ stays in between the curves $M$ and $L$ at the points $\zeta_{n-1}$ and $\zeta_{n+1}$ respectively. If $\mathsf{F}_n$ occurs, then $\mathfrak{h}^{(1)}_t$ touches the curve $U$ at some point in the interval $[\zeta_{n}, \zeta_{n+1}]$. Therefore, on the event $\mathsf{B}_n$, the curve $\mathfrak{h}^{(1)}_t$ hits $U$ somewhere in the interval $(\zeta_n, \zeta_{n+1})$, whereas it stays in between $M$ and $L$ at the points $\zeta_{n-1}$ and $\zeta_{n+1}$. The solid black curve in Figure \ref{fig:Figure2} is an instance of  $\mathfrak{h}_t(\cdot)$  on the event $\mathsf{B}_n$.

Let us define $\sigma_n := \inf \Big\{y\in (\zeta_n, \zeta_{n+1}): \mathfrak{h}^{(1)}_t(y)\geq U(y) \Big\}.$
Consider the crossing event
   \begin{align}\label{eq:DefMathfrakB}
    \mathsf{C}_n := \Big\{\mathfrak{h}^{(1)}_t(\zeta_n)\geq \frac{\sigma_n-\zeta_n}{\sigma_n -\zeta_{n-1}}L(\zeta_{n-1})+ \frac{\zeta_n-\zeta_{n-1}}{\zeta_n- \zeta_{n-1}}U(\sigma_n)\Big\} \, ,
   \end{align}
 which in Figure \ref{fig:Figure2} is the event that the solid black curve stays above the solid bullet at time $\zeta_n$.
   We adopt the shorthand
   \begin{align}
\mathbb{P}_{\mathbf{H}_{t}}:=\mathbb{P}^{1,1,(\zeta_{n-1}, \zeta_n), \mathfrak{h}^{(1)}_t(\zeta_{n-1}), \mathfrak{h}^{(1)}_t(\sigma_{n}), +\infty, \mathfrak{h}^{(2)}_t}_{\mathbf{H}_{t}}, \quad \widetilde{\mathbb{P}}_{\mathbf{H}_{t}}:=\mathbb{P}^{(\zeta_{n-1}, \sigma_n), \mathfrak{h}^{(1)}_t(\zeta_{n-1}), \mathfrak{h}^{(1)}_t(\sigma_{n})}_{\mathrm{free}}.\label{eq:DefP2Up}
\end{align}

Recalling Definition~\ref{LineEns}, note that,   since  $(\zeta_{n-1}, \sigma_n)$ is a $\{1\}$-stopping domain   for the KPZ line ensemble, the strong $\mathbf{H}_{t}$-Brownian Gibbs property Lemma \ref{lem:strongBGP} implies that
 \begin{align}\label{eq:TildeB}
 \mathbb{E}\big[\mathbf{1}_{\mathsf{B}_n} \,\mathbf{1}_{\mathsf{C_n}}|\mathfrak{F}_{\mathrm{ext}}\big(\{1\}\times (\zeta_{n-1}, \sigma_{n})\big)\big] = \mathbf{1}_{\mathsf{B}_n}\,  \mathbb{P}_{\mathbf{H}_{t}}(\mathsf{C}_n) \, .
\end{align}

By Lemma~\ref{Coupling1}, there exists a monotone coupling 
between the laws $\mathbb{P}_{\mathbf{H}_{t}}$ and $\widetilde{\mathbb{P}}_{\mathbf{H}_{t}}$; in Figure~\ref{fig:Figure2}, the curve $B$ is supposed to represent a sample from $\widetilde{\mathbb{P}}_{\mathbf{H}_{t}}$ coupled to $\mathfrak{h}^{(1)}_t$ distributed according to $\mathbb{P}_{\mathbf{H}_{t}}$.
 Using this and that the probability of $\mathsf{C}_n$ increases under pointwise increase of its sample paths, we find that $\mathbb{P}_{\mathbf{H}_{t}}(\mathsf{C}_n)\geq \widetilde{\mathbb{P}}_{\mathbf{H}_{t}}(\mathsf{C}_n)$. Since $\widetilde{\mathbb{P}}_{\mathbf{H}_{t}}$ is the law of a Brownian bridge, 
 there is probability one-half that it stays above the line joining the two endpoints at a given intermediate time such as~$\zeta_n$. Since that linear interpolation value at time $\zeta_n$ (the empty circle in Figure \ref{fig:Figure2}) sits below the value considered in $\mathsf{C}_n$ (the solid bullet in Figure \ref{fig:Figure2}), we see that $\widetilde{\mathbb{P}}_{\mathbf{H}_{t}}(\mathsf{C}_n)\geq 1/2$. Substituting this into \eqref{eq:TildeB} and taking expectation yields
    \begin{align}
\mathbb{P} (\mathsf{B}_n)\leq 2\,\mathbb{E}\big[\mathbf{1}_{\mathsf{B}_n}  \mathbf{1}_{\mathsf{C}_n}\big] \, . \label{eq:AllIsAbove}
\end{align}

  To bound the right-hand side of \eqref{eq:AllIsAbove}, observe that
  \begin{align}\label{eq:Ubounds}
  \frac{\sigma_n-\zeta_n}{\sigma_n -\zeta_{n-1}}L(\zeta_{n-1})+ \frac{\zeta_n-\zeta_{n-1}}{\sigma_n- \zeta_{n-1}}U(\sigma_n)\geq -\frac{(1-\tfrac{\nu}{8})}{2}\zeta^2_n -\frac{(4+34\nu)}{2s^{2}}+ \frac{1}{2}\Big(\frac{s}{2} - s^{2/3}\Big).
  \end{align}
 In order to demonstrate \eqref{eq:Ubounds}, we use the equalities and bounds:
\begin{align}
  \frac{(\sigma_n -\zeta_{n})\zeta^2_{n-1}+(\zeta_n -\zeta_{n-1})\sigma^2_{n}}{\sigma_n - \zeta_{n-1}}- \zeta^2_n &= (\sigma_n -\zeta_{n})(\zeta_n -\zeta_{n-1})\leq \frac{1}{s^{2}} \, , &&&\label{eq:Series1}\\
  \frac{-\frac{1}{2}(\sigma_n -\zeta_{n})\zeta^2_{n-1}+(\zeta_n -\zeta_{n-1})\sigma^2_{n}}{\sigma_n - \zeta_{n-1}}  +\frac{1}{2}\zeta^2_n &= -\frac{1}{2}(\sigma_n -\zeta_{n})(\zeta_n -\zeta_{n-1}) + \frac{3}{2}\frac{\zeta_n -\zeta_{n-1}}{\sigma_n -\zeta_{n-1}}\sigma^2_{n} \, , &&&\label{eq:Series2}\\
  \frac{3}{2}\frac{\zeta_n -\zeta_{n-1}}{\sigma_n -\zeta_{n-1}}\sigma^2_{n}-\frac{1}{2}\zeta^2_n &\geq \frac{1}{4}\zeta^2_n - 2\frac{|\zeta_n|}{s^{1}} \geq \frac{1}{8}\zeta^2_n - \frac{32}{s^{2}} \, .&&&\label{eq:Series3}
\end{align}
The first inequality of \eqref{eq:Series3} uses that $(\zeta_n-\zeta_{n-1})/(\sigma_n-\zeta_{n-1})\geq \frac{1}{2}$ and $\sigma^2_n\geq \zeta^2_n -2|\zeta_n|s^{-(1+\delta)}$; the second inequality of that line follows from $8^{-1}\zeta^2_n- 2|\zeta_n|s^{-(1+\delta)}\geq 0$ for all $|n|\geq 16$.

%

Owing to \eqref{eq:Ubounds}, when $\mathsf{C}_n$ occurs, $\mathfrak{h}^{(1)}_t(\zeta_{n}) $ will be greater than the right-hand side of \eqref{eq:Ubounds}. The latter quantity is bounded below by $-2^{-2/3}(1-\nu/8)+s/8$ when $s$ is large enough. Hence, omitting the indicator $\mathbf{1}_{\mathsf{B}_n}$ from the right-hand side of \eqref{eq:AllIsAbove}, we learn that
   \begin{align}
   \mathbb{P} (\mathsf{B}_n)\leq 2\,\mathbb{P}(\mathsf{C}_n)\leq 2 \,\mathbb{P}\Big(\mathfrak{h}^{(1)}_t(\zeta_n)\geq -\frac{(1-\tfrac{\nu}{8})}{2}\zeta^2_n + \frac{s}{8}\Big) \, .\label{eq:EachBd}
   \end{align}
  The claim \eqref{eq:CplexBd} now follows by recalling from Proposition \ref{ppn:Stationarity} that $\mathfrak{h}^{(1)}_t(\zeta^2_n)+\frac{\zeta^2_n}{2}\stackrel{d}{=}\mathfrak{h}^{(1)}_t(0)$.
   \end{proof}

Whereas Propositions \ref{ppn:LowerTail} and \ref{ppn:UpperTail} deal with the lower and upper tail of the infimum and supremum of the entire spatial process $\mathfrak{h}_t$, Proposition \ref{p.locreg} addresses the tail behaviour of small spatial increments of $\mathfrak{h}_t$. This proposition asserts that conditioned on a good -- or typical -- event~\eqref{eq:GEvent}, the tails of the increments are roughly the same as for that of Brownian motion; the result's proof is a brief but powerful application of the Brownian Gibbs technique which runs in parallel to the derivation of its zero-temperature cousin~\cite[Proposition~$3.5$]{Hammond18b}. The good event with which Proposition~\ref{p.locreg} deals has lighter than Gaussian tails, so that, without conditioning, the power law in the exponential becomes $3/2$ instead of $2$. The result that thus arises was recorded earlier as Theorem~\ref{t.locreg}, and is proved a little later in this section, along with another consequence of Proposition \ref{p.locreg} -- namely, Proposition \ref{cor:HitParabola}, concerning tails of spatial increments.


\bp\label{p.locreg}
For any $x\in \RR$, $\e\in(0,1]$ and $s\geq 0$, define
\begin{align}\label{eq:GEvent}
\boundgood_{\e,s}(x) =\!\!\! \!\!\!\bigcap_{y \in \{ x+\e-2,x,x+\e, x+2 \}}\!\!\!\!\!\! \big\{ -s/4 \leq  \mathfrak{h}_t(y) + \tfrac{y^2}{2} \leq s/4  \, \big\} \, .
\end{align}
There exist constants $c_1,c_2>0$ such that, for $\e \in (0,1]$, $t>0$, $s \geq 4$ and $x \in \R$, 
\begin{align}\label{eq:hBound}
 \PP \, \bigg( \, \sup_{z \in [x,x+\e]} \Big\vert \big(\mathfrak{h}_{t}(z)+\tfrac{z^2}{2}\big) - \big(\mathfrak{h}_{t}(x)+\tfrac{x^2}{2}\big)  \Big\vert \, \geq \,  s\e^{1/2} \, , \, \boundgood_{\e, s}(x) \, \bigg)\leq c_1\exp( - c_2 s^2) \, .
\quad \end{align}


\ep


\begin{proof}
By the stationarity in $x$ of $\mathfrak{h}_{t}(x)+\tfrac{x^2}{2}$,  we may suppose that $x =0$. 
By the equality in distribution  of $\mathfrak{h}_t(\cdot)$ with $\mathfrak{h}^{(1)}_t(\cdot)$ stated in Proposition \ref{NWtoLineEnsemble},
we may substitute $\mathfrak{h}^{(1)}_t$ for $\mathfrak{h}_t$ in the statement and proof of the desired result. Now, define the events
\begin{align*}
\down_{\e,s}&= \Big\{  \inf_{z \in [0,\e]}\big\{\mathfrak{h}^{(1)}_{t}(z)+\tfrac{z^2}{2}\big\}  \leq \mathfrak{h}^{(1)}_{t}(0) - s \e^{1/2} \Big\},\\
\up_{\e,s} &= \Big\{  \sup_{z \in [0,\e]}\big\{ \mathfrak{h}^{(1)}_{t}(z) +\tfrac{z^2}{2} \big\}\geq \mathfrak{h}^{(1)}_{t}(0) + s \e^{1/2} \Big\} \, .
\end{align*}
To obtain Proposition~\ref{p.locreg} it is enough to verify two bounds (here and below, let $\boundgood_{\e, s}=\boundgood_{\e, s}(0)$):
\begin{align}\label{e.downboundtext}
\PP \Big(  \down_{\e,s},  \boundgood_{\e, s} \Big)
\leq  c_1\exp (- c_2 s^2),
\qquad
\PP \Big(  \up_{\e,s}  ,  \boundgood_{\e, s}\Big)
\leq  c_1\exp(- c_2 s^2).
\end{align}
We start by proving the $\down$ bound in \eqref{e.downboundtext}. By the Brownian Gibbs property and stochastic monotonicity (Lemma \ref{Coupling1}) for the KPZ line ensemble, we bound above the probability of a large fall by the corresponding probability for a suitable Brownian bridge. On the event $\boundgood_{\e, s}$ we can control this Brownian fall event. Figure~\ref{fig:Figure3} shows an illustration of this idea.

\begin{figure}[t]
\includegraphics[width=.5\linewidth]{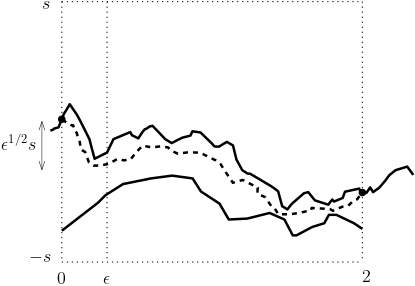}
\caption{Illustration for the proof of the $\down$ event. The top and bottom solid black curves are $\mathfrak{h}^{(1)}_t(\cdot)$ and $\mathfrak{h}^{(2)}_t(\cdot)$ respectively.
The event that $\mathfrak{h}^{(1)}_t(y)-\mathfrak{h}^{(1)}_t(0)$ is less than $-\e^{1/2}s$ for some $y\in (0,\e)$ is $\down_{\epsilon,s}$, and $\mathbb{P}^{1,2,(0,2),(\mathfrak{h}^{(1)}_{t}(0), \mathfrak{h}^{(1)}_{t}(2)),(+\infty, \mathfrak{h}^{(2)}_{t})}_{\mathbf{H}_{t}}\big(\down_{\epsilon,s}\big)$ is the probability of $\down_{\epsilon,s}$ conditioned on the sigma algebra $\mathfrak{F}_{\mathrm{ext}}\big(\{1\}\times (0, 2)\big)$ which is generated by everything outside of the first curve on the interval $(0,2)$.  The dashed black curve is a free Brownian bridge $\tilde{B}$ with law $\mathbb{P}^{(0,2),(\mathfrak{h}^{(1)}_t(0),\mathfrak{h}^{(2)}_t(2))}_{\mathrm{free}}$ coupled to $\mathfrak{h}^{(1)}_t(\cdot)$ so that it stays below $\mathfrak{h}^{(1)}_t(\cdot)$  on the interval $(0,2)$. Owing to this coupling, $\mathbb{P}^{1,2,(0,2),(\mathfrak{h}^{(1)}_{t}(0), \mathfrak{h}^{(1)}_{t}(2)),(+\infty, \mathfrak{h}^{(2)}_{t})}_{\mathbf{H}_{t}}(\down_{\e,s})\leq\mathbb{P}^{(0,2),\mathfrak{h}^{(1)}_t(0),\mathfrak{h}^{(1)}_t(2)}_{\mathrm{free}}(\widetilde{\down}_{\e,s})$ where $\widetilde{\down}_{\e,s}$ is the $\down$ event with $\tilde{B}$ in place of $\mathfrak{h}^{(1)}_t$. The probability of $\widetilde{\down}_{\e,s}$ under $\mathbb{P}^{(0,2),\mathfrak{h}^{(1)}_t(0),\mathfrak{h}^{(1)}_t(2)}_{\mathrm{free}}$ conditioned on $\mathfrak{h}^{(1)}_t(0)\leq s/4$ and $\mathfrak{h}^{(1)}_t(2)+2\geq -s/4$ is maximized when $\mathfrak{h}^{(1)}_t(0)= s/4$, and $\mathfrak{h}^{(1)}_t(2)+2= -s/4$ and the maximum value is equal to the probability of $\mathrm{inf}_{z\in[0,\e]}\{\tilde{B}(z)+z^2/2\}$ being less than $-s\e^{1/2}$ where $\tilde{B}$ is a Brownian bridge in $[0,2]$ with $\tilde{B}(0)=0$ and $\tilde{B}(2)=-2-2^{-1}s$. The reflection principle bounds this probability.
}
\label{fig:Figure3}
\end{figure}

Recalling the notation from Definition \ref{LineEns}, we will argue that
\begin{align}\label{e:longstring}
\PP \Big(  \down_{\e,s},  \boundgood_{\e, s} \Big) & \leq  \PP \Big( \down_{\e,s}  \, , \,
  \mathfrak{h}^{(1)}_{t}(0) \leq s/4   \, , \, \mathfrak{h}^{(1)}_{t}(2) + 2 \geq - s/4 \, \Big) \\
  &= \EE\Big[\mathbf{1}_{\mathfrak{h}^{(1)}_{t}(0) \leq s/4}\cdot \mathbf{1}_{\mathfrak{h}^{(1)}_{t}(2) + 2 \geq - s/4}\cdot \mathbb{P}^{1,1,(0,2), (\mathfrak{h}^{(1)}_t(0), \mathfrak{h}^{(1)}_t(2)), (+\infty, \mathfrak{h}^{(2)}_t)}_{\mathbf{H}_t}\big(\down_{\e,s}\big) \Big]\\
  &\leq  \EE\Big[\mathbf{1}_{\mathfrak{h}^{(1)}_{t}(0) \leq s/4}\cdot \mathbf{1}_{\mathfrak{h}^{(1)}_{t}(2) + 2 \geq - s/4}\cdot \mathbb{P}^{(0,2), \mathfrak{h}^{(1)}_t(0), \mathfrak{h}^{(1)}_t(2)}_{\mathrm{free}}\big( \widetilde{\down}_{\e,s}\big) \Big]\\
  &\leq \sup\Big\{ \mathbb{P}^{(0,2), y, z}_{\mathrm{free}}\big( \widetilde{\down}_{\e,s}\big): y\leq s/4 , z+2\geq -s/4\big\} \\
  &=  \mathbb{P}^{(0,2), s/4, -2-4^{-1}s}_{\mathrm{free}}\big( \widetilde{\down}_{\e,s}\big) \, ,
\end{align}
where we allow $\widetilde{\down}_{\e,s}$ to denote the $\down$ event with respect to the concerned Brownian bridge law.
Here, the first line follows by dropping two of the conditions in $ \boundgood_{\e, s}$. The equality with the second line uses the Brownian Gibbs property for the KPZ line ensemble and conditional expectations with respect to the $\sigma$-field external to $\mathfrak{h}^{(1)}_t$ on the interval $[0,2]$. The inequality with the  third line uses the monotone coupling (Lemma \ref{Coupling1}) to  remove the second curve in the conditioning. Under this coupling, the curve $\mathfrak{h}^{(1)}_t$ is bounded below by a Brownian bridge $\tilde{B}$ connecting the same endpoints. The inequality follows since the probability of the event $\down$ may only increase in response to the lowering of the curve. The inequality with the fourth line comes from maximizing over all values of $\mathfrak{h}^{(1)}_t(0)$ and $\mathfrak{h}^{(1)}_t(2)$ that satisfy the pair of conditions in the indicator functions of the preceding line. This maximum is achieved when the boundaries are maximally displaced; thus, the final equality.

We may rewrite the final term in \eqref{e:longstring} as
$\mathbb{P}^{(0,2), s/4, -2- s/4}_{\mathrm{free}}\big(\widetilde{\down}_{\e,s}\big)= \mathbb{P}\Big(\inf\limits_{z\in[0,\e]}\big\{ B^\prime(z) +\tfrac{z^2}{2} \big\}\leq - s\e^{1/2}\Big)$
where $B^\prime:[0,2]\to \RR$ is a Brownian bridge with $B^\prime(0)=0$ and $B^\prime(2)= -2-\tfrac{s}{2}$; here, we shifted the whole system down by $s/4$ to relocate the starting value to zero. Removing the parabola from the infimum only increases the probability. Now set $B(z) = B^\prime(z) -\tfrac{z}{2} (-2-\tfrac{s}{2})$, so that $B$ is a Brownian bridge with $B(0)=B(2)=0$. Taking into account the maximal effect of this linear shift shows that
\begin{align*}
\mathbb{P}\Big(\!\inf_{z\in[0,\e]} \!\big( B^\prime(z) +\tfrac{z^2}{2}\Big) \leq - s\e^{1/2}\big) &\leq \mathbb{P}\Big(\!\inf_{z\in[0,\e]}\! B(z)  \leq - s\e^{1/2}+\tfrac{\e}{2} (-2-\tfrac{s}{2})\Big)\\&\leq \mathbb{P}\big(\!\inf_{z\in[0,\e]} \!B(z)  \leq - \tfrac{s}{2}\e^{1/2} \big).
\end{align*}
The latter inequality is due to $\tfrac{\e}{2}(2+\tfrac{s}{2}) \leq s\e^{1/2}/2$, a bound that holds for $s\geq 4$ -- assuming $\e\in (0,1]$, as we do. The right-hand probability can be estimated via the reflection principle, which yields the desired bound of the form  $c_1\exp (- c_2 s^2)$ for  suitable constants $c_1, c_2>0$.

\begin{figure}[t]
\includegraphics[width=.6\linewidth]{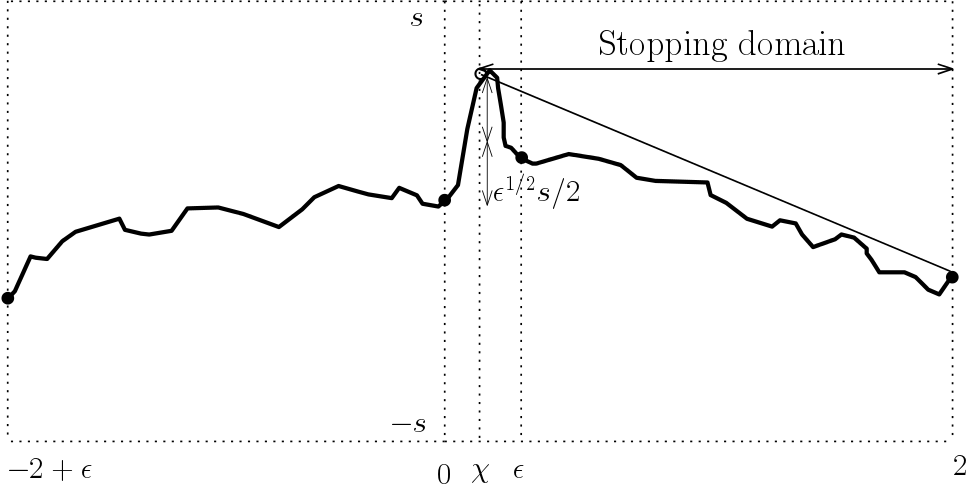}
\caption{Illustration for the proof concerning the $\up$ event. The solid black curve is $\mathfrak{h}^{(1)}_t(\cdot)$. We start by rewriting the event  $\up_{\epsilon,s}$ as $\big\{\chi \leq \e\big\}$ where $\chi$ is the smallest $x\in (0,\e]$ such that $\mathfrak{h}^{(1)}_t(x)-\mathfrak{h}^{(1)}_t(0)\geq \e^{1/2}s$; if no such $x$ satisfies this inequality, then $\chi=+\infty$. On the event $\chi\leq \e$, the interval $(\chi,2)$ forms a stopping domain, thus by the strong Gibbs property the law of $\mathfrak{h}^{(1)}_t(\cdot)$ on $(\chi,2)$ is that of a Brownian bridge conditioned to stay above $\mathfrak{h}^{(2)}_t$ on that interval. We may couple $\mathfrak{h}^{(1)}_t(\cdot)$ to lie above a free Brownian bridge. At any time in $(\chi,2)$, a Brownian bridge has probability one-half of being above the linear interpolation between its endpoints (see the diagonal line in the figure). On the event $\boundgood_{\e, s}$, we can control the slope of such a linear interpolation by a constant $c$ time $s$. Thus, between $\chi$ and $\e$, the linear interpolation may have dropped at most $cs\e$ which, for $\e$ small enough is bounded by $\e^{1/2}s/2$. The means that, on the event that $\chi\leq \e$, there is at least a one-half chance of  $\mathfrak{h}^{(1)}_t(\e)-\mathfrak{h}^{(1)}_t(0)\geq \e^{1/2}s/2$. However, viewed with a backwards' glance from position~$\e$, this occurrence  is a $\down$ event; so that the just given  argument concerning $\down$ yields a Gaussian tail for the probability of the event in question.
}
\label{fig:Figure4}
\end{figure}

We now turn to prove the $\up$ bound in \eqref{e.downboundtext} (see also Figure~\ref{fig:Figure4} and its caption).
Define $\chi$ to be the infimum of $x\in (0,\e]$ such that $\mathfrak{h}^{(1)}_{t}(x)- \mathfrak{h}^{(1)}_{t}(0)\geq s\e^{1/2}$; if no such points exist, $\chi$ is set equal to~$+\infty$. Since
$$\PP(\up_{\e,s},  \boundgood_{\e, s})= \PP(\chi \leq \e, \boundgood_{\e, s}),$$
it suffices to bound the probability on the right-hand side, which we write in the form
$$
 \PP\Big(\chi \leq \e, \boundgood_{\e, s}, \mathfrak{h}^{(1)}_{t}(\chi) - \mathfrak{h}^{(1)}_{t}(\e) <s \e^{1/2}/2 \Big) \, + \, \PP\Big(\chi \leq \e, \boundgood_{\e, s}, \mathfrak{h}^{(1)}_{t}(\chi) - \mathfrak{h}^{(1)}_{t}(\e) \geq s \e^{1/2}/2 \Big) \, .
$$
In the first displayed term, the occurrence of $\big\{\chi\leq \e\big\} \cap \big\{\mathfrak{h}^{(1)}_{t}(\chi) - \mathfrak{h}^{(1)}_{t}(\e) <s \e^{1/2}/2\big\}$ entails that  $\big\{\mathfrak{h}^{(1)}_{t}(\e) - \mathfrak{h}^{(1)}_{t}(0) \geq s \e^{1/2}/2\big\}$ holds, where note that continuity of $\mathfrak{h}^{(1)}_{t}(\cdot)$ implies that $\mathfrak{h}^{(1)}_{t}(\chi)$ equals $ \mathfrak{h}^{(1)}_{t}(0) + s\e^{1/2}$. By the argument concerning the $\down$ event, we can bound the probability of this occurrence by $c_1\exp(-c_2 s^2)$ for suitable constants $c_1$ and $c_2$; this is why $\boundgood_{\e, s}$ involves control on $\mathfrak{h}^{(1)}_t(-2+\e)$ and $\mathfrak{h}^{(1)}_t(\e)$.
Since $(\chi,2)$ is a strong stopping domain, by the strong Gibbs property (Lemma \ref{lem:strongBGP}), we deduce the first  bound in the next display, in whose second line, $B$ is  $\mathbb{P}^{1,1,(\chi,2), (\mathfrak{h}^{(1)}_t(\chi), \mathfrak{h}^{(1)}_t(2)), (+\infty, \mathfrak{h}^{(2)}_t)}_{\mathbf{H}_t}$-distributed; in whose third line, $B$ is  $\mathbb{P}^{(\chi,2), \mathfrak{h}^{(1)}_t(\chi), \mathfrak{h}^{(1)}_t(2)}_{\mathrm{free}}$-distributed; and in whose fourth line, $\bar{B}$ is the linear interpolation between $\big(\chi,B(\chi)\big)$ and $\big(2,B(2)\big)$.

\begin{eqnarray*}
& & \PP\big(\chi \leq \e, \boundgood_{\e, s}, \mathfrak{h}^{(1)}_{t}(\chi) - \mathfrak{h}^{(1)}_{t}(\e) \geq s \e^{1/2}/2 \big)\\
&= & \EE\Big[\mathbf{1}_{\chi \leq \e, \boundgood_{\e, s}} \mathbb{P}^{1,1,(\chi,2), (\mathfrak{h}^{(1)}_t(\chi), \mathfrak{h}^{(1)}_t(2)), (+\infty, \mathfrak{h}^{(2)}_t)}_{\mathbf{H}_t}\big(B(\chi) - B(\e) \geq s \e^{1/2}/2 \big) \Big]\\
& \leq & \EE\Big[\mathbf{1}_{\chi \leq \e, \boundgood_{\e, s}} \mathbb{P}^{(\chi,2), \mathfrak{h}^{(1)}_t(\chi), \mathfrak{h}^{(1)}_t(2)}_{\mathrm{free}}\big(B(\chi) - B(\e) \geq s \e^{1/2}/2 \big) \Big]\\
& \leq & \EE\Big[\mathbf{1}_{\chi \leq \e, \boundgood_{\e, s}} \mathbb{P}^{(\chi,2), \mathfrak{h}^{(1)}_t(\chi), \mathfrak{h}^{(1)}_t(2)}_{\mathrm{free}}\big(B(\e) \leq \bar{B}(\e) \big) \Big]\\
& = &  \tfrac{1}{2}  \, \EE\Big[\mathbf{1}_{\chi \leq \e, \boundgood_{\e, s}}\Big] \, . \\
\end{eqnarray*}
The inequality between the second and third line comes from applying the monotone coupling  (Lemma \ref{Coupling1}). To justify the inequality between the third and fourth lines takes a bit more work. We claim that,  on the event $\big\{\chi \leq \e, \boundgood_{\e, s}\big\}$, and as long as $\e$ is small enough, the condition $B(\chi) - B(\e) \geq s \e^{1/2}/2$ entails that $B(\e) \leq \bar{B}(\e)$.  This is because the maximal slope of $\bar{B}$ is a constant multiple~$c$ of~$s$, and the maximal displacement of $\e-\chi$ is $\e$. Thus  $\bar{B}(\e)\geq B(\chi) - cs\e$. For $\e$ small enough, $cs\e\leq s\e^{1/2}/2$; thus, $\bar{B}(\e)\geq B(\chi) - s\e^{1/2}/2$. This can be rewritten as  $\bar{B}(\e)+ s\e^{1/2}/2\geq B(\chi)$, and $B(\chi) - B(\e) \geq s \e^{1/2}/2$ can be rewritten as $B(\chi)  \geq  B(\e)+ s \e^{1/2}/2$. Combining these two yields $\bar{B}(\e) \geq B(\e)$, as claimed. The final equality is simply from Brownian bridge having probability one-half of being below its linear interpolation at any given time.

Putting this all together, we have proved that
$$
\PP\big( \chi \leq \e, \boundgood_{\e, s} \big)  \leq  c_1\exp(-c_2 s^2) + \tfrac{1}{2} \PP\big( \chi \leq \e, \boundgood_{\e, s} \big) ,
$$
which implies that
$$
\PP \big( \chi \leq \e, \boundgood_{\e, s} \big)  \leq  2 c_1\exp(-c_2 s^2),
$$
as desired to prove the $\up$ bound in \eqref{e.downboundtext}.
\end{proof}

\begin{proof}[Proof of Theorem~\ref{t.locreg}]
Applying Propositions~\ref{ppn:LowerTail} and~\ref{ppn:UpperTail}, there exist $s_0=s_0(t_0)$ and $c=c(t_0)$ such that, for $x\in \RR$, $t\geq t_0$ and $s\geq s_0$,
$\mathbb{P}\big(\neg\boundgood_{\e,s}(x) \big)\leq \exp\big(-cs^{3/2}\big)$. This inference and Proposition~\ref{p.locreg} yield \eqref{eq:SuphBd}; thus is the theorem proved.
\end{proof}

We will need one more result, concerning tails of increments.



\bp\label{cor:HitParabola} 
For any $t_0>0$, $\nu>0$ and $\epsilon\in (0,1)$, there exist $s=s(t_0,\nu,\epsilon)$ and $c=c(t_0,\nu,\epsilon)$ such that, for $t\geq t_0$, $s\geq s_0$ and $\theta>1$,
\begin{align}
\mathbb{P}\Big(\sup_{x\in [0,\theta^{1/3}]}\big\{\mathfrak{h}_t( \theta^{-2/3}x) -\mathfrak{h}_t(0)-\theta^{-1/3}\tfrac{\nu x^2}{2}\big\} \geq  \theta^{-1/3}s \Big)\leq \exp(-cs^{9(1-\epsilon)/8}) \, , \label{eq:HitParabola}\\
\mathbb{P}\Big(\inf_{x\in [0,\theta^{1/3}]}\big\{\mathfrak{h}_t( \theta^{-2/3}x) -\mathfrak{h}_t(0)+\theta^{-1/3}\tfrac{\nu x^2}{2}\big\} \leq  -\theta^{-1/3}s \Big)\leq \exp(-cs^{9(1-\epsilon)/8}) \, . \label{eq:HitParabola2}
\end{align}
\ep
\begin{proof}
We prove only \eqref{eq:HitParabola}, since a very similar argument yields~\eqref{eq:HitParabola2}. Let $\mathsf{E}$ denote the event on the left-hand side of \eqref{eq:HitParabola}.
Suppose first that $\theta<2^3$. Then
\begin{align}\label{eq:SplitEv}
\mathbb{P}\big(\mathsf{E}\big)\leq \mathbb{P} \Big(\sup_{x\in [0, 2]}\big\{\mathfrak{h}_t( x) - \mathfrak{h}_t(0)\big\}\geq \frac{s}{2}\Big)\leq \mathbb{P}\Big(\sup_{x\in \RR}\mathfrak{h}_t( x)\geq \frac{s}{4}\Big) + \mathbb{P}\Big(\mathfrak{h}_t( 0)\leq -\frac{s}{4}\Big) \, .
\end{align}
The first bound is due to the supremum increasing in response to the omission of the parabola and extension of the range of $x$ to $[0,2]$. The second bound uses  $\sup_{x\in [0,2^{1/3}]}\mathfrak{h}_t(x)\leq \sup_{x\in \RR}\mathfrak{h}_t(x)$ and
  \begin{align}
  \Big\{\sup_{x\in [0, 2]}\big\{\mathfrak{h}_t( x) - \mathfrak{h}_t(0)\big\}\geq s/2 \Big\}\subset \Big\{\sup_{x\in \RR}\mathfrak{h}_t( x)\geq  s/4\Big\}\cup \Big\{ \mathfrak{h}_t(0)\leq - s/4\Big\} \, .
  \end{align}

Propositions~\ref{ppn:onepointlowertail} and~\ref{ppn:UpperTail} provide bounds on the right-hand probabilities in~\eqref{eq:SplitEv} of the form $\exp(-cs^{3/2})$, for some $c=c(t_0)>0$ when $s>s_0(t_0)$ is large enough. This proves \eqref{eq:HitParabola} for $\theta<2^3$.

Now suppose that $\theta\geq 2^3$. We may partition $[1,\theta^{1/3}]=\bigsqcup_{\ell=1}^{3K_0+1}\mathcal{I}_\ell$ where $K_0:=\lfloor \frac{1}{3} \log_2\theta\rfloor$ and $\mathcal{I}_{\ell}:= [2^{(\ell-1)/3}, 2^{\ell/3}]$ for $1\leq \ell\leq 3K_0$ and $\mathcal{I}_{3K_0+1}:= [2^{K_0}, \theta^{1/3}]$. For $1\leq \ell\leq 3K_0+1$, define
 \begin{align}
 \mathcal{Q}_{\ell}:= \bigg\{\sup_{x\in \mathcal{I}_{\ell}} \big\{\mathfrak{h}_t( \theta^{-2/3}x) -\mathfrak{h}_t(0)\big\} \geq \theta^{-1/3}\Big(s+ \frac{\nu 2^{2(\ell-1)/3}}{2}\Big) \bigg\} \, .
\end{align}
We seek to apply Proposition~\ref{p.locreg}. This proposition will bound $\mathbb{P}(\mathcal{Q}_{\ell}\cap \boundgood_\ell)$ for an event $\boundgood_\ell$ that we will define soon. Before doing that, let us observe that if, in the $\mathcal{Q}_{\ell}$ event, we replace $x\in \mathcal{I}_{\ell}$ by $x\in [0,2^{\ell/3}]$, then the probability may only rise. Similarly, if we replace $ \theta^{-1/3}\Big(s+ \frac{\nu 2^{2(\ell-1)/3}}{2}\Big)$ by a smaller number, then the probability may also only rise. In realizing this second action, we will use the weighted arithmetic-geometric mean inequality to bound (see below for notational choices)
\begin{align}
 s+ \frac{\nu 2^{2(\ell-1)/3}}{2} =w_1\cdot\frac{s}{w_1}+ w_2\cdot\frac{\nu}{2w_2} 2^{2(\ell-1)/3} &\geq  2^{2w_2(\ell-1)/3}\Big(\frac{\nu}{2w_2}\Big)^{w_2}\Big(\frac{s}{w_1}\Big)^{w_1}\\& \geq a_\ell^{(1+3\epsilon)/4} \Big(\frac{s}{3}\Big)^{(3/4)(1-\epsilon)}
 \end{align}
where $\epsilon\in (0,1)$ is arbitrary, $w_1=\tfrac{3(1-\epsilon)}{4}$, $w_2=\tfrac{(1+3\epsilon)}{4}$ and $a_\ell = \frac{\nu}{2(1+3\epsilon)} 2^{2(\ell-1)/3}$. The last inequality of the above display follows by substituting $a_{\ell}$ and noting that $s/w_1\geq s/3$.
This means that
\begin{align}
\mathbb{P}(\mathcal{Q}_{\ell}\cap \boundgood_\ell)&\leq \mathbb{P}\Big(\Big\{\sup_{x\in [0,2^{\ell/3}]}\big\{\mathfrak{h}_t( \theta^{-2/3}x) -\mathfrak{h}_t(0)\big\} \geq  \theta^{-1/3}a^{(1+3\epsilon)/4}_\ell (s/3)^{(3/4)(1-\epsilon)}\Big\}\cap \boundgood_\ell\Big)\\\label{eq:Reverse}
&= \mathbb{P}\Big(\Big\{\sup_{x\in [0,\epsilon_\ell]}\big\{\mathfrak{h}_t(x) -\mathfrak{h}_t(0)\big\} \geq s_\ell \epsilon_\ell^{1/2}\Big\}\cap \boundgood_\ell\Big),
\end{align}
where in the second line we have used the notation $\epsilon_\ell = 2^{\ell/3}\theta^{-2/3}$ and
$$
s_\ell = 2^{\tfrac{\epsilon\ell-(1+3\epsilon)}{6}}\Big(\frac{\nu}{2(1+3\epsilon)}\Big)^{(1+3\epsilon)/4} \Big(\frac{s}{3}\Big)^{(3/4)(1-\epsilon)} ,
$$
and performed a simple change over variables.

By the definition of $\mathcal{Q}_{\ell}$, $\mathsf{E}\subset \bigcup^{3K_0+1}_{\ell=1} \mathcal{Q}_{\ell}$. Choosing $\boundgood_\ell=\boundgood_{\epsilon_\ell,s_\ell}$ from~\eqref{eq:GEvent}, by the union bound,
\begin{align}
\mathbb{P}\big(\mathsf{E}\big)\leq \sum_{\ell=1}^{3K_0+1} \Big(\mathbb{P}(\mathcal{Q}_{\ell}\cap \boundgood_\ell)+ \mathbb{P}(\neg \boundgood_\ell)\Big) \, .\label{eq:SplitProbBd}
\end{align}
Propositions \ref{ppn:onepointlowertail} and \ref{ppn:onepointuppertail} furnish $s_0=s_0(t_0,\nu,\epsilon)$ and $c=c(t_0,\nu,\epsilon)$ such that, for $\ell\geq 1$ and $s>s_0$,
\begin{align}\label{eq:GBound}
\mathbb{P}(\neg \boundgood_\ell)\leq \exp\big(- cs_{\ell}^{3/2}\big) = \exp\big(-c 2^{\epsilon(\ell-3)/4} s^{\frac{9}{8}(1-\epsilon)}\big).
\end{align}
The sum over $\ell$ of such a bound produces an upper bound of the desired form $\exp(-c s^{9(1-\epsilon)/8})$. Turning to $\mathbb{P}\big(\mathcal{Q}_{\ell}\cap \boundgood_\ell \big)$, we may apply Proposition~\ref{p.locreg} with $\epsilon=\epsilon_\ell$ and $s=s_{\ell}$. The result is precisely the same sort of bound as on $\mathbb{P}(\neg \boundgood_\ell)$ above; hence, by summing over $\ell$, we again recover the sought upper bound, of the form $\exp(-c s^{9(1-\epsilon)/8})$. This completes the proof of \eqref{eq:HitParabola}.
%
%
%
\end{proof}

\section{Remote correlation: Proof of Theorem~\ref{thm:Main1}}\label{sec:remote}

The composition law from Proposition~\ref{ppn:TimeEvl} will figure prominently in this argument since it permits us to describe the two-time distribution in terms of spatial processes that we understand well. Recall that this result shows that, for $\alpha > 1$ given,  we may write  $\mathfrak{h}_t(\alpha,0)$ as the composition of independent spatial processes $\mathfrak{h}_t(1,\cdot)$ and $\mathfrak{h}_{\alpha t\downarrow t}(\cdot)$. The latter is distributed as $\mathfrak{h}_t(\alpha-1,\cdot)$. We will use a suggestive shorthand, in the spirit of $\mathfrak{h}_{\alpha t\downarrow t}$:
\begin{align}\label{eq:shorthand}
\mathfrak{h}_{0\uparrow t}(\cdot) :=\mathfrak{h}_t(1,\cdot),\qquad \mathfrak{h}_{0\uparrow \alpha t}(\cdot) := \mathfrak{h}_t(\alpha,\cdot) \, ;
\end{align}
and, when we write $\mathfrak{h}_{0\uparrow t}$, we mean $\mathfrak{h}_{0\uparrow t}(0)$; and likewise for $\mathfrak{h}_{0\uparrow \alpha t}$ and $\mathfrak{h}_{\alpha t\downarrow t}$.
The shorthand is intended to suggest that the value of $\mathfrak{h}_{0\uparrow \alpha t}(0)$ from time zero to $\alpha t$ is obtained via the integral operation $I_t$ by composing the function $\mathfrak{h}_{0\uparrow t}(\cdot)$ with $\mathfrak{h}_{\alpha t\downarrow t}(\cdot)$, the latter under the opposite direction of time. The shorthand will cause confusion if $\mathfrak{h}_{0\uparrow t}$ or $\mathfrak{h}_{0\uparrow \alpha t}$ are regarded as functions of $t$ or $\alpha$. These two parameters are, however, given: the variable `$\cdot$' in~\eqref{eq:shorthand} is spatial.

%

This section's goal is to prove the bounds \eqref{eq:Corr1} in Theorem~\ref{thm:Main1} which in our new notation are:
\begin{align}\label{eq:Corr1rewrite}
c_1 \alpha^{-1/3} \leq \mathrm{Corr}\Big(\mathfrak{h}_{0\uparrow t}, \mathfrak{h}_{0\uparrow \alpha t}\Big) \leq c_2 \alpha^{-1/3} \, .
\end{align}
The derivation of Theorem~\ref{thm:Main1}
depends on two principal results,
Propositions \ref{lem:UpTail1} and~\ref{lem:LowTail}. We state these; use them to prove the theorem; and then prove them in turn.

Define
$\mathcal{A}_{\mathrm{high}} = \big\{\mathfrak{h}_{0\uparrow \alpha t}- \mathfrak{h}_{\alpha t\downarrow t}- \mathfrak{h}_{0\uparrow t}\geq s \big\}$ and
$\mathcal{A}_{\mathrm{low}} = \big\{ \mathfrak{h}_{0\uparrow \alpha t}-  \mathfrak{h}_{\alpha t\downarrow t}- \mathfrak{h}_{0\uparrow t}\leq - t^{-1/3}s \big\}$.

\bp\label{lem:UpTail1}
For any $t_0>0$ and $\alpha_0>1$, there exist $s_0=s_0(t_0,\alpha_0)>0$ and $c=c(t_0,\alpha_0)>0$ such that, for $s\geq s_0$, $t>t_0$ and $\alpha>\alpha_0$,
\begin{align}\label{eq:SdTail}
\mathbb{P}\big(\mathcal{A}_{\mathrm{high}}\big) \leq \exp(-cs^{3/4}) \, .
\end{align}
\ep

\bp\label{lem:LowTail}
For any $t_0>0$, there exist $s_0=s_0(t_0)>0$ and $c = c(t_0)>0$ such that, for $s\geq s_0$, $t>t_0$ and $\alpha>2$,
\begin{align}\label{eq:AlowBd1}
\mathbb{P}\big(\mathcal{A}_{\mathrm{low}}\big) \leq \exp\big(- cs^{3/2}\big) \, .
\end{align}
Further, for any $t_0>0$, $\alpha_0>1$ and $\delta>0$, there exist $s_0=s_0(t_0,\alpha_0,\delta)>0$ and $c = c(t_0,\alpha_0,\delta)>0$ such that, for $s\geq s_0$, $t>t_0$ and $\alpha>\alpha_0$,
\begin{align}\label{eq:AlowBd2}
\mathbb{P}\, \Big(\mathcal{A}_{\mathrm{low}} \, \Big\vert \, \mathfrak{h}_{0\uparrow t}\geq \mathbb{E}\big[\mathfrak{h}_{0\uparrow t}\big]+\delta\Big) \, \leq \, \exp\big(- cs^{3/2}\big) \, .
\end{align}
\ep

The two results address the rarity of the events that $\mathfrak{h}_{0\uparrow \alpha t}$ significantly exceeds, or falls below,  $\mathfrak{h}_{0\uparrow t} + \mathfrak{h}_{\alpha t\downarrow t}$. Note that, in the latter case -- via the definition of $\mathcal{A}_{\mathrm{low}}$ -- deviation is measured on scale $t^{-1/3}$.
To interpret the two results and the latter choice of scaling, it may be helpful to consider the composition law and the supremum variational problem obtained in the high $t$ limit which was a focus of attention in Section~\ref{sec:KPZfixedpoint}. In the high~$t$ case, a counterpart to Proposition \ref{lem:UpTail1} would examine the tail of the {\em difference} between the variational problem's solution and the value of the pair sum associated to the choice of location zero at the intermediate time. In this $t\nearrow \infty$ limit, Proposition~\ref{lem:LowTail} would become trivial, because the difference in question can never be negative. Back in our finite~$t$ world, the softening of the variational problem leads to a degree of violation to this strict ordering.
The $t^{-1/3}$ tail that we study probes the extent of this violation.

\smallskip

\begin{proof}[Proof of Theorem~\ref{thm:Main1}]
We first show the following stronger version of the upper bound on the correlation in \eqref{eq:Corr1rewrite}: for any $t_0>0$, there exist $c_2=c_2(0)>0$ such that, for  $t>t_0$, and $\alpha>2$, $|\mathrm{Corr}(\mathfrak{h}_{0\uparrow \alpha t}, \mathfrak{h}_{0\uparrow  t})|\leq c_2 \alpha^{-1/3}$. Recall that, by definition,
\begin{align}\label{eq:CorrDef}
\mathrm{Corr}\big(\mathfrak{h}_{0\uparrow t}, \mathfrak{h}_{0\uparrow \alpha t}\big) = \frac{\mathrm{Cov}\big(\mathfrak{h}_{0\uparrow \alpha t},\mathfrak{h}_{0\uparrow t}\big)}{\sqrt{\mathrm{Var}(\mathfrak{h}_{0\uparrow \alpha t})}\sqrt{\mathrm{Var}(\mathfrak{h}_{0\uparrow t})}} \, .
\end{align}
Under the coupling provided by Proposition \ref{ppn:TimeEvl}, $\mathfrak{h}_{\alpha t\downarrow t}$ and $\mathfrak{h}_{0\uparrow t}$ are independent. Hence,
\begin{align}\label{eq:Cov}
\big|\mathrm{Cov}&\big(\mathfrak{h}_{0\uparrow \alpha t}, \mathfrak{h}_{0\uparrow  t}\big)\big| = \mathrm{Cov}\big(\mathfrak{h}_{0\uparrow \alpha t} - \mathfrak{h}_{\alpha t\downarrow t}- \mathfrak{h}_{0\uparrow t}, \mathfrak{h}_{0\uparrow t}\big)+ \mathrm{Var}\big(\mathfrak{h}_{0\uparrow t}\big) \, .
\end{align}
Applying the Cauchy-Schwarz inequality to the first term yields
\begin{align}\label{eq:CSIneq}
\big|\mathrm{Cov}\big(\mathfrak{h}_{0\uparrow \alpha t} &- \mathfrak{h}_{\alpha t\downarrow t}- \mathfrak{h}_{0\uparrow t},  \mathfrak{h}_{0\uparrow t}\big)\big| \leq \sqrt{\mathrm{Var}\big(\mathfrak{h}_{0\uparrow \alpha t}-\mathfrak{h}_{\alpha t\downarrow t} -\mathfrak{h}_{0\uparrow t}\big)\,\mathrm{Var}(\mathfrak{h}_{0\uparrow t})} \, .
\end{align}
Substituting \eqref{eq:Cov} into and applying \eqref{eq:CSIneq} to \eqref{eq:CorrDef} then leads to
\begin{align}\label{eq:Corrbd}
\big|\mathrm{Corr}&(\mathfrak{h}_{0\uparrow \alpha t}, \mathfrak{h}_{0\uparrow  t})\big|\leq \frac{\sqrt{\mathrm{Var}(\mathfrak{h}_{0\uparrow \alpha t}-\mathfrak{h}_{\alpha t\downarrow t} -\mathfrak{h}_{0\uparrow  t})}}{\sqrt{\mathrm{Var}(\mathfrak{h}_{0\uparrow \alpha t})}} + \frac{\sqrt{\mathrm{Var}(\mathfrak{h}_{0\uparrow  t})}}{\sqrt{\mathrm{Var}(\mathfrak{h}_{0\uparrow \alpha t})}} \, .
\end{align}
The tail bounds on $\mathfrak{h}_t(0)$ (or on $\mathfrak{h}_{0\uparrow t}$) from Propositions \ref{ppn:onepointlowertail} and \ref{ppn:onepointuppertail} imply via Lemma~\ref{ppn:VarBound} that there exist $c=c(t_0)>0$ and $C=C(t_0)>0$ such that
\begin{align}
c \leq \mathrm{Var}\big(\mathfrak{h}_{0\uparrow  t}\big)\leq C \, , \qquad \textrm{and}\quad c \alpha^{2/3}\leq \mathrm{Var}\big(\mathfrak{h}_{0\uparrow \alpha t}\big)\leq C \alpha^{2/3} \, . \label{eq:VarBd1}
\end{align}
Here, the second bound uses $\mathfrak{h}_{0\to \alpha t}\eqdist \mathfrak{h}_{\alpha t}(1,0) \alpha^{1/3}$.
Similarly, Propositions ~\ref{lem:UpTail1} and ~\ref{lem:LowTail} imply that
\begin{align}
\mathrm{Var}\big(\mathfrak{h}_{0\uparrow \alpha t}-\mathfrak{h}_{\alpha t\downarrow t} -\mathfrak{h}_{0\uparrow  t}\big)\leq C \, . \label{eq:VarBd2}
\end{align}
%
Substituting \eqref{eq:VarBd1} and \eqref{eq:VarBd2} into the right-hand side of \eqref{eq:Corrbd} yields
\begin{align}
\big|\mathrm{Corr}(\mathfrak{h}_{0\uparrow \alpha t}, \mathfrak{h}_{0\uparrow  t})\big|\leq c_2 \alpha^{-1/3}
\end{align}
for a constant $c_2=c_2(t_0)>0$. This is the strengthened upper bound that we have sought.
\smallskip

Now we turn to prove the  lower bound in \eqref{eq:Corr1rewrite}. Assume for now that $t>1$, though we will eventually need to impose that $t>t_0$ for $t_0$ sufficiently large.  The lower bound will arise from an appeal to Corollary~\ref{cor:CovLow} with $X :=\mathfrak{h}_{0\uparrow \alpha t}$ and $Y :=\mathfrak{h}_{0\uparrow  t}$. The other two parameters in the corollary, $C_1$ and $C_2$, will be specified after the next calculation, which will inform our choice of their value. Observe that, for $y := \mathbb{E}[\mathfrak{h}_{0\uparrow  t}]+\delta$ with $\delta>0$,
$$
 \mathbb{E}\Big[ \mathfrak{h}_{0\uparrow \alpha t}\big| \mathfrak{h}_{0\uparrow  t} \geq y	 \Big] \geq \mathbb{E}\big[\mathfrak{h}_{\alpha t\downarrow t}\big] + y+\mathbb{E}\Big[\mathfrak{h}_{0\uparrow \alpha t} - \mathfrak{h}_{\alpha t\downarrow t}-\mathfrak{h}_{0\uparrow t}\big|\mathfrak{h}_{0\uparrow  t}\geq y\Big] \, .
$$
This bound follows from $\mathbb{E}\big[\mathfrak{h}_{\alpha t \downarrow t}|\mathfrak{h}_{0 \uparrow t}\geq y\big] = \mathbb{E}[\mathfrak{h}_{\alpha t \downarrow t}]$ -- a consequence of the independence of $\mathfrak{h}_{0 \uparrow t}$ and  $\mathfrak{h}_{\alpha t \downarrow t}$ -- and the trivial $\mathbb{E}\big[\mathfrak{h}_{0\uparrow t}|\mathfrak{h}_{0 \uparrow t}\geq y\big]\geq y$. Next note that
 \begin{align}
 &\mathbb{E}\, \Big[ \, \mathfrak{h}_{0\uparrow \alpha t} - \mathfrak{h}_{\alpha t\downarrow t}-\mathfrak{h}_{0\uparrow t} \, \Big\vert \, \mathfrak{h}_{0\uparrow  t}\geq y\, \Big] \,
 \geq \, \mathbb{E}\, \Big[ \, \min\{0,\mathfrak{h}_{0\uparrow \alpha t} - \mathfrak{h}_{\alpha t\downarrow t}-\mathfrak{h}_{0\uparrow t}\, \Big\vert \, \mathfrak{h}_{0\uparrow  t}\geq y \, \Big]\\
 = & - t^{-1/3}\int^{\infty}_{0}\mathbb{P}\Big(\mathfrak{h}_{0\uparrow \alpha t} - \mathfrak{h}_{\alpha t\downarrow t}-\mathfrak{h}_{0\uparrow t}\leq -t^{-1/3}s \, \Big\vert \,\mathfrak{h}_{0\uparrow  t}\geq y\Big)ds
 \geq -t^{-1/3}c(\delta)\label{eq:LstLine}
\end{align}
 for some $c(\delta)>0$. The last inequality is due to  an application of \eqref{eq:AlowBd2} in Proposition~\ref{lem:LowTail}.

 We now apply Corollary~\ref{cor:CovLow} with $X :=\mathfrak{h}_{0\uparrow \alpha t}$, $Y :=\mathfrak{h}_{0\uparrow  t}$, $C_1:=y=\mathbb{E}[\mathfrak{h}_{0\uparrow  t}]+\delta$ and $C_2:= \delta - t^{-\tfrac{1}{3}}c(\delta)$. Notice that $C_1$ and $C_2$ both depend on the parameter $\delta>0$, which is as yet unspecified. By \eqref{eq:CovBd2},
 \begin{align}
\mathrm{Cov}\big(\mathfrak{h}_{0\uparrow \alpha t}, \mathfrak{h}_{0\uparrow  t}\big) \geq \big(\delta-t^{-1/3}c(\delta)\big) \cdot \mathbb{P}\big(\mathfrak{h}_{0\uparrow  t}\geq y \big)\cdot 
\Big(\mathbb{E}\big[\mathfrak{h}_{0\uparrow  t}\big|\mathfrak{h}_{0\uparrow  t} \geq y\big]- \mathbb{E}\big[\mathfrak{h}_{0\uparrow  t}\big|\mathfrak{h}_{0\uparrow  t} < y\big]\Big) \, , \label{eq:CovLowBd}
\end{align}
where we recall that $y=\mathbb{E}[\mathfrak{h}_{0\uparrow  t}] +\delta$.
Observe that
\begin{align}
\mathbb{E}\big[\mathfrak{h}_{0\uparrow  t}\big|\mathfrak{h}_{0\uparrow  t} \geq y\big]- \mathbb{E}\big[\mathfrak{h}_{0\uparrow  t}\big|\mathfrak{h}_{0\uparrow  t} < y\big]
= \frac{\mathbb{E}\big[\mathfrak{h}_{0\uparrow  t}\big|\mathfrak{h}_{0\uparrow  t} \geq y\big]- \mathbb{E}\big[\mathfrak{h}_{0\uparrow  t}\big]}{\mathbb{P}(\mathfrak{h}_{0\uparrow  t} < y)}
\geq \frac{\delta}{\mathbb{P}(\mathfrak{h}_{0\uparrow  t} < \mathbb{E}[\mathfrak{h}_{0\uparrow  t}] + \delta)}\geq \delta \, .
 \label{eq:ExpBd}
\end{align}
Substituting this into \eqref{eq:CovLowBd}, we arrive at
$$
\mathrm{Cov}\big(\mathfrak{h}_{0\uparrow \alpha t}, \mathfrak{h}_{0\uparrow  t}\big) \geq \big(\delta-t^{-1/3}c(\delta)\big) \cdot \mathbb{P}\big(\mathfrak{h}_{0\uparrow  t}\geq y \big)\cdot\delta \, .
$$
Fix any $\delta>0$. Observe that for $t>t_0 := \big(\tfrac{2c(\delta)}{\delta}\big)^3$, $\delta-t^{-1/3}c(\delta)\geq \delta/2$.  By the upper bound in Proposition \ref{ppn:onepointuppertail}, for $t\geq t_0$, we may bound  $\mathbb{E}[\mathfrak{h}_{0\uparrow t}]<C$ and hence $y<C+\delta$, for $C=C(t_0)>0$. Using the lower bound in Proposition \ref{ppn:onepointuppertail}, we further infer that $\mathbb{P}\big(\mathfrak{h}_{0\uparrow  t}\geq y \big)>C'$ for $C'=C'(t_0,\delta)>0$. This shows that $\mathrm{Cov}\big(\mathfrak{h}_{0\uparrow \alpha t}, \mathfrak{h}_{0\uparrow  t}\big) \geq \big(\delta-t^{-1/3}c(\delta)\big) \cdot C'\cdot\delta$. This right-hand side is a strictly positive constant which holds uniformly over $t>t_0$. Substituting this and the upper bounds of \eqref{eq:VarBd1} into the right-hand side of \eqref{eq:CorrDef} produces the desired lower bound in \eqref{eq:Corr1rewrite} on $\mathrm{Corr}( \mathfrak{h}_{0\uparrow \alpha t}, \mathfrak{h}_{0\uparrow  t})$.
\end{proof}

\subsection{Proof of Proposition~\ref{lem:UpTail1}}\label{sec:UpTail1}

For this proof and Proposition \ref{lem:LowTail}'s, we will return to writing $\mathfrak{h}_{0\uparrow t}(0)$ in place of $\mathfrak{h}_{0\uparrow t}$, because we will utilize $\mathfrak{h}_{0\uparrow t}(x)$ for various values of $x$.
Define
%
 \begin{align}
 \mathcal{E} &:=\Big\{\sup_{|x|<  t^{2/3}} \big\{\mathfrak{h}_{0\uparrow t}(t^{-2/3}x) + \tfrac{x^2}{4t^{4/3}} -\mathfrak{h}_{0\uparrow t}(0)\big\} \geq \frac{s}{2}  \Big\},\\
 \mathcal{G} &:= \Big\{\sup_{|x|\geq   t^{2/3}}\big\{\mathfrak{h}_{0\uparrow t}(t^{-2/3}x)+\tfrac{x^2}{4t^{4/3}}\big\}\geq \frac{s}{4}\Big\}\cup \Big\{\mathfrak{h}_{0\uparrow t}(0) \leq - \frac{s}{4}\Big\},\\
 \mathcal{B} &:= \Big\{ \int^{\infty}_{-\infty} e^{t^{1/3}(\mathfrak{h}_{\alpha t\downarrow t}(t^{-2/3} x)-\frac{x^2}{4t^{4/3}})} dx \geq  e^{t^{1/3}(\mathfrak{h}_{\alpha t\downarrow t}(0)+\frac{s}{2})}\Big\}.
\end{align}


The derivation has three steps. \textbf{Step I} shows that $\mathcal{A}_{\mathrm{high}} \cap \neg\mathcal{E} \cap \neg\mathcal{G} \subset \mathcal{B}$. The desired bound on $\mathbb{P}(\mathcal{A}_{\mathrm{high}})$ will thus result from bounds on  $\mathbb{P}(\mathcal{B})$ and $\mathbb{P}(\mathcal{E} \cup \mathcal{G})$ which are provided in \textbf{Steps II} and~\textbf{III}.

\smallskip
\textbf{Step I:} To show that $\mathcal{A}_{\mathrm{high}} \cap \neg\mathcal{E} \cap \neg\mathcal{G} \subset \mathcal{B}$,  we will argue that, on the event $\neg\mathcal{E}\cap \neg\mathcal{G}$,
\begin{align*}
\int_{|x|< t^{2/3}}
e^{t^{1/3}\big(\mathfrak{h}_{\alpha t\downarrow t}(t^{-2/3} x)+\mathfrak{h}_{0\uparrow t}(-t^{-2/3}x)\big)}dx &\leq
e^{t^{1/3}(\mathfrak{h}_{0\uparrow t}(0)+\frac{s}{2})}
\int_{|x|< t^{2/3}} e^{t^{1/3}\big(\mathfrak{h}_{\alpha t\downarrow t}(t^{-2/3} x)-\frac{x^2}{4t^{4/3}}\big)}  dx  \\
\int_{|x|\geq   t^{2/3}} e^{t^{1/3}\big(\mathfrak{h}_{\alpha t\downarrow t}(t^{-2/3} x)+\mathfrak{h}_{0\uparrow t}(-t^{-2/3}x)\big)} dx&\leq
e^{t^{1/3}(\mathfrak{h}_{0\uparrow t}(0)+\frac{s}{2})}
\int_{|x|\geq  t^{2/3}} e^{t^{1/3}\big(\mathfrak{h}_{\alpha t\downarrow t}(t^{-2/3} x)-\frac{x^2}{4t^{4/3}}\big)} dx \, .
\end{align*}
Indeed, the first bound holds because, on the event $\neg\mathcal{E}$, $\sup_{|x|\leq  t^{2/3}}\mathfrak{h}_{0\uparrow t}(t^{-2/3}x)\leq \mathfrak{h}_{0\uparrow t}(0)-\tfrac{x^2}{4 t^{4/3}}+ s/2$; while the second bound is due to the validity on the event $\neg\mathcal{G}$, and for  $|x|> t^{2/3}$, of the bound
$\mathfrak{h}_{0\uparrow t}(t^{-2/3}x)\leq s/4-\frac{x^2}{4t^{4/3}}\leq \mathfrak{h}_{0\uparrow t}(0)+s/2-\frac{x^2}{4t^{4/3}}$.

 Summing the displayed bounds and using the composition law Proposition \ref{ppn:TimeEvl} yields
 \begin{align}\label{eq:SplitComb}
 e^{t^{1/3}\mathfrak{h}_{0\uparrow \alpha t}(0)} \leq e^{t^{1/3}\big(\mathfrak{h}_{0\uparrow t}(0)+ s/2\big)} \int^{\infty}_{-\infty} e^{t^{1/3} \big(\mathfrak{h}^{(2)}_{\alpha t\downarrow t}(t^{-2/3} x) - \frac{x^2}{4t^{4/3}}\big)} dx \, .
 \end{align}

On the event $\mathcal{A}_{\mathrm{high}}$, the left-hand side of \eqref{eq:SplitComb} is at least $e^{t^{1/3}(\mathfrak{h}_{\alpha t\downarrow t}(0)+ \mathfrak{h}_{0\uparrow t}(0)+s)}$. Applying this, and cancelling  $e^{t^{1/3}(\mathfrak{h}_{0\uparrow t}(0)+ s/2)}$, we arrive at the inequality which defines the event $\mathcal{B}$; hence, we conclude that  $\mathcal{A}_{\mathrm{high}}\cap \neg\mathcal{E} \cap \neg\mathcal{G} \subset \mathcal{B}$, as desired.

\smallskip

\textbf{Step II:} Here, we prove that, for $t_0>0$ and $\alpha_0>1$, there exist $s_0= s_0(t_0, \alpha_0)$ and $c=c(t_0, \alpha_0)$ such that, for $t> t_0$, $\alpha>\alpha_0$ and $s>s_0$,
\begin{align}\label{eq:BcalBd}
\mathbb{P}(\mathcal{B})\leq \exp(-cs^{3/4}) \, .
\end{align}

Set $x_0 =  ((\alpha-1)t^2)^{1/3}$ and introduce two events
\begin{align}
\widetilde{\mathcal{E}} &:=\bigg\{\sup_{|x|<x_0}\big\{\mathfrak{h}_{\alpha t\downarrow t}(t^{-2/3} x) - \mathfrak{h}_{\alpha t\downarrow t}(0)-\tfrac{x^2}{8t^{4/3}}\big\}\geq s/4  \bigg\} \, ,\label{eq:tildeE1}\\
\widetilde{\mathcal{G}} &:=\bigg\{\sup_{|x|\geq x_0} \big\{\mathfrak{h}_{\alpha t\downarrow t}(t^{-2/3} x)\big\}\geq \frac{s}{8} +\frac{k^2_0}{8 t^{4/3}}\bigg\}\cup \bigg\{\mathfrak{h}_{\alpha t\downarrow t}(0)\leq - \frac{s}{8} -\frac{k^2_0}{8t^{4/3}}\bigg\} \, .\label{eq:tildeG}
\end{align}

We first show that there exists $s_0=s_0(t_0)$ such that, for $s\geq s_0$, $\mathcal{B}\subset \widetilde{\mathcal{E}}\cup \widetilde{\mathcal{G}}$. We show the contrapositive: $\neg\widetilde{\mathcal{E}} \cap \neg\widetilde{\mathcal{G}} \subset \neg\mathcal{B}$ when $s$ is large enough. On the event $\neg\widetilde{\mathcal{E}}\cap \neg\widetilde{\mathcal{G}}$, we have that
\begin{align}
\int_{|x|<x_0}  e^{t^{1/3} \big(\mathfrak{h}_{\alpha t\downarrow t}(t^{-2/3} x) - \frac{x^2}{4t^{4/3}}\big)} dx &\leq
 e^{t^{1/3}\big(\mathfrak{h}_{\alpha t\downarrow t}(0)+ s/4 \big)} \int_{|x|<x_0} e^{- \frac{x^2}{8t}} dx \, , \label{eq:Inclusion1}\\
 \int_{|x|\geq x_0}   e^{t^{1/3} \big(\mathfrak{h}_{\alpha t\downarrow t}(t^{-2/3} x) - \frac{x^2}{4t^{4/3}}\big)} dx &\leq
 e^{t^{1/3}\big(\mathfrak{h}_{\alpha t\downarrow t}(0)+ s/4 \big)} \int_{|x|\geq x_0} e^{-\frac{(|x|- x_0)^2}{4t}} dx \, . \label{eq:Inclusion3}
\end{align}
Indeed, the first bound is due to the event $\neg\widetilde{\mathcal{E}}$ entailing that $\mathfrak{h}_{\alpha t\downarrow t}(x)\leq \mathfrak{h}_{\alpha t\downarrow t}(0)+\tfrac{s}{4}+ \frac{x^2}{8t^{4/3}}$ for all $x\in [-x_0,x_0]$.  The second bound follows by combining the inequality $\exp(-(|x|^2-x_0^2)/4t)\leq \exp(-(|x|-x_0)^2/4t)$ for all $|x|\geq x_0$ with the fact that, on the event $\neg\widetilde{\mathcal{G}}$,
$$\sup_{|x|\geq x_0}t^{1/3}\mathfrak{h}_{\alpha t\downarrow t}(t^{-2/3} x)\leq t^{1/3} \big(\frac{s}{8}+\frac{k^2_0}{8 t^{4/3}}\big)\leq t^{1/3}\big(\mathfrak{h}_{\alpha t\downarrow t}(0)+\frac{s}{4}+\frac{k^2_0}{4 t^{4/3}}\big).$$

Summing \eqref{eq:Inclusion1} and \eqref{eq:Inclusion3}, and using  $\int_{-\infty}^{\infty}\big( e^{- \frac{x^2}{8t}} +  e^{-\frac{(|x|- x_0)^2}{4t}} \big)dx \leq c t^{1/2}$ for some $c>0$, yields
\begin{align}\label{eq:Inclusion}
\int^{\infty}_{-\infty} e^{t^{1/3}\big( \mathfrak{h}_{\alpha t\downarrow t}(t^{-2/3} x) - \frac{x^2}{4t^{4/3}}\big)} dx\leq c t^{1/2} e^{t^{1/3}\big(\mathfrak{h}_{\alpha t\downarrow t}(0)+ 4^{-1}s\big)}\leq e^{t^{1/3}\big(\mathfrak{h}_{\alpha t\downarrow t}(0)+ 2^{-1}s\big)} \, ,
\end{align}
where the second inequality holds provided that $s>s_0(t_0)$ for some suitably large $s_0(t_0)>0$. This shows that, for $s\geq s_0(t_0)$, $\mathcal{B}\subset \widetilde{\mathcal{E}}\cup \widetilde{\mathcal{G}}$, as claimed.

Returning to the proof of \eqref{eq:BcalBd}, by the above proved claim, $\mathcal{B}\subset \widetilde{\mathcal{E}}\cup \widetilde{\mathcal{G}}$ for $s\geq s_0(t_0)$; therefore,
\begin{align}\label{eq:ViaSumBd}
\mathbb{P}\big(\mathcal{B}\big)\leq \mathbb{P}\big(\widetilde{\mathcal{E}}\big)  + \mathbb{P}\big(\widetilde{\mathcal{G}}\big) \, .
\end{align}
We first bound $\mathbb{P}(\widetilde{\mathcal{E}})$. Owing to the definition of $\mathfrak{h}$ in  \eqref{eq:Upsilon} and $\mathfrak{h}_{\alpha t\downarrow t}$ in Proposition~\ref{ppn:TimeEvl}, we have
\begin{equation}\label{eq:EqDist}
\big\{(\alpha-1)^{-1/3}\mathfrak{h}_{\alpha t\downarrow t}(x):x\in \RR\big\} \stackrel{d}{=} \big\{\mathfrak{h}_{(\alpha-1)t}((\alpha-1)^{-2/3}x): x\in \RR\big\}.
\end{equation}
Using this and the change of variables $y= t^{-2/3}x$ for all $x\in [-x_0,x_0]$ implies that, for $\theta := \alpha -1$,
\begin{align}
\mathbb{P}\big(\widetilde{\mathcal{E}}\big) = \mathbb{P}\bigg(\sup_{|y|< \theta^{1/3}} \big\{\mathfrak{h}_{\theta t}\big(\theta^{-2/3}y\big)- \mathfrak{h}_{\theta t}(0)- \theta^{-1/3}\tfrac{y^2}{8}\big\}\geq \theta^{-1/3}\frac{s}{4}\bigg) \, .
\end{align}
Applying Proposition~\ref{cor:HitParabola} with this $\theta$, $\nu=1/4$, $s$ replaced by $s/4$, and $t$ replaced by $\theta t$ (which is still bounded below by some $t'_0>0$ since $\alpha\geq\alpha_0>1$), we see that there exist $s_0=s_0(t_0,\alpha_0)$ and $c=c(t_0,\alpha_0)$ such that, for $t\geq t_0$, $s\geq s_0$ and $\alpha>\alpha_0$, $\mathbb{P}\big(\widetilde{\mathcal{E}}\big)\leq \exp(-cs^{3/4})$.

We now bound $\mathbb{P}\big(\widetilde{\mathcal{G}}\big)$. By the arithmetic-geometric mean inequality, $(\alpha-1)^{-1/3} \Big(\frac{s}{8}+ \frac{k^2_0}{8t^{4/3}}\Big)\geq \frac{ s^{1/2}}{4}$. This, in conjunction with \eqref{eq:EqDist} and the union bound, shows that
 \begin{align}
 \mathbb{P}\big(\widetilde{\mathcal{G}}\big)\leq \mathbb{P}\Big(\sup_{x\in \RR} \mathfrak{h}_{(\alpha-1) t}(x)\geq \frac{s^{1/2}}{4}\Big)+ \mathbb{P}\Big(\mathfrak{h}_{(\alpha-1) t}(0)\leq - \frac{s^{1/2}}{4}\Big). \label{eq:tildeGbBd1}
 \end{align}
Applying Propositions~\ref{ppn:UpperTail} and \ref{ppn:onepointlowertail}, there exist $s_0=s_0(t_0,\alpha_0)$ and $c=c(t_0,\alpha_0)$ such that, for $t>t_0$, $s\geq s^{\prime}_0$ and $\alpha>\alpha_0$, $\mathbb{P}\big(\widetilde{\mathcal{G}}\big)\leq \exp\big(-cs^{3/4}\big)$.
 Substituting the upper bounds on $\mathbb{P}(\widetilde{\mathcal{E}})$ and $\mathbb{P}(\widetilde{\mathcal{G}})$ into \eqref{eq:ViaSumBd} and summing, we arrive at \eqref{eq:BcalBd}.

\smallskip

\textbf{Step III:} By Step I, $\mathcal{A}_{\mathrm{high}}\cap \neg\mathcal{E} \cap \neg\mathcal{G} \subset \mathcal{B}$; hence,
\begin{align}\label{eq:FinalAupBd}
\mathbb{P}(\mathcal{A}_{\mathrm{high}})\leq \mathbb{P}(\mathcal{B}) + \mathbb{P}(\mathcal{E}) + \mathbb{P}(\mathcal{G}) \, .
\end{align}
 We may bound $\mathbb{P}(\mathcal{E}) \leq \exp(-cs^{3/2})$, using Theorem~\ref{t.locreg}); and we may obtain the same bound for~$\mathbb{P}(\mathcal{G})$, using Propositions~\ref{ppn:UpperTail} and~\ref{ppn:onepointlowertail}. Here, we have assumed that  $s>s_0(t_0)$ and $c=c(t_0)$. Combining these bounds with the Step II  bound \eqref{eq:BcalBd}, we arrive at \eqref{eq:SdTail}, and thus complete the proof of Proposition~\ref{lem:UpTail1}. \qed

\subsection{Proof of Proposition~\ref{lem:LowTail}}\label{sec:LowTail} We separately address the two claims \eqref{eq:AlowBd1} and \eqref{eq:AlowBd2}.

\noindent {\bf Proof of \eqref{eq:AlowBd1}.}
 Define events
\begin{align}
\mathcal{W}_1 &:= \Big\{\inf_{x\in [0,1]}\big(\mathfrak{h}_{\alpha t\downarrow t}(t^{-2/3} x)  - \mathfrak{h}_{\alpha t\downarrow t}(0)\big) \leq - t^{-1/3}s/2 \Big\} \, , \\
\mathcal{W}_2 &:= \Big\{\inf_{x\in [0,1]}\big( \mathfrak{h}_{0\uparrow t}(t^{-2/3} x)  - \mathfrak{h}_{0\uparrow t}(0)\big) \leq - t^{-1/3}s/2 \Big\} \, .
\end{align}
On the event $\neg\mathcal{W}_1 \cap \neg\mathcal{W}_2$,
\begin{align}
\mathfrak{h}_{\alpha t\downarrow t}(t^{-2/3}x)\geq \mathfrak{h}_{\alpha t\downarrow t}(0) - t^{-1/3}2^{-1} s \, , \qquad  \mathfrak{h}_{0\uparrow t}(t^{-2/3}x)\geq \mathfrak{h}_{0\uparrow t}(0) - t^{-1/3} 2^{-1}s
\end{align}
for all $x\in [0,1]$. It follows from the composition law Proposition \ref{ppn:TimeEvl} and these inequalities that
\begin{align}\label{eq:IntLowBd}
e^{t^{1/3} \mathfrak{h}_{0\uparrow \alpha t}(0)} \geq \int_{[0,1]} e^{t^{1/3}\big(\mathfrak{h}_{\alpha t\downarrow t}(t^{-2/3} x) + \mathfrak{h}_{0\uparrow t}(t^{-2/3} x)\big)} dx\geq e^{t^{1/3} \mathfrak{h}_{\alpha t\downarrow t}(0)+t^{1/3}\mathfrak{h}_{0\uparrow t}(0) -s} \, .
\end{align}
This shows that $\neg\mathcal{W}_1 \cap \neg\mathcal{W}_2 \subset \neg\mathcal{A}_{\mathrm{low}}$. Hence,
\begin{align}\label{eq:ProbSplit1}
\mathbb{P}\big(\mathcal{A}_{\mathrm{low}}\big)\leq \mathbb{P}(\mathcal{W}_1)+ \mathbb{P}(\mathcal{W}_2) \, .
\end{align}

To bound $\mathbb{P}(\mathcal{W}_1)$ and $ \mathbb{P}(\mathcal{W}_2)$, we set  $\epsilon_1 := ((\alpha-1)t)^{-2/3}$ and $\epsilon_2:=t^{-2/3}$, and rewrite via \eqref{eq:EqDist}
\begin{align}
\mathbb{P}(\mathcal{W}_1) &= \mathbb{P}\Big(\inf_{x\in [0,\epsilon_1]}\big(\mathfrak{h}_{(\alpha-1)t}(x)  - \mathfrak{h}_{(\alpha-1)t}(0)\big) \leq - \frac{s}{2} \epsilon_1^{1/2}\Big),\label{eq:W1Pr}\\
\mathbb{P}(\mathcal{W}_2) &= \mathbb{P}\Big(\inf_{x\in [0,\epsilon_2]}\big(\mathfrak{h}_{t}(x)  - \mathfrak{h}_{t}(0)\big) \leq - \frac{s}{2} \epsilon_2^{1/2}\Big) \, .\label{eq:W2Pr}
\end{align}
We will bound these probabilities via Theorem~\ref{t.locreg}, but two aspects of this application require mention. First, $\epsilon_1$ and $\epsilon_2$ may exceed the upper bound of one assumed on $\epsilon$ in the theorem. These quantities are, however, bounded above by a constant depending on $t_0$ and $\alpha_0$; hence, at the price of degrading the values of $s_0$ and $c$ in the theorem,  variants of the result concerning shifts may be satisfactorily applied. Second, the theorem involves a parabolic shift, whereas the above expressions do not. The parabolic shift can be absorbed by changing the value of $s$. By applying Theorem~\ref{t.locreg} in this manner, the right-hand sides of \eqref{eq:W1Pr} and \eqref{eq:W2Pr} may be bounded above by $\exp(-cs^{3/2})$ for some $c=c(t_0,\alpha_0)>0$, provided that $s>s_0$ for some $s_0=s_0(t_0,\alpha_0)$. Substituting this upper bound into \eqref{eq:ProbSplit1} completes the proof of \eqref{eq:AlowBd1}.

\smallskip

\noindent {\bf Proof of \eqref{eq:AlowBd2}.} In the proof of \eqref{eq:AlowBd1}, we showed that  $\mathcal{A}_{\mathrm{low}} \subset\mathcal{W}_1\cup \mathcal{W}_2$. Since $\mathfrak{h}_{0\uparrow t}$ and $\mathfrak{h}_{\alpha t\downarrow t}$ are independent, we may bound
\begin{align}\label{eq:ProbSplit2}
\mathbb{P}\big(\mathcal{A}_{\mathrm{low}}\big|\mathfrak{h}_{0\uparrow t}(0)\geq y \big)\leq \mathbb{P}\big(\mathcal{W}_1\big)+ \mathbb{P}\big(\mathcal{W}_2\big|\mathfrak{h}_{0\uparrow t}(0)\geq y\big)\leq  \exp(-cs^{3/2}) \big(1+ \mathbb{P}\big(\mathfrak{h}_{0\uparrow t}(0)\geq y\big)^{-1}\big) \, ,
\end{align}
where we set $y :=  \mathbb{E}[\mathfrak{h}_{0\uparrow t}(0)]+\delta$ and where we have used the bounds  on $\mathbb{P}(\mathcal{W}_1)$ and $\mathbb{P}(\mathcal{W}_2)$ established in the proof of \eqref{eq:AlowBd1}.
By the upper bound in Proposition \ref{ppn:onepointuppertail}, we find that, for $t\geq t_0$, $\mathbb{E}[\mathfrak{h}_{0\uparrow t}]<C$ and hence $y<C+\delta$ for $C=C(t_0)>0$. Thus, by the lower bound in Proposition \ref{ppn:onepointuppertail}, $\mathbb{P}\big(\mathfrak{h}_{0\uparrow  t}(0)\geq y \big)>C'$ for $C'=C'(t_0,\delta)>0$. Thus we may bound above the term $1+ \mathbb{P}\big(\mathfrak{h}_{0\uparrow t}(0)\geq y\big)^{-1}$ by a constant. By choosing suitable values of  $c=c(t_0,\alpha_0,\delta)$ and $s_0=s_0(t_0,\alpha_0,\delta)$, we may absorb this constant and thus show that, for $t>t_0$, $s>s_0$ and $\alpha>\alpha_0$,
$
\mathbb{P}\big(\mathcal{W}_2\big|\mathfrak{h}_{0\uparrow t}(0)\geq  y_0 \big) \leq \exp(-cs^{3/2}),
$
as desired to demonstrate \eqref{eq:AlowBd2}.

\section{Adjacent Correlation: Proof of Theorem~\ref{thm:Main2}}\label{sec:adjacent}

We are now concerned with the correlation between $\mathfrak{h}_{t}(1+\beta,0)$ and $\mathfrak{h}_{t}(1,0)$. As in Section \ref{sec:remote}, we will rely on the composition law Proposition \ref{ppn:TimeEvl} to realize $\mathfrak{h}_{t}(1+\beta,0)$ in terms of $\mathfrak{h}_{t}(1,\cdot)$ and $\mathfrak{h}_{(1+\beta)t\downarrow t}(\cdot)$. In this section, we will use a variation on the shorthand from Section \ref{sec:remote}:
\begin{align}
\mathfrak{h}_{0\uparrow t}(\cdot):=\mathfrak{h}_t(1,\cdot),\qquad \mathfrak{h}_{0\uparrow (1+\beta) t}(\cdot):=\mathfrak{h}_t(1+\beta,\cdot) \, ;
\end{align}
and by $\mathfrak{h}_{0\uparrow t}$ is meant  $\mathfrak{h}_{0\uparrow t}(0)$ -- and likewise for $\mathfrak{h}_{0\uparrow (1+\beta) t}$ and $\mathfrak{h}_{(1+\beta) t\downarrow t}$.
Note may continue to be taken of the caveat that warns of the confusion that this shorthand may provoke.  Theorem~\ref{thm:Main2}'s conclusion~\eqref{eq:Corr2} in the present notation asserts that
\begin{align}\label{eq:Corr2prime}
c_1 \beta^{2/3}\leq 1-\mathrm{Corr}\big( \mathfrak{h}_{0\uparrow t}, \mathfrak{h}_{0\uparrow (1+\beta)t}\big)\leq c_2\beta^{2/3} \, .
\end{align}
Theorem \ref{ppn:HDivTail} and the next stated Proposition \ref{cl:VarAux} will yield these two bounds. We state the new proposition; prove Theorem~\ref{thm:Main2}; and, in two ensuing subsections, prove Theorem~\ref{ppn:HDivTail} and Proposition~\ref{cl:VarAux}.

%
%
%

  \bp\label{cl:VarAux}
For any $t_0>0$, there exist $s_0=s_0(t_0)$ and $c=c(t_0)$  such that, for $s>s_0$, $t>t_0$ and  $\beta\in (0,\frac{1}{2})$ satisfying $\beta t>t_0$,
 \begin{align}\label{eq:ShiftDiff}
 \mathbb{P}\bigg(\mathfrak{h}_{0\uparrow (1+\beta)t}-\mathfrak{h}_{0\uparrow t} - \frac{\mathrm{Cov}\big(\mathfrak{h}_{0\uparrow (1+\beta)t}-\mathfrak{h}_{0\uparrow t}, \mathfrak{h}_{0\uparrow t}\big)}{\mathrm{Var}(\mathfrak{h}_{0 \uparrow t})} \mathfrak{h}_{0 \uparrow t}\geq \beta^{1/3}s\bigg) \, \geq \, \exp\big(-cs^{3/2}\big) \, .
 \end{align}
 \ep

\begin{proof}[Proof of Theorem~\ref{thm:Main2}]
Aiming to apply Lemma~\ref{lem:GenCorrBd}, suppress $t$ and $\beta$ in notation that sets  $X:= \mathfrak{h}_{0\uparrow (1+\beta)t}-\mathfrak{h}_{0\uparrow t}$, $Y:= \mathfrak{h}_{0 \uparrow t}$ and $Z:= \mathfrak{h}_{0 \uparrow (1+\beta)t}$; and, as in Lemma~\ref{lem:GenCorrBd}, set
\begin{align}
\chi := \frac{\mathrm{Var}(X)}{\mathrm{Var}(Y)},\qquad \Psi := \frac{\mathrm{Cov}(X,Y)}{\mathrm{Var}(Y)}, \qquad  \Theta := \frac{\mathrm{Var}(X)\mathrm{Var}(Y)- (\mathrm{Cov}(X,Y))^2}{(\mathrm{Var}(Y))^2} \, .
\end{align}
We may argue
that there exist $c_1=c_1(t_0)<c_2=c_2(t_0)$ such that, for all $t> t_0/\beta$,
\begin{align}\label{eq:VarUpLow}
c_1\leq  \mathrm{Var}(Y)\leq  c_2, \qquad c_1 \beta^{2/3}\leq \mathrm{Var}(X)\leq c_2 \beta^{2/3} \, .
\end{align}
Indeed, the bounds on $\mathrm{Var}(Y)$ follow by combining the tail bounds from Propositions \ref{ppn:onepointlowertail} and \ref{ppn:onepointuppertail} with the later presented tool, Lemma~\ref{ppn:VarBound}, that translates tail bounds into bounds on variance.  For the bounds on $\mathrm{Var}(X)$, it is  Theorem~\ref{ppn:HDivTail} that instead furnishes the needed tail bounds.

In view of these bounds and the Cauchy-Schwarz inequality in the guise $\Psi^2\leq \chi$, there exists $\beta_0\in (0,\frac{1}{2})$ such that, for $\beta<\beta_0$ and $t> t_0/\beta$,
$\max\big\{\chi, |\Psi+\frac{1}{2}\chi|\big\} <1$.
This verifies the hypothesis of Lemma~\ref{lem:GenCorrBd} and thus demonstrates the existence of constants $C_1<C_2$ for which $\beta<\beta_0$ and $t>t_0/\beta$ imply that
\begin{align}\label{eq:CorrUpLow}
1- \Theta/2 +C_1\chi^{3/2}\leq \mathrm{Corr}(Z,Y) \leq 1- \Theta/2 +C_2\chi^{3/2} \, .
\end{align}
We {\em claim} that
\begin{align}\label{e.threeclaim}
c_1 \beta^{2/3} \leq \Theta \leq \chi \leq c_2 \beta^{2/3} \, ,
\end{align}
where here we assume that  $\beta<\beta_0$ and $0< c_1=c_1(t_0,\beta_0)<c_2=c_2(t_0,\beta_0)$. This claim alongside~\eqref{eq:CorrUpLow} proves \eqref{eq:Corr2prime} for $\beta<\beta_0$. The result is extended to $\beta \in [\beta_0,1/2)$ by possibly changing the constants and recognizing that the correlations are bounded above by one; thus the proof of Theorem~\ref{thm:Main2} is complete subject to verifying the three bounds in the claim~\eqref{e.threeclaim}.

Consider the first claimed bound  $c_1 \beta^{2/3} \leq \Theta$.
Observe that $\Theta= \frac{\mathbb{E}[(X-\Psi Y)^2]}{\mathrm{Var}(Y)}$. We may bound $\mathrm{Var}(Y)$ above by a constant using \eqref{eq:VarUpLow}. Thus it remains to bound $\mathbb{E}\big[(X-\Psi Y)^2\big]\geq c \beta^{2/3}$ for some $c=c(t_0,\beta_0)$. To do this, we will appeal to the second part of  Lemma~\ref{ppn:VarBound}. There are two hypotheses that we must verify for that application.  The first is that $\big|\EE[X-\Psi Y]\big|\leq C\beta^{1/3}$ for some $C=C(t_0,\beta_0)>0$, and the second is the upper-tail lower bound on $X-\Psi Y$ already given in Proposition \ref{cl:VarAux}. Predicated upon showing these bounds, the second part of  Lemma~\ref{ppn:VarBound} immediately implies $\mathbb{E}\big[(X-\Psi Y)^2\big]\geq c \beta^{2/3}$ as desired.
To show that $\big|\EE[X-\Psi Y]\big|\leq C\beta^{1/3}$ we demonstrate that $\EE[X]$ converges to 0; $\Psi\leq C\beta^{1/3}$; and $\EE[Y]<C$ for some $C=C(t_0,\beta_0)>0$. That $\EE[X]\to 0$ is shown by combining the convergence of $\mathbb{E}[\mathfrak{h}_{0\uparrow t}]$ and $\mathbb{E}[\mathfrak{h}_{0\uparrow (1+\beta) t}]$ to a common limit (this due to the weak convergence result of Proposition \ref{ppn:Stationarity}) with the tail bounds afforded in Propositions  \ref{ppn:onepointlowertail} and \ref{ppn:onepointuppertail}. The Cauchy-Schwarz inequality  yields $\Psi \leq \sqrt{\frac{\mathrm{Var}(X)}{\mathrm{Var}(Y)}}$, and the inference $\Psi \leq C\beta^{1/3}$ is then made by applying~\eqref{eq:VarUpLow}; these bounds in \eqref{eq:VarUpLow} also show that $\EE[Y]<C$.

The second claimed bound $\Theta \leq \chi$ in \eqref{e.threeclaim} follows since $\Theta = \chi - \Psi^2$, and the third bound $\chi\leq c_2 \beta^{2/3}$ is due to \eqref{eq:VarUpLow}. Thus the claim~\eqref{e.threeclaim} is verified and the proof of Theorem~\ref{thm:Main2} completed.
\end{proof}

\subsection{Proof of Theorem~\ref{ppn:HDivTail}}\label{sec:HDivTail}

There are three inequalities to prove, and the proof is divided into three numbered {\em stages} accordingly. The notation  $\mathfrak{h}_{0\uparrow t}$ in the theorem's statement becomes   $\mathfrak{h}_{0\uparrow t}(0)$ in the proof, because we have cause to consider $\mathfrak{h}_{0\uparrow t}(x)$ for general values of $x$.

\begin{proof}[Stage $1$. Proof of $\mathbb{P}\big(\mathfrak{h}_{0\uparrow (1+\beta)t}-\mathfrak{h}_{0\uparrow t} \geq \beta^{1/3}s\big)\leq \exp\big(-c_2 s^{3/4}\big)$]

%

Seeking to bound $\mathbb{P}(\mathcal{A}_{\textrm{high}})$, define
\begin{align}\label{eq:DefA2}
\mathcal{A}_{\textrm{high}}&:=\Big\{\mathfrak{h}_{0\uparrow (1+\beta)t}(0)- \mathfrak{h}_{0\uparrow t}(0) \geq\beta^{1/3}s \Big\} \, ,\\
\mathcal{E}&:=\Big\{\sup_{x\in \RR}\Big\{\mathfrak{h}_{(1+\beta)t\downarrow t}\big(t^{-2/3}x\big) +\frac{x^2}{4\beta t^{4/3}} \Big\}\geq \beta^{1/3}\frac{s}{2}\Big\} \, .
\end{align}
By the union bound, $\mathbb{P}\big(\mathcal{A}_{\textrm{high}}\big)\leq \mathbb{P}\big(\mathcal{E}\big) + \mathbb{P}\big(\mathcal{A}_{\textrm{high}}\cap \neg\mathcal{E}\big)$; it thus suffices to bound the probabilities of the two right-hand events by expressions of the form $\exp\big(-c s^{3/4}\big)$.
Bounding $\mathbb{P}\big(\mathcal{E}\big)$ is easier:
 by the definition of~$\mathfrak{h}$,
\begin{align}
\Big\{\mathfrak{h}_{(1+\beta)t\downarrow t}(\beta^{2/3}x):x\in \RR\Big\}\stackrel{d}{=}\Big\{\beta^{1/3}\mathfrak{h}_{\beta t}(x):x\in \RR\Big\} \, .\label{eq:DistEq2}
\end{align}
Via the change of variables $y= (\beta t)^{-2/3}x$  and this distributional identity,
\begin{align}
\mathbb{P}(\mathcal{E}) = \mathbb{P}\Big(\sup_{y\in \RR}\Big\{\mathfrak{h}_{\beta t}\big(y\big) +\frac{x^2}{4} \Big\}\geq \frac{s}{2}\Big) \, .
\end{align}
By Proposition~\ref{ppn:UpperTail}, there exist $s_0=s_0(t_0)>0$ and $c_0=c_0(t_0)$ such that, for $s>s_0$, $t>t_0$ and  $\beta\in (0,\frac{1}{2})$ satisfying $\beta t>t_0$,
$\mathbb{P}(\mathcal{E})\leq \exp(-cs^{3/2})$ -- thus obtaining the desired bound on $\PP(\mathcal{E})$.

We turn to demonstrating that  $\mathbb{P}\big(\mathcal{A}_{\textrm{high}}\cap \neg\mathcal{E} \big)\leq \exp\big(-c s^{3/4}\big)$. Set  $y_0:= (\beta t^2)^{1/3}$ and
\begin{align}\label{eq:DefB2}
\mathcal{D}&:=\bigg\{ \int^{\infty}_{-\infty}\exp\Big(t^{1/3}\mathfrak{h}_{0\uparrow t}(t^{-2/3}x)-\frac{x^2}{4\beta t}\Big) dx\geq \exp\Big(t^{1/3} \mathfrak{h}_{0\uparrow t}(0)+(\beta t)^{1/3}\frac{s}{2}\Big)\bigg\} \, ,\\
\mathcal{Q}&:=\bigg\{\sup_{|x|< y_0}\Big\{\mathfrak{h}_{0\uparrow t}(t^{-2/3}x) - \mathfrak{h}_{0\uparrow t}(0) - \frac{x^2}{8\beta t^{4/3}}\Big\}\geq \beta^{1/3}\frac{s}{4} \bigg\} \, , \\
\mathcal{G}&:=\bigg\{\sup_{|x|\geq y_0 }\mathfrak{h}_{0\uparrow t}(t^{-2/3}x)\geq \frac{y^2_0}{16\beta t^{4/3}}+\frac{\beta^{1/3}s}{8}\bigg\}\bigcup \bigg\{\mathfrak{h}_{0\uparrow t_1}(0)\leq -\frac{y^2_0}{16\beta t^{4/3}}-\frac{\beta^{1/3}s}{8}\bigg\} \, .
\end{align}
To obtain the presently sought bound, we show in three steps that $\mathcal{A}_{\textrm{high}}\cap \neg\mathcal{E} \subset \mathcal{D}$; that $\mathcal{D}\subset \mathcal{Q}\cup \mathcal{G}$; and that $\mathbb{P}(\mathcal{Q}\cup \mathcal{G}) \leq \exp(-cs^{3/4})$.

{\bf Step I:} We seek to show that $\mathcal{A}_{\textrm{high}}\cap \neg\mathcal{E} \subseteq \mathcal{D}$. When $\neg\mathcal{E}$ occurs, we have that, for $x\in \RR$,
  \begin{align}\label{eq:Eneg}
  t^{1/3}\mathfrak{h}_{(1+\beta)t\downarrow t}(t^{-2/3}x)\leq -\frac{x^2}{4\beta t} +(\beta t)^{1/3} \frac{s}{2} \, .
  \end{align}
Applying  the composition law Proposition \ref{ppn:TimeEvl} to write $\mathfrak{h}_{0\uparrow (1+\beta)t}(0)$ in terms of $\mathfrak{h}_{0\uparrow t}(\cdot)$ and $\mathfrak{h}_{(1+\beta)t\downarrow t}(\cdot)$,
\begin{align}
\exp\big(t^{1/3}\mathfrak{h}_{0\uparrow (1+\beta)t}(0)\big)\leq \exp\Big((\beta t)^{1/3}\frac{s}{2}\Big)\int^{\infty}_{-\infty}\exp\Big(t^{1/3}\mathfrak{h}_{0\uparrow t}(t^{-2/3} x) -  \frac{x^2}{4\beta t}\Big) dx \, .\label{eq:CoupPlug}
\end{align}
On the event $\mathcal{A}_{\textrm{high}}$, this left-hand side is bounded below by $\exp\big(t^{1/3}\mathfrak{h}_{0\uparrow t}(0)+(\beta t)^{1/3}s\big)$. Since the bound specifying $\mathcal{D}$ has been verified, the desired containment has been shown.
\smallskip

{\bf Step II:} We seek to show that there exists $s_0>0$ such that, for $s\geq s_0$, $\mathcal{D}\subset \mathcal{Q}\cup \mathcal{G}$. We will show the contrapositive $\neg \mathcal{Q} \cap \neg\mathcal{G} \subset \neg\mathcal{D}$.
Note first that
\begin{align}
\int^{\infty}_{-\infty} &\exp\Big(t^{1/3}\mathfrak{h}_{0\uparrow t}(t^{-2/3}x)- \frac{x^2}{4\beta t} \Big) dx\\
& \, \, \leq \exp\Big(t^{1/3}\sup_{|y|\leq y_0}\Big\{\mathfrak{h}_{0\uparrow t}(t^{-2/3} y)-\frac{y^2}{8\beta t^{4/3}}\Big\}\Big)\int_{|y|\leq y_0} \exp\big(-\frac{x^2}{8\beta t}\big) dx \\
& \quad + \, \,  \exp \Big(t^{1/3}\sup_{|y|\geq y_0}\mathfrak{h}_{0\uparrow t}(t^{-2/3} y)\Big) \int_{|y|\geq y_0} \exp\big(-\frac{x^2}{8\beta t}\big) dx\, .  \label{eq:Ineq1}
\end{align}
Here, the first right-hand term is present because, when $|x|\leq y_0$, $\mathfrak{h}_{0\uparrow t}(t^{-2/3}x)-\frac{x^2}{4\beta t}$ is at most $\sup_{|y|\leq y_0} \big( \mathfrak{h}_{0\uparrow t}(t^{-2/3}y)-\frac{y^2}{8\beta t} \big) - \frac{x^2}{8\beta t}$. The second right-hand term appears because  $\mathfrak{h}_{0\uparrow t}(t^{-2/3}x)-\frac{x^2}{4\beta t}$ is at most $\sup_{|y|\geq y_0} \big\{\mathfrak{h}_{0\uparrow t}(t^{-2/3}y)\big\}-\frac{x^2}{8\beta t}$ when  $|x|\geq y_0$.

On the event $\neg \mathcal{Q} \cap \neg\mathcal{G}$, we have the bounds
\begin{align}
\exp\Big(t^{1/3}\sup_{|x|< y_0}\Big(\mathfrak{h}_{0\uparrow t}(t^{-2/3} x)-\frac{x^2}{8\beta t^{4/3}}\Big)\Big) &\leq \exp\Big(t^{1/3}\mathfrak{h}_{0\uparrow t}(0)+ (\beta t)^{1/3}\frac{s}{4}\Big) \, ,\label{eq:Theta1}\\
\exp\Big(t^{1/3}\sup_{|x|\geq y_0}\mathfrak{h}_{0\uparrow t}(t^{-2/3} x)\Big)&\leq \exp\Big(t^{1/3}\mathfrak{h}_{0\uparrow t}(0) + \frac{y^2_0}{8\beta t} + \frac{(\beta t)^{1/3}s}{4}\Big) \, . \label{eq:Theta3}
\end{align}
Substitute these into \eqref{eq:Ineq1} and note that $\int_{|x|\leq y_0}\exp\Big(-\frac{x^2}{8\beta t}\Big) dy$ and $ \int_{|x|\geq y_0} \exp\Big(-\frac{x^2-y^2_0}{8\beta t}\Big) dx $ are bounded above by $c(\beta t)^{\frac{1}{2}}$ for some $c>0$.
What we learn is that, for some constant $c>0$,
\begin{align}
\int^{\infty}_{-\infty} &\exp\Big(t^{1/3}\mathfrak{h}_{0\uparrow t}(t^{-2/3}x)- \frac{x^2}{4\beta t} \Big) dx \leq c (\beta t)^{\frac{1}{2}} \exp\Big(t^{1/3}\mathfrak{h}_{0\uparrow t}(0)+ (\beta t)^{1/3}\frac{s}{4}\Big).\label{eq:rhsBd}
\end{align}
If  $s> 4\log\big(c (\beta t)^{\frac{1}{2}}\big) (\beta t)^{-1/3}$ then $c (\beta t)^{1/2} \leq \exp\big((\beta t)^{1/3} \frac{s}{4}\big)$. Since $r^{-1} \log r$ is bounded for $r>0$, we find that there exists $s_0>0$ such that, for $s>s_0$, the right-hand side of ~\eqref{eq:rhsBd} is at most $\exp\big(t^{1/3}\mathfrak{h}_{0\uparrow t}(0)+ (\beta t)^{1/3}\frac{s}{2}\big)$. Since this bound specifies the event $\neg\mathcal{D}$, the conclusion of this second step has been obtained.
%


{\bf Step III:} We seek to show that there exist $s_0=s_0(t_0)>0$ and $c=c(t_0)>0$ such that, for $s\geq s_0$,
$
\mathbb{P}\big( \mathcal{Q} \cup \mathcal{G}\big) \leq  \exp(-cs^{3/4})$.
Recall that the spatial process $\mathfrak{h}_{0\uparrow t}(\cdot)$ has same law as $\mathfrak{h}_{t}(\cdot)$. Apply the change of variables $x= (\beta t)^{2/3} y$ in the definition of $\mathcal{Q}$, and set $\theta =\beta^{-1}$, to obtain
$$
\mathcal{Q} = \bigg\{\sup_{|y|<\theta^{1/3}} \Big\{\mathfrak{h}_t(\theta^{-2/3} y) - \mathfrak{h}_t(0) - \theta^{-1/3} \frac{y^2}{8}\Big\}  \geq \theta^{-1/3} \frac{s}{4}\bigg\} \, .
$$
Since $\beta<1$ by assumption, $\theta>1$. Owing to this, we can apply Proposition~\ref{cor:HitParabola} with $\nu=1/4$. This result  controls the supremum only over positive $\theta$. However, using the reflection invariance of $\mathfrak{h}_t$ given in Proposition \ref{ppn:Stationarity}, we may likewise control the supremum over negative $\theta$. Allying the resulting inferences with the union bound, we conclude that there exist $s_0=s_0(t_0)$ and $c=c(t_0)$ such that, for $s\geq s_0$,
$\mathbb{P}(\mathcal{Q}) \leq  \exp(-cs^{3/4})$.

 Using the upper bound on the upper tail for the law of $\sup_{x\in \RR}\mathfrak{h}_{0\uparrow t}(x)$ from Proposition~\ref{ppn:LowerTail} and the upper bound on the lower tail of the law of  $\mathfrak{h}_{0\uparrow t}(0)$ from Proposition \ref{ppn:onepointlowertail} shows that there exist $s_0=s_0(t_0)>0$, $c=c(t_0)>0$ and $c^\prime=c^\prime(t_0)>0$ so that
\begin{align}\label{eq:GBd}
\mathbb{P}(\mathcal{G})\leq \exp\Big(-c\Big(\frac{1}{16\beta^{1/3}}+ \frac{\beta^{1/3} s}{8}\Big)^{3/2}\Big)\leq \exp(-c^{\prime}s^{3/4}) \, ,
\end{align}
where the latter inequality uses  $\frac{1}{16\beta^{1/3}}+ \frac{\beta^{1/3} s}{8}\geq \sqrt{\frac{s}{32}}$ (via the arithmetic-geometric mean inequality).

Combining this bound on $\mathbb{P}(\mathcal{G})$ with the above bound on $\mathbb{P}(\mathcal{Q})$ yields the conclusion of Step III. Thus do we obtain the sought bound $\mathbb{P}\big(\mathcal{A}_{\textrm{high}}\cap \neg\mathcal{E} \big)\leq \exp\big(-c s^{3/4}\big)$ completing the derivation of the first of three inequalities asserted by Theorem~\ref{ppn:HDivTail}, namely  $\mathbb{P}\big(\mathfrak{h}_{0\uparrow (1+\beta)t}-\mathfrak{h}_{0\uparrow t} \geq \beta^{1/3}s\big)\leq \exp\big(-c_2 s^{3/4}\big)$.
%
\end{proof}

\begin{proof}[Stage $2$. Proof of $\exp\big(-c_1s^{3/2}\big) \leq \mathbb{P}\big(\mathfrak{h}_{0\uparrow (1+\beta)t}-\mathfrak{h}_{0\uparrow t} \geq \beta^{1/3}s\big)$]

Fix $\theta := (\beta t)^{2/3}s^{-1}$ and define events
\begin{align}
\mathcal{C}&:= \bigcap_{\epsilon\in \{-1,1\}}\Big\{\mathfrak{h}_{(1+\beta)t\downarrow t}(\epsilon t^{-2/3}\theta)\geq 2\beta^{1/3}s\Big\} \, ,
&&\mathcal{D}:=  \Big\{\inf_{|y|\leq \theta} \mathfrak{h}_{(1+\beta)t\downarrow t}( t^{-2/3}y)\geq \beta^{1/3}\frac{ 7s}{4}\Big\}\, , \\
\mathcal{E}&:= \Big\{\inf_{|y|\leq \theta}\mathfrak{h}_{0\uparrow t}(t^{-2/3} y)-\mathfrak{h}_{0\uparrow t}(0)\geq -\beta^{1/3}\frac{s}{2}\Big\} \, ,
&&\mathcal{W}: = \mathcal{C} \cap \mathcal{D}\cap \mathcal{E} \, .
\end{align}
The proof has two steps. First, we show that, for $s$ large enough, $\mathcal{W}\subset \mathcal{A}_{\textrm{high}}$; second, we find a lower bound of the form $\mathbb{P}(\mathcal{W})\geq \exp(-c s^{3/2})$.

{\bf Step I:} We seek to show that $\mathcal{W}\subset \mathcal{A}_{\textrm{high}}$ for $s>s_0$ for $s_0=s_0(t_0)>0$.  On the event $\mathcal{W}$, we have
\begin{align}
t^{1/3}\mathfrak{h}_{(1+\beta)t\downarrow t}( t^{-2/3}y)+ t^{1/3}\mathfrak{h}_{0\uparrow t}(t^{-2/3} y) \geq t^{1/3} \mathfrak{h}_{0\uparrow t}(0)+ (\beta t)^{1/3}\frac{5 s}{4}, \quad \forall y\in [-\theta, \theta].
\end{align}
Substituting this inequality into the composition law Proposition~\ref{ppn:TimeEvl} yields
\begin{align}\label{eq:Containment4}
\exp\big(t^{1/3}\mathfrak{h}_{0\uparrow t}(0)\big) & \geq \int_{-\theta}^{\theta} \exp\Big(t^{1/3}\mathfrak{h}_{0\uparrow t}(0)+(\beta t)^{1/3}\frac{5s}{4}\Big) dx
\\&= 2\exp\Big(t^{1/3}\mathfrak{h}_{0\uparrow t}(0)+(\beta t)^{1/3}\frac{5s}{4} -\log s + \tfrac{2}{3} \log(\beta t)\Big).
\end{align}
The quantity in the preceding line is greater than $\exp(t^{1/3}\mathfrak{h}_{0\uparrow t}(0)+(\beta t)^{1/3} s)$ provided that  $s>s_0$ for some $s_0(t_0)>0$, because $\beta t >t_0$ by assumption. The event $\mathcal{A}_{\textrm{high}}$ is specified by the resulting inequality. This shows  that $\mathcal{W}\subset \mathcal{A}_{\textrm{high}}$.

{\bf Step II:} We seek to prove that $\mathbb{P}(\mathcal{W})\geq \exp\big(-cs^{3/2}\big)$ where $c=c(t_0)>0$ and $s>s_0$ for $s_0=s_0(t_0)>0$. Since $\mathfrak{h}_{(1+\beta)t\downarrow t}$ and $\mathfrak{h}_{0\uparrow t}$ are independent,
\begin{align}
\mathbb{P}(\mathcal{W}) =\mathbb{P}(\mathcal{C}\cap \mathcal{D}) \, \mathbb{P}(\mathcal{E}) \, .
%
\end{align}
By setting $\epsilon = \beta^{2/3}s^{-1}$, we may write  $\mathcal{E}=\big\{\inf_{|y|\leq \epsilon}\mathfrak{h}_{0\uparrow t}(t^{-2/3} y)-\mathfrak{h}_{0\uparrow t}(0)\geq -\epsilon^{1/2}\frac{s^{3/2}}{2}\big\}$.  Thus, we may apply Theorem \ref{t.locreg} to find that
there exist $s_0=s_0(t_0)>0$ and $c=c(t_0)>0$ such that $\mathbb{P}(\mathcal{E}) \geq 1-\exp\big(-c s^{9/4}\big)$, provided that $s>s_0$ and $\beta t \geq t_0$. We see then that the lower bound on $\mathbb{P}(\mathcal{W})$ will follow once we derive $\mathbb{P}(\mathcal{C}\cap \mathcal{D}) \geq \exp\big(-cs^{3/2}\big)$. Combining the distributional equality~\eqref{eq:DistEq2}, the stationarity in Proposition \ref{ppn:Stationarity} and the lower bound in Proposition \ref{ppn:onepointuppertail}, we find that $\mathbb{P}(\mathcal{C}) \geq \exp\big(-cs^{3/2}\big)$.  By decomposing $\mathbb{P}(\mathcal{C}\cap \mathcal{D} )  = \mathbb{P}(\mathcal{C})  - \mathbb{P}(\mathcal{C}\cap \neg \mathcal{D} )$ we see that the present Step II will be completed by obtaining $\mathbb{P}(\mathcal{C}\cap \neg \mathcal{D} )\leq   \exp\big(-cs^{\frac{3}{2}+\eta}\big)$ for any fixed $\eta>0$. We will use the KPZ line ensemble to derive such a bound with the choice $\eta=1$.

Using \eqref{eq:DistEq2}, Proposition \ref{NWtoLineEnsemble}, and changing variables, we express $\PP \big( \mathcal{C}\cap \neg \mathcal{D} \big)$ as a probability concerning the lowest indexed curve $\mathfrak{h}^{(1)}_{\beta t}(\cdot)$ of the KPZ$_{\beta t}$ line ensemble, and bound the resulting probability:
\begin{align}
\mathbb{P}\big(\mathcal{C}\cap \neg \mathcal{D} \big) & = \mathbb{P}\bigg(\bigcap_{\epsilon\in\{-1,1\}}\Big\{ \mathfrak{h}^{(1)}_{\beta t}(\epsilon s^{-1})\geq 2 s\Big\} \bigcap \Big\{\inf_{|y|\leq s^{-1}} \mathfrak{h}^{(1)}_{\beta t}(y) \leq \frac{7 s}{4}\Big\}\bigg) \\
& \leq \mathbb{P}_{B(\pm s^{-1}) = 2s}\Big(\inf_{|y|\leq s^{-1}} B(y) \geq \frac{7 s}{4}\Big) \leq \exp\big(-c^{\prime} s^{5/2}\big) \, . \label{eq.cdneg}
\end{align}
Here, $B$ is a Brownian bridge on $[-s^{-1},s^{-1}]$ with starting and ending values $2s$, and $c'>0$ is a constant. The first upper bound in \eqref{eq.cdneg} is due to an appeal
to the Brownian Gibbs property in view of the coupling in Lemma \ref{Coupling1}. That is, we write the concerned probability as the expectation of the indicator of the first event $\bigcap_{\epsilon\in\{-1,1\}}\big\{ \mathfrak{h}^{(1)}_{\beta t}(\epsilon s^{-1})\geq 2 s\big\}$ multiplied by the conditional expectation (with respect to the sigma field generator by everything outside of the first curve on the interval $[-s^{-1},s^{-1}]$)  of the indicator for the second event. By Lemma \ref{Coupling1} and the nature of the second event, this conditional expectation only increases when the second curve drops to $-\infty$ and the starting and ending points drop to $2s$. Doing this, and then bounding the indicator for the first event by $1$, yields the first bound in~\eqref{eq.cdneg}.
 The second bound in~\eqref{eq.cdneg} is due to the reflection principle.

 As we have noted, the bound~\eqref{eq.cdneg} yields Step II, namely the bound $\mathbb{P}(\mathcal{W})\geq \exp\big(-cs^{3/2}\big)$. Steps I and II completed, the second of the three stages of the proof of Theorem~\ref{ppn:HDivTail} is also finished.
\end{proof}

\begin{proof}[Stage~$3$. Proof of $\mathbb{P}\big(\mathfrak{h}_{0\uparrow (1+\beta)t}-\mathfrak{h}_{0\uparrow t}\leq -\beta^{1/3}s\big)\leq \exp\big(-c_3s^{3/2}\big)$]
In order to bound $\mathbb{P}\big(\mathcal{A}_{\mathrm{low}} \big)$, set
\begin{align}\label{eq:DefALow}
\mathcal{A}_{\mathrm{low}} &:=\Big\{\mathfrak{h}_{0\uparrow (1+\beta)t} - \mathfrak{h}_{0\uparrow t}(0) \leq-\beta^{1/3}s\Big\} \, ,\\
\mathcal{\widetilde E} &:= \Big\{\inf_{x\in \RR}\big\{\mathfrak{h}_{(1+\beta)t\downarrow t}(t^{-2/3}x)+ \frac{(1+\nu)}{2\beta t^{4/3}}\big\}\leq -\beta^{1/3}\frac{s}{4}  \Big\} \, ,\label{eq:DefELow}
\end{align}
where $\nu\in (0,1)$ is arbitrary and fixed. The union bound shows that $\mathbb{P}(\mathcal{A}_{\mathrm{low}})\leq \mathbb{P}(\mathcal{\widetilde E}) + \mathbb{P}(\mathcal{A}_{\mathrm{low}}\cap \neg \mathcal{\widetilde E})$, so that our task is to bound the two right-hand terms. Regarding the first (and easier) of the two, recall the distributional identity in \eqref{eq:DistEq2} and change variables $y=(\beta t)^{-2/3}x$ to find that
\begin{align}
\mathbb{P}(\mathcal{\widetilde E}) &=\,  \mathbb{P}\, \bigg(\, \inf_{y\in \RR}\Big(\mathfrak{h}_{\beta t}(y)+\frac{(1+\nu)}{2}x^2\Big)\leq -\frac{s}{4}\, \bigg)  \, .
\end{align}
By Proposition~\ref{ppn:LowerTail}, there exist $s_0=s_0(t_0)>0$ and $c=c(t_0)>0$ such that, for $s> s_0$ and $ t >\frac{t_0}{\beta}$,
 \begin{align}
 \mathbb{P}\big(\mathcal{\widetilde E} \big)\leq \exp(-cs^{5/2}) \, ,
 \end{align}
 which is the first of the two bounds that we seek.
The second is an upper bound on $\mathbb{P}\big(\mathcal{A}_{\mathrm{low}}\cap \neg \mathcal{\widetilde E} \big)$. Fixing $y_0:= (\beta t^2)^{1/3}$, define events
\begin{align}
\mathcal{\widetilde D}&:= \bigg\{\int^{\infty}_{-\infty} \exp\Big(t^{1/3} \mathfrak{h}_{0\uparrow t}(t^{-2/3}x)- (1+\nu)\frac{x^2}{2\beta t}\Big)dx \leq \exp\Big(t^{1/3}\mathfrak{h}_{0\uparrow t}(0)- (\beta t)^{1/3}\frac{3s}{4}\Big)\bigg\} \, ,\\
\mathcal{\widetilde Q}&:= \bigg\{\inf_{|x|< y_0}\Big(\mathfrak{h}_{0\uparrow t}(t^{-2/3}x)- \mathfrak{h}_{0\uparrow t}(0)+\frac{(1+\nu)x^2}{2\beta t^{4/3}}\Big)\leq -\beta^{1/3}\frac{s}{2}\bigg\} \, ,\\
\mathcal{\widetilde G}&:= \bigg\{\inf_{|y|\geq y_0}\mathfrak{h}_{0\uparrow t}(t^{-2/3}y)\leq -\frac{(1+\nu)y^2_0}{2\beta t^{4/3}} -\frac{\beta^{1/3}s}{4}\bigg\} \, \bigcup\,  \bigg\{\mathfrak{h}_{0\uparrow t}(0)\geq \frac{(1+\nu)y^2_0}{2\beta t^{4/3}} +\frac{\beta^{1/3}s}{4}\bigg\} \, .
\end{align}
These events are naturally similar to $\mathcal{D}$, $\mathcal{Q}$, and $\mathcal{G}$ from the preceding proof of the upper bound on  $\mathbb{P}(\mathcal{A}_{\mathrm{high}})$, but the inequalities are each of opposite type. The same three steps govern the proof: in \textbf{Step I}, we show that  $\mathcal{A}_{\mathrm{low}}\cap \neg \mathcal{\widetilde E} \subset \mathcal{\widetilde D}$; in \textbf{Step II}, that  $\mathcal{\widetilde D}\subset \mathcal{\widetilde Q}\cup\mathcal{\widetilde G}$; and in \textbf{Step III}, we find an upper bound on $\mathbb{P}\big(\mathcal{\widetilde Q}\cup\mathcal{\widetilde G}\big)$. The logic of the argument in each step is unchanged from its counterpart's in the derivation of an upper bound on $\mathbb{P}(\mathcal{A}_{\mathrm{low}})$, and we do not record the arguments again.
\end{proof}


\subsection{Proof of Proposition~\ref{cl:VarAux}}\label{sec:VarAux}
The derivation runs in parallel with Stage~$2$ of Theorem \ref{ppn:HDivTail}'s proof.
Propositions~\ref{ppn:onepointlowertail} and ~\ref{ppn:onepointuppertail}, Theorem ~\ref{ppn:HDivTail} and Lemma~\ref{ppn:VarBound} imply that there exist $0\!<\!C_2\!=\!C_2(t_0)\!\leq\!C_1\!=\!C_1(t_0)$ for which
\begin{align}
C_2\leq \mathrm{Var}(\mathfrak{h}_{0\uparrow t})\leq C_1 \, , \qquad C_2\beta^{2/3}\leq \mathrm{Var}\big(\mathfrak{h}_{0\uparrow (1+\beta)t}- \mathfrak{h}_{0\uparrow t}\big)\leq C_1\beta^{2/3} \, . \label{eq:VarUpLowBd2}
\end{align}
By \eqref{eq:VarUpLowBd2} and the Cauchy-Schwarz inequality, there exists $ C^{\prime}=C^{\prime}(t_0)>1$ such that
 \begin{align}
 \rho:= \beta^{-1/3}
 \frac{\mathrm{Cov}\big(\mathfrak{h}_{0\uparrow (1+\beta)t}- \mathfrak{h}_{0\uparrow t}, \mathfrak{h}_{0\uparrow t}(0)\big)}{\mathrm{Var}(\mathfrak{h}_{0\uparrow t}(0))}  \leq C^{\prime} \, .
 \end{align}
 Fix $\delta \in (0,1)$ and denote $\theta := (\beta t)^{2/3}s^{-1}$. Define $\widetilde{\mathcal{W}} = \widetilde{\mathcal{C}}\cap \widetilde{\mathcal{D}}\cap \widetilde{\mathcal{E}} \cap \widetilde{\mathcal{F}}$, where
\begin{align}
\widetilde{\mathcal{C}} &:= \bigcap_{\epsilon\in \{-1,1\}}\Big\{\mathfrak{h}_{(1+t)\beta\downarrow t}(\epsilon t^{-2/3}\theta)\geq 3\beta^{1/3}s\Big\} \, ,\\
\widetilde{\mathcal{D}} &:=\Big\{\inf_{|y|\leq \theta} \mathfrak{h}_{(1+\beta)t\downarrow t}(t^{-2/3} y)\geq (2+\delta)\beta^{1/3}s\Big\} \, , \\
\widetilde{\mathcal{E}}&:= \Big\{\inf_{|y|\leq \theta}\mathfrak{h}_{0\uparrow t}(t^{-2/3} y) - \mathfrak{h}_{0\uparrow t}(0)\geq - \beta^{1/3}\frac{s}{2}\Big\} \, ,
&&\widetilde{\mathcal{F}}:= \Big\{\mathfrak{h}_{0\uparrow t}(0)\leq  \beta^{1/3}\frac{ s}{2 \rho}\Big\} \, .
\end{align}
On the event $\widetilde{\mathcal{W}}$, it follows that, for $|y|\leq \theta$,
\begin{align}
 & t^{1/3}\mathfrak{h}_{(1+\beta)t\downarrow t}( t^{-2/3}y)+ t^{1/3}\mathfrak{h}_{0\uparrow t}(t^{-2/3} y) \nonumber \\
 \geq & \, (\beta t)^{1/3}\big(\tfrac{3}{2}+\delta\big)s +t^{1/3} \mathfrak{h}_{0\uparrow t}(0)
 \, \geq \, (\beta t)^{1/3}(1+\delta)s+ (1+\rho\beta^{1/3})t^{1/3} \mathfrak{h}_{0\uparrow t}(0) \, , \label{eq:Steps}
\end{align}
where the first bound is due to $\mathfrak{h}_{(1+\beta)t\downarrow t}(t^{-2/3}y)$ and $\mathfrak{h}_{0\uparrow t}(t^{-2/3} y)$ being respectively greater than $\beta^{1/3}(2+\delta)s$ and $\mathfrak{h}_{0\uparrow t}(0)- \beta^{1/3}s/2$; while the second depends on $\mathfrak{h}_{0\uparrow t}(0)\leq \beta^{1/3}\frac{s}{2\rho}$.

Substituting~\eqref{eq:Steps} into the composition law Proposition \ref{ppn:TimeEvl} yields
\begin{align}
\exp\big(t^{1/3} \mathfrak{h}_{0\uparrow (1+\beta)t}\big)&\geq \int^{\theta}_{-\theta} \exp\Big((\beta t)^{1/3}(1+\delta)s+ (1+\rho\beta^{1/3})t^{1/3} \mathfrak{h}_{0\uparrow t}(0) \Big) dx \\
&\geq 2 \exp\Big( (1+\delta)(\beta t)^{1/3}s+ (1+\rho\beta^{1/3})t^{1/3}\mathfrak{h}_{0\uparrow t}(0) +\tfrac{2}{3}\log (\beta t) - \log s\Big) \, .
\end{align}
Since we have assumed that $\beta t \geq t_0$, there exists $s_0=s_0(t_0)>0$ such that the quantity in the preceding line is bounded below by $\exp\big((\beta t)^{1/3}s +(1+\rho \beta^{1/3})t^{1/3}\mathfrak{h}_{0\uparrow t}(0) \big)$  for $s>s_0$. Thus, we see that the occurrence of $\widetilde{\mathcal{W}}$ entails that
$\big\{\mathfrak{h}_{0\uparrow (1+\beta)t}- \mathfrak{h}_{0\uparrow t} - \rho \mathfrak{h}_{0\uparrow t}(0)\geq (\beta t)^{1/3}s\big\}$,
 which is nothing other than the event in \eqref{eq:ShiftDiff} whose probability we seek to bound below. The proof of Proposition~\ref{sec:VarAux}  has thus been reduced to the derivation of a suitable lower bound on $\mathbb{P}(\widetilde{\mathcal{W}})$.

Owing to the independence between $\mathfrak{h}_{(1+\beta)t\downarrow t}$ and $\mathfrak{h}_{0\uparrow t}$, $ \mathbb{P}(\widetilde{\mathcal{W}}) =  \mathbb{P}(\widetilde{\mathcal{C}}\cap \widetilde{\mathcal{D}}) \, \mathbb{P}(\widetilde{\mathcal{E}}\cap \widetilde{\mathcal{F}})$.


First let us bound $\mathbb{P}(\widetilde{\mathcal{C}}\cap \widetilde{\mathcal{D}})\geq \exp\big(-c s^{3/2}\big)$, where $c=c(t_0)>0$ and $s>s_0$ for some $s_0(t_0)>0$. Combining the distributional equality~\eqref{eq:DistEq2}, the stationarity in Proposition \ref{ppn:Stationarity} and the lower bound in Proposition \ref{ppn:onepointuppertail}, we find that $\mathbb{P}(\widetilde{\mathcal{C}}) \geq \exp\big(-cs^{3/2}\big)$.
By decomposing $\mathbb{P}(\widetilde{\mathcal{C}}\cap \widetilde{\mathcal{D}}) = \mathbb{P}(\widetilde{\mathcal{C}})- \mathbb{P}(\widetilde{\mathcal{C}}\cap \neg \widetilde{\mathcal{D}})$ we see that our bound on  $\mathbb{P}(\widetilde{\mathcal{C}}\cap \widetilde{\mathcal{D}})$ will be established if we can prove an upper bound $\mathbb{P}(\mathcal{C}\cap \neg \mathcal{D} )\leq   \exp\big(-cs^{\frac{3}{2}+\eta}\big)$ for any $\eta>0$ fixed. This can be established for $\eta=1$ in  precisely the manner of~\eqref{eq.cdneg}.

To obtain a lower bound on $\mathbb{P}\big(\widetilde{\mathcal{W}} \big)$ -- and hence to finish the proof of  Proposition~\ref{cl:VarAux} --  it suffices to show that $\mathbb{P}(\widetilde{\mathcal{E}}\cap \widetilde{\mathcal{F}}) \geq \exp\big(-cs^{3/2}\big)$. Observe that
$\mathbb{P}(\widetilde{\mathcal{E}}\cap \widetilde{\mathcal{F}}) \geq  \mathbb{P}(\widetilde{\mathcal{F}}) - \mathbb{P}(\neg\widetilde{\mathcal{E}})$. From the tail bounds in Propositions \ref{ppn:onepointlowertail} and \ref{ppn:onepointuppertail}, there exists some constant $C>0$ such that $\mathbb{P}(\widetilde{\mathcal{F}})>C$.
The event~$\widetilde{\mathcal{E}}$  may be written in the form $\big\{\inf_{|y|\leq \epsilon}\mathfrak{h}_{0\uparrow t}(t^{-2/3} y)-\mathfrak{h}_{0\uparrow t}(0)\geq -\epsilon^{1/2}\frac{s^{3/2}}{2}\big\}$ by setting $\epsilon = \beta^{2/3}s^{-1}$; an application of Theorem \ref{t.locreg} then implies
that $s_0=s_0(t_0)>0$ and $c=c(t_0)>0$ exist so that $\mathbb{P}\big(\widetilde{\mathcal{E}}\big) \geq 1-\exp\big(-c s^{9/4}\big)$, provided that $s>s_0$ and $\beta t \geq t_0$.
Thus we find that, for $s$ large enough, $\mathbb{P}(\widetilde{\mathcal{E}}\cap \widetilde{\mathcal{F}}) \geq C - \exp\big(-c s^{9/4}\big) \geq \exp\big(-c s^{3/2}\big)$; and thus is the proof of Proposition~\ref{cl:VarAux} concluded. \qed

\section{Spatial modulus of continuity}\label{sec:Holderspace}
\begin{proof}[Proof of Theorem~\ref{thm:SpaceHolder}]
The cases of small and large $|x_1-x_2|$ will separately occupy our attention.
The shorthand  $\widetilde{\mathfrak{h}}_t(x):= \mathfrak{h}_t(x) + \tfrac{x^2}{2}$ will be employed to prevent repetitious parabolic shifting.
Define
\begin{align}
\mathcal{C}_{\ll} &:=\sup_{\substack{x_1\neq x_2\in [a,b]\\|x_1-x_2|<1}} \, \bigg\{ \, |x_1-x_2|^{-1/2}\Big( \log\frac{|b-a|}{|x_1-x_2|}\Big)^{-2/3}\big|\widetilde{\mathfrak{h}}_t(x_1) - \widetilde{\mathfrak{h}}_t(x_2) \big| \, \bigg\} \, , \\
\mathcal{C}_{\gg} &:=\sup_{\substack{x_1\neq x_2\in [a,b]\\ |x_1-x_2|>1/2}}\, \bigg\{ \, |x_1-x_2|^{-1/2}\Big( \log\frac{|b-a|}{|x_1-x_2|}\Big)^{-2/3}\big|\widetilde{\mathfrak{h}}_t(x_1) - \widetilde{\mathfrak{h}}_t(x_2) \big| \, \bigg\} \, ;
\end{align}
and note that $\mathcal{C}$ from \eqref{eq:sholder} is given by $\mathcal{C} = \max(\mathcal{C}_{\ll},\mathcal{C}_{\gg})$. Upper-tail bounds for $\mathcal{C}_{\ll}$ and $\mathcal{C}_{\gg}$ will thus extend to $\mathcal{C}$. Some notation will aid in the derivation of such bounds. Let $n_0 = \inf\big\{n\geq 1: |b-a|< 2^{n-1}\big\}$. For $n\geq n_0$, dyadically partition $[a,b]$ as
\begin{align}\label{eq:dyadic}
[a,b]=\bigcup^{2^n}_{k=1}\mathcal{J}^{(n)}_{k} \qquad \textrm{s.t.} \quad\mathcal{J}^{(n)}_k:= \big[\alpha^{(n)}_{k-1}, \alpha^{(n)}_{k}\big], \quad \alpha^{(n)}_{k}:=a+\frac{k}{2^n}(b-a),  \text{ for }k =0, \ldots , 2^n \, .
\end{align}

First we bound $\mathcal{C}_{\gg}$. Unless $|b-a|>1/2$, there is nothing to prove. Supposing then the latter bound, consider the dyadic partition with $n=n_0$. Each interval in the partition has length less than $1/2$; thus, $x_1$ and $x_2$ lie in distinct intervals. Label by $k_1$ and $k_2$ the respective indices~$k$ of the right endpoint~$\alpha^{(n_0)}_k$ of the interval containing $x_1$ and $x_2$.  By the triangle inequality,
\begin{align}
\big|&\widetilde{\mathfrak{h}}_t(x_1)-\widetilde{\mathfrak{h}}_t(x_2)\big|\\&\,\,\leq 4 \max\Big(\big|\widetilde{\mathfrak{h}}_t(x_1)-\widetilde{\mathfrak{h}}_t(\alpha^{(n_0)}_{k_1})\big|,
\big|\widetilde{\mathfrak{h}}_t(\alpha^{(n_0)}_{k_1})\big|,
\big|\widetilde{\mathfrak{h}}_t(\alpha^{(n_0)}_{k_2})\big|,
\big|\widetilde{\mathfrak{h}}_t(x_1)-\widetilde{\mathfrak{h}}_t(\alpha^{(n_0)}_{k_2})\big|\Big).
\end{align}
Introducing $B^{(n)}_k:= \sup_{x\in \mathcal{J}^{(n)}_{k}}\big|\widetilde{\mathfrak{h}}_t(x)-\widetilde{\mathfrak{h}}_t(\alpha^{(n)}_{k})\big|$ and $C^{(n)}_k := \big|\widetilde{\mathfrak{h}}_t(\alpha^{(n_0)}_{k})\big|$, we learn that
\begin{align}
\mathcal{C}_{\gg} \leq (1/2)^{-1/2}\Big(\log \frac{|b-a|}{1/2}\Big)^{-2/3}  4 \max_{k\in \{0,\ldots, 2^{n_0}\}}\big\{B^{(n_0)}_{k},C^{(n_0)}_{k}\big\}.
\end{align}
In total, the maximum is taken over $2^{n_0}+1$ possible values of $k$. By Theorem \ref{t.locreg}, we see that $\mathbb{P}\big(B^{(n_0)}_k\geq s\big)\leq \exp\big(-cs^{3/2}\big)$ for some $c=c(t_0)>0$, provided that $s>s_0$ for $s_0=s_0(t_0)$. Proposition \ref{ppn:onepointuppertail} provides a similar bound on the upper tails of the $C^{(n_0)}_k$. A union bound transmits this inference to an upper-tail bound on $\mathcal{C}_{\gg}$, with dependences of the form $c=c(t_0,|b-a|)$ and $s_0=s_0(t_0,|b-a|)$. The dependence on $|b-a|$ arises by the absorption into $c$ of the term $(1/2)^{-\frac{1}{2}}\Big(\log \frac{|b-a|}{1/2}\Big)^{-2/3}$ and of the term $2(2^{n_0}+1)$ arising from the union bound.

As we turn to bounding above $\mathcal{C}_{\ll}$, we may impose that $|x_1-x_2|<1$. Let $n\geq n_0$ be the smallest integer such that $|x_1-x_2|\geq |b-a| 2^{-n-1}$. Either $x_1$ and $x_2$ lie in a given interval in the dyadic partition of level $n$, or they lie in consecutive intervals. When $x_1,x_2\in \mathcal{J}^{(n)}_k$ for some $k\in \{0,\ldots, 2^{n}\}$, recall of the event $B^{(n)}_k$ and use of the triangle inequality yield
\begin{align}
&|x_1-x_2|^{-1/2}\bigg( \log\frac{|b-a|}{|x_1-x_2|}\bigg)^{-2/3}\big|\widetilde{\mathfrak{h}}_t(x_1) - \widetilde{\mathfrak{h}}_t(x_2) \big|\\\leq &\,
\big(|b-a| 2^{-n-1}\big)^{-1/2} \Big( \log(2^{n+1})\Big)^{-2/3} 2 B^{(n)}_k \, .
\end{align}
When $x_1\in \mathcal{J}^{(n)}_{k}$ and $x_2\in \mathcal{J}^{(n)}_{k+1}$, we may
set $\widetilde{B}^{(n)}_k:= \sup_{x\in \mathcal{J}^{(n)}_{k}}\big|\widetilde{\mathfrak{h}}_t(x)-\widetilde{\mathfrak{h}}_t(\alpha^{(n)}_{k-1})\big|$ to obtain
\begin{align*}
 & |x_1-x_2|^{-1/2}\bigg( \log\frac{|b-a|}{|x_1-x_2|}\bigg)^{-2/3}\big|\widetilde{\mathfrak{h}}_t(x_1) - \widetilde{\mathfrak{h}}_t(x_2) \big| \\
   \leq & \,
\big(|b-a| 2^{-n-1}\big)^{-1/2} \Big( \log(2^{n+1})\Big)^{-2/3} 2 \max\big(B^{(n)}_k,\widetilde{B}^{(n)}_k\big) \, .
\end{align*}
From these bounds, we arrive at
\begin{align}
\mathcal{C}_{\ll} \, \leq \, \sup_{n\geq n_0} \sup_{k\in \{0,\ldots, 2^n\}} \big(|b-a| 2^{-n-1}\big)^{-1/2} \Big( \log(2^{n+1})\Big)^{-2/3} 2 \max\big(B^{(n)}_k,\widetilde{B}^{(n)}_k\big)  \, .
\end{align}
In contrast to the analysis of $\mathcal{C}_{\gg}$, infinitely many terms must now be considered. Indeed, we have
\begin{align}
\big\{\mathcal{C}_{\ll}\geq s\big\} = \bigcup_{n=n_0}^{\infty}\bigcup_{k=0}^{ 2^n}\, \bigg\{ \, \max\big(B^{(n)}_k,\widetilde{B}^{(n)}_k\big) \geq \big(|b-a| 2^{-n-1}\big)^{1/2} \Big( \log(2^{n+1})\Big)^{2/3} 2^{-1} s \, \bigg\} \, .
\end{align}
Apply the union bound and the tail bound in Theorem \ref{t.locreg} with $\epsilon = 2^{-n}|b-a|$ -- which parameter is at most one since $n\geq n_0$. We thus find that there exist $s_0=s_0(t_0,|b-a|)>0$ and $c=c(t_0,|b-a|)>0$ such that, for $t>t_0$ and $s>s_0$,
\begin{align}\label{eq:TailOfA}
\mathbb{P}\big(\mathcal{C}_{\ll}\geq s\big)\leq \sum_{n=n_0}^{\infty}\sum_{k=0}^{ 2^n} \exp\big(- n c s^{3/2}\big)\leq \sum_{n=n_0}^{\infty} \exp\big(-n cs^{3/2}\big)\leq \exp\big(-c s^{3/2}\big) \, ;
\end{align}
here, the values of $c$ and $s_0$ change between each inequality in order to absorb the higher indexed terms in the two sums. Indeed, in the second bound, the sum over $k$ contributes a factor $2^{n}+1$ which is absorbed by decreasing the constant $c$ provided that $s_0$ is high enough. The third bound arises by computing the geometric sum, expressing $n_0$ in terms of $|b-a|$, and absorbing the resulting constant into~$c$.

These bounds $\mathbb{P}\big(\mathcal{C}_{\ll}\geq s\big)$ and $\mathbb{P}\big(\mathcal{C}_{\gg}\geq s\big)$ obtained, the proof is Theorem~\ref{thm:SpaceHolder} is complete.
\end{proof}

\section{Appendix: tail probabilities and covariances}\label{sec:Appendix}

Several tools are stated and proved here. It is tempting to rewrite in \eqref{eq:Monotone} and \eqref{eq:Monotone1} in terms of conditional probabilities.
The latter may, however, not exist, so we use a different formulation.
\bl\label{lem:Aux}
Let $X$ and $Y$ be two real-valued integrable random variables and suppose that, for  $r \in \RR$ and $u>v\in \RR$,
\begin{align}
\mathbb{P}(Y> v)\mathbb{P}(X> r, Y> u)\geq \mathbb{P}(Y > u)\mathbb{P}(X> r, Y> v) \, , \label{eq:Monotone}\\
\mathbb{P}(Y\leq u)\mathbb{P}(X>  r, Y\leq  v)\geq \mathbb{P}(Y \leq v)\mathbb{P}(X> r, Y\leq u) \, .\label{eq:Monotone1}
\end{align}
\begin{enumerate}
\item Then $\mathrm{Cov}(X,Y)\geq 0$.
\item Moreover, for any $a\in \RR$,
\begin{align}\label{eq:CovBd1}
\mathrm{Cov}(X, Y)\geq \mathbb{P}(Y\geq a)\cdot\Big(\mathbb{E}\big[X|Y\geq a\big] - \mathbb{E}[X]\Big)\cdot \Big(\mathbb{E}\big[Y|Y> a\big] - \mathbb{E}\big[Y|Y\leq a\big]\Big).
\end{align}
\end{enumerate}
\el

\begin{proof}
\noindent{\bf $(1)$:} Since \eqref{eq:Monotone} holds after replacing $X$ by $X-E[X]$ and $Y$ by $Y-\mathbb{E}[Y]$, we may assume that $E[X]= E[Y]=0$. Denote $X_{+}:= \max\{X,0\}$ and $X_{-}:= \max\{-X,0\}$, so that $X= X_{+}-X_{-}$; and likewise for $Y$. Thus,
$\mathbb{E}[X]= \mathbb{E}[X_{+}]- \mathbb{E}[X_{-}]$; $\mathbb{E}[Y]= \mathbb{E}[Y_{+}]- \mathbb{E}[Y_{-}]$; and
\begin{align}
\mathrm{Cov}(X,Y) = \mathbb{E}[X_{+}Y_{+}]- &\mathbb{E}[X_{+}Y_{-}] -\mathbb{E}[X_{-}Y_{+}] + \mathbb{E}[X_{-}Y_{-}] \, .
\end{align}
To prove that $\mathrm{Cov}(X,Y)\geq 0$, it suffices to show that
\begin{align}\mathbb{E}[X_{+}Y_{+}]\geq \mathbb{E}[X_{+}]\mathbb{E}[Y_{+}] \, , &\qquad \mathbb{E}[X_{+}Y_{-}]\leq \mathbb{E}[X_{+}]\mathbb{E}[Y_{-}] \, ,\label{eq:1stIneq}\\
\mathbb{E}[X_{-}Y_{+}]\leq \mathbb{E}[X_{-}]\mathbb{E}[Y_{+}] \, , &\qquad \mathbb{E}[X_{-}Y_{-}]\geq \mathbb{E}[X_{-}]\mathbb{E}[Y_{-}] \, . \label{eq:2ndIneq}
\end{align}
We prove \eqref{eq:1stIneq}; the bounds in the following line are derived similarly.
Taking $v\to -\infty$ in \eqref{eq:Monotone} yields
\begin{align}\label{eq:1stLine}
\mathbb{P}(X> r, Y> u)\geq \mathbb{P}(X>r)\mathbb{P}(Y>u), \quad \forall\, r\in \RR, u\in \RR.
\end{align}
Subtracting \eqref{eq:1stLine} from $\mathbb{P}(X>r)$ yields
\begin{align}\label{eq:2ndLine}
\mathbb{P}(X> r, Y\leq u) \leq \mathbb{P}(X> r)\mathbb{P}(Y\leq u), \quad \forall r\in \RR, u \in \RR.
\end{align}
Integrating \eqref{eq:1stLine} with respect to $(r,u)$ over $(0,\infty)\times (0,\infty)$, we see that
\begin{align}
\int_{0}^{\infty}\!\!\int^{\infty}_{0}\!\!\mathbb{P}(X> r, Y> u) dr du= \mathbb{E}[X_{+}Y_{+}]\geq \int_{0}^{\infty}\!\!\int^{\infty}_{0}\!\! \mathbb{P}(X> r)\mathbb{P}(Y\leq u)dr du= \mathbb{E}[X_{+}]\mathbb{E}[Y_{+}] \, .\qquad\label{eq:Int1}
\end{align}
 Integrating \eqref{eq:2ndLine} with respect to $(r,u)$ over $(0,\infty)\times (-\infty,0]$ yields
 \begin{align}
 \int_{0}^{\infty}\!\!\int^{0}_{-\infty}\!\!\mathbb{P}(X> r, Y\leq u) dr du= \mathbb{E}[X_{+}Y_{+}]\leq \int_{0}^{\infty}\!\!\int^{0}_{-\infty}\!\! \mathbb{P}(X> r)\mathbb{P}(Y\leq u)dr du= \mathbb{E}[X_{+}]\mathbb{E}[Y_{-}] \, .\qquad \label{eq:Int2}
 \end{align}
 The bounds in line~\eqref{eq:1stIneq} follow from  \eqref{eq:Int1} and \eqref{eq:Int2}.
 \smallskip

\noindent{\bf $(2)$:}
Fix $a\in \RR$. If $\mathbb{P}(Y\geq a) \in \{ 0 ,1 \}$, the right-hand side of \eqref{eq:CovBd1} equals zero; in which case, \eqref{eq:CovBd1} follows from the just derived part $(1)$. Suppose then that $\mathbb{P}(Y\geq a) \in (0,1)$. Let the pair $(X^{\prime},Y^{\prime})$ be an independent copy of $(X,Y)$. One may write
\begin{align}\label{eq:ProdBd}
\mathbb{E}[X(Y-a)] \, = \, \frac{\mathbb{E}\big[ \mathbf{1}_{Y^{\prime}> a} X \,(Y-a)_+ \big]}{\mathbb{P}(Y> a)} - \frac{\mathbb{E}\big[\mathbf{1}_{Y^{\prime}\leq a} X\, (Y-a)_{-}\big]}{\mathbb{P}(Y\leq a)} \, .
\end{align}

We will next argue that
\begin{align}
\mathbb{E}\big[ \mathbf{1}_{Y^{\prime}> a} X (Y-a)_+ \big]&\geq \mathbb{E}\big[X^{\prime}\mathbf{1}_{Y^{\prime}> a}\big]\mathbb{E}\big[(Y-a)_+\big] \, , \label{eq:MyEq1}\\
\mathbb{E}\big[\mathbf{1}_{Y^{\prime}\leq a} X (Y-a)_{-}\big]&\leq \mathbb{E}\big[X^{\prime}\mathbf{1}_{Y^{\prime}\leq a}\big]\mathbb{E}\big[(Y-a)_{-}\big] \, . \label{eq:MyEq2}
\end{align}
In fact, we will merely show how \eqref{eq:MyEq1} follows from \eqref{eq:Monotone};  \eqref{eq:MyEq2} follows in a similar way from \eqref{eq:Monotone1}. By ~\eqref{eq:Monotone}, we see that, for $r\in \RR$ and $u >0$,
\begin{align}
\mathbb{P}(Y^{\prime}>a, X> r, Y>a+u)\geq \mathbb{P}(Y^{\prime}>a,X^{\prime}> r)\mathbb{P}(Y>a+u) \, .
\end{align}
Integrating this inequality with respect to $(r,u)$ over $[0,\infty)^2$, we obtain
\begin{align}\label{eq:Grt1}
\mathbb{E}\big[\mathbf{1}_{Y^{\prime}> a}X_+\,(Y-a)_+\big] \, \geq \, \mathbb{E}\big[\mathbf{1}_{Y^{\prime}> a}X^{\prime}_+\big]\, \mathbb{E}\big[(Y-a)_+\big] \, .
\end{align}
Subtracting \eqref{eq:Monotone} from $\mathbb{P}(Y>u)\mathbb{P}(Y>v)$ with $u=a$ yields
\begin{align}
\mathbb{P}(Y^{\prime}>a, X\leq r,  Y>u)\leq \mathbb{P}(Y^{\prime}>a,X^{\prime}\leq r )\, \mathbb{P}(Y>u) \, ;
\end{align}
which, after integrating with respect to $(r,u)$ over $(- \infty,0)\times [a,\infty)$, gives
\begin{align}\label{eq:Leq1}
\mathbb{E}\big[\mathbf{1}_{Y^{\prime}> a} X_{-}(Y-a)_{+}\big]\leq \mathbb{E}\big[\mathbf{1}_{Y^{\prime}> a} X^{\prime}_{-}\big]\, \mathbb{E}\big[(Y-a)_+\big] \, .
\end{align}
Subtracting  \eqref{eq:Leq1} from  \eqref{eq:Grt1} yields \eqref{eq:MyEq1}.

Pursuing the proof of \eqref{eq:CovBd1}, we substitute \eqref{eq:MyEq1} and \eqref{eq:MyEq2} into the right-hand side of \eqref{eq:ProdBd} to learn that
\begin{align}
\mathbb{E}\big[X(Y-a)\big]\geq \mathbb{E}\big[X\mathbf{1}_{Y> a}\big] \Big(\mathbb{E}\big[Y|Y> a\big] -a\Big) + \mathbb{E}\big[X\mathbf{1}_{Y\leq a}\big] \Big(\mathbb{E}\big[Y|Y\leq a\big] -a\Big) \, .
\end{align}
 Subtracting $\mathbb{E}[X]\mathbb{E}[Y-a]$ from this inequality, and simplifying, yields
 \begin{align}
 \mathrm{Cov}( X,Y)\geq & \, \, \Big(\mathbb{E}\big[X\mathbf{1}_{Y>a}\big] - \mathbb{E}[X]\,\mathbb{P}(Y> a)\Big)\Big(\mathbb{E}\big[Y|Y> a\big] -a\Big) \\
 & \, \, + \, \Big(\mathbb{E}\big[X\mathbf{1}_{Y\leq a}\big] - \mathbb{E}[X]\mathbb{P}(Y\leq a)\Big)\Big(\mathbb{E}\big[Y|Y\leq a\big] -a\Big) \, . \label{eq:Simpl}
 \end{align}
Now, we obtain \eqref{eq:CovBd1}: it follows from \eqref{eq:Simpl} by noting that
$$ \mathbb{E}\big[X\mathbf{1}_{Y\leq a}\big] - \mathbb{E}[X]\,\mathbb{P}(Y\leq a)=-\Big( \mathbb{E}\big[X\mathbf{1}_{Y> a}\big] - \mathbb{E}[X]\,\mathbb{P}(Y> a)\Big) $$
and rewriting $\mathbb{E}[X\mathbf{1}_{Y> a}]= \mathbb{P}(Y> a)\mathbb{E}[X|{Y> a}]$.
\end{proof}

\bc\label{cor:CovLow}
Let $X$ and $Y$ be integrable random variables that satisfy \eqref{eq:Monotone} and \eqref{eq:Monotone1}. Suppose that there exist $C_1, C_2>0$ such that
\begin{align}\label{eq:CondExBd}
\mathbb{E}\big[X\big|Y \geq C_1\big] \geq \mathbb{E}[X] +C_2 \, .
\end{align}
 Then
\begin{align}\label{eq:CovBd2}
\mathrm{Cov}(X,Y) \, \geq \,  C_2\cdot \mathbb{P}(Y\geq C_1) \Big(\mathbb{E}[Y|Y\geq C_1] - \mathbb{E}[Y|Y\leq C_1]\Big) \, .
\end{align}
\ec
\begin{proof}
We obtain \eqref{eq:CovBd2} from \eqref{eq:CovBd1} by taking $a=C_1$ and applying \eqref{eq:CondExBd}.
\end{proof}

\bl\label{lem:GenCorrBd}
Let $X$, $Y$ and $Z$ be  non-degenerate real-valued random variables with finite second moments such that $Z= X+Y$. Define
\begin{align}
\chi := \frac{\mathrm{Var}(X)}{\mathrm{Var}(Y)} \, , \qquad \Psi:= \frac{\mathrm{Cov}(X,Y)}{\mathrm{Var}(Y)} \, ,\qquad \Theta:= \frac{\mathrm{Var}(X)\mathrm{Var}(Y)- (\mathrm{Cov}(X,Y))^2}{(\mathrm{Var}(Y))^2} \, .
\end{align}
Suppose that $\max\{\chi,|2\Psi+\chi|\}<1$. There exist two constants $C_2 > C_1 > 0$ such that
\begin{align}\label{eq:UpLowCorrBd}
1- \Theta/2 +C_1 \chi^{3/2} \leq \mathrm{Corr}(Z,Y) \leq 1- \Theta/2 +C_2 \chi^{3/2} \, .
\end{align}
\el

\begin{proof}  Since $\mathrm{Var}(Z)=(1+2\Psi+\chi)\mathrm{Var}(Y)$, we may rewrite
\begin{align}\label{eq:CorDefRe}
\mathrm{Corr}\big(Z,Y\big)& = \frac{\mathrm{Cov}\big(Z,Y\big)}{\sqrt{\mathrm{Var}(Z)}\sqrt{\mathrm{Var}(Y)}} = \frac{1+ \Psi}{\sqrt{1+2\Psi+\chi}} \, .
\end{align}
Taylor expanding $(1+2\Psi +\chi)^{-1/2}$ with respect to $\Xi:=\Psi+\frac{1}{2}\chi$ yields
\begin{align}
 \mathrm{Corr}\big(Z,Y\big)= \big(1+ \Psi\big) \Big(1- \Xi + \tfrac{3}{2}\Xi^2 + \mathfrak{R}\Big) \quad \textrm{where}\quad |\mathfrak{R}|\leq C |\Xi|^3 \label{eq:CorrExpan}
\end{align}
for some constant $C>0$.  Simplifying the product and substituting $\Xi= \Psi+\frac{1}{2}\chi$, we find that
 \begin{align}\label{eq:ExpPsi}
 \mathrm{Corr}\big(Z,Y\big)= 1-\frac{1}{2}\chi + \frac{1}{2}\Psi^2+\frac{1}{2}\Psi\chi+ \frac{3}{8}\chi^2 + \frac{3}{2}\Psi\big(\Psi +\frac{1}{2}\chi\big)^2 + \mathfrak{R}^{\prime} \, , \quad \textrm{where} \quad \vert \mathfrak{R}^{\prime} \vert \leq C^{\prime}|\Xi|^3  \qquad
 \end{align}
 for some constant $C^{\prime}>0$.  By the Cauchy-Schwarz inequality, $\Psi^2\leq \chi$;
note also that $\chi\in [0,1]$. Thus, there exist $C_2 > C_1 > 0$ such that
\begin{align}\label{eq:ResEst}
C_1\chi^{3/2}\leq \frac{1}{2}\Psi\chi+ \frac{3}{8}\chi^2 + \frac{3}{2}\Psi\big(\Psi +\frac{1}{2}\chi\big)^2\leq C_2\chi^{3/2} \, .
\end{align}
Substituting \eqref{eq:ResEst} into \eqref{eq:ExpPsi} and noting that $\Theta = \chi -\Psi^2$ proves \eqref{eq:UpLowCorrBd}.
\end{proof}
\smallskip

\bl\label{ppn:VarBound}
Fix $\theta>0$. Let $X$ be a real-valued random variable.
\begin{enumerate}
\item Suppose there exist $s_0$, $\alpha>0$ and $c>0$ such that, for $s\geq s_0$,
\begin{align}\label{eq:PUpLowBd}
\mathbb{P}(X\leq - s\theta)\leq e^{-cs^{\alpha}} &\qquad \textrm{and}\qquad\mathbb{P}(X\geq s\theta)\leq e^{-cs^{\alpha}} \, .
\end{align}
Then there exists $C=C(s_0, c, \alpha)>0$ such that $\mathrm{Var}(X)\leq C \theta^2$.
\item Suppose that $\big\vert \mathbb{E}[X] \big\vert \leq C_1\theta$ for some $C_1>0$, and that there exist $s_0$, $\alpha>0$ and $c>0$ such that, for $s\geq s_0$,
\begin{align}\label{eq:PLowBd}
\mathbb{P} \big(X \geq  s\theta \big)\geq e^{-cs^{\alpha}} \, .
\end{align}
Then there exists $C=C(C_1, s_0,c, \alpha)>0$ such that $\mathrm{Var}(X)\geq C\theta^2$.
\end{enumerate}
\el
\begin{proof} To prove $(1)$, observe that
\begin{align}\label{eq:Vx}
\mathrm{Var}(X) \leq  \mathbb{E}[X^2]=\theta^2\int^{\infty}_{0}2s\Big(\mathbb{P}(X>s\theta)+ \mathbb{P}(X<-s\theta)\Big) ds \, .
\end{align}
Substituting \eqref{eq:PUpLowBd} into the right-hand side and integrating, we obtain the sought bound on  $\mathrm{Var}(X)$.

To prove $(2)$, suppose first that $\mathbb{E}[X]\leq 0$. Then
\begin{align}
\mathrm{Var}(X)\geq \mathbb{E}\big[(X_+)^2\big] = \theta^2\int^{\infty}_0 2s\mathbb{P}(X>s\theta) d\theta \, . \label{eq:VarLow1}
\end{align}
Substituting \eqref{eq:PLowBd} into this right-hand side, and integrating,  results in the desired bound on $\mathrm{Var}(X)$. Suppose instead that $\mathbb{E}[X]>0$. Since $E[X]\leq C_1\theta$, we have
\begin{align}
\mathrm{Var}(X)\geq \mathbb{E}\big[\big( X-C_1\theta_+\big)^2\big]\geq \theta^2 \int^{\infty}_02s\mathbb{P}\big(X>(C_1+s)\theta\big) ds \, .
\end{align}
Similarly, we now substitute \eqref{eq:PLowBd} into this right-hand side, and integrate, to obtain the   bound on $\mathrm{Var}(X)$ that we seek.
\end{proof}

 \bibliographystyle{imsart-number}
\bibliography{Reference}

\begin{thebibliography}{63}

\bibitem{alberts2014}
\begin{barticle}[author]
\bauthor{\bsnm{Alberts},~\bfnm{T.}\binits{T.}},
  \bauthor{\bsnm{Khanin},~\bfnm{K.}\binits{K.}} \AND
  \bauthor{\bsnm{Quastel},~\bfnm{J.}\binits{J.}}
(\byear{2014}).
\btitle{The intermediate disorder regime for directed polymers in dimension
  $1+1$}.
\bjournal{Ann. Probab.}
\bvolume{42}
\bpages{1212--1256}.
\bdoi{10.1214/13-AOP858}
\end{barticle}
\endbibitem

\bibitem{Amir11}
\begin{barticle}[author]
\bauthor{\bsnm{Amir},~\bfnm{G.}\binits{G.}},
  \bauthor{\bsnm{Corwin},~\bfnm{I.}\binits{I.}} \AND
  \bauthor{\bsnm{Quastel},~\bfnm{J.}\binits{J.}}
(\byear{2011}).
\btitle{Probability distribution of the free energy of the continuum directed
  random polymer in {$1+1$} dimensions}.
\bjournal{Commun. Pure Appl. Math.}
\bvolume{64}
\bpages{466--537}.
\bdoi{10.1002/cpa.20347}
\bmrnumber{2796514}
\end{barticle}
\endbibitem

\bibitem{BZ17}
\begin{barticle}[author]
\bauthor{\bsnm{{Baik}},~\bfnm{J.}\binits{J.}} \AND
  \bauthor{\bsnm{{Liu}},~\bfnm{Z.}\binits{Z.}}
(\byear{2017}).
\btitle{{Multi-point distribution of periodic TASEP}}.
\bjournal{arXiv:1710.03284, to appear in J. Amer. Math. Soc.}
\end{barticle}
\endbibitem

\bibitem{BQS11}
\begin{barticle}[author]
\bauthor{\bsnm{Bal\'azs},~\bfnm{M.}\binits{M.}},
  \bauthor{\bsnm{Quastel},~\bfnm{J.}\binits{J.}} \AND
  \bauthor{\bsnm{Sepp\"{a}l\"{a}inen},~\bfnm{T.}\binits{T.}}
(\byear{2011}).
\btitle{Fluctuation exponent of the {KPZ}/stochastic {B}urgers equation.}
\bjournal{J. Amer. Math. Soc.}
\bvolume{24}
\bpages{683--708}.
\end{barticle}
\endbibitem

\bibitem{BG18}
\begin{barticle}[author]
\bauthor{\bsnm{{Basu}},~\bfnm{R.}\binits{R.}} \AND
  \bauthor{\bsnm{{Ganguly}},~\bfnm{S.}\binits{S.}}
(\byear{2018}).
\btitle{{Time Correlation Exponents in Last Passage Percolation}}.
\bjournal{arXiv:1807.09260}.
\end{barticle}
\endbibitem

\bibitem{HammondSarkarBaus}
\begin{barticle}[author]
\bauthor{\bsnm{{Basu}},~\bfnm{R.}\binits{R.}},
  \bauthor{\bsnm{{Ganguly}},~\bfnm{S.}\binits{S.}} \AND
  \bauthor{\bsnm{{Hammond}},~\bfnm{A.}\binits{A.}}
(\byear{2019}).
\btitle{{Fractal geometry of Airy$_2$ processes coupled via the Airy sheet.}}
\bjournal{arXiv:1904.01717}.
\end{barticle}
\endbibitem

\bibitem{BGZ19}
\begin{barticle}[author]
\bauthor{\bsnm{{Basu}},~\bfnm{R.}\binits{R.}},
  \bauthor{\bsnm{{Ganguly}},~\bfnm{S.}\binits{S.}} \AND
  \bauthor{\bsnm{{Zhang}},~\bfnm{L.}\binits{L.}}
(\byear{2019}).
\btitle{{Temporal Correlation in Last Passage Percolation with Flat Initial
  Condition via Brownian Comparison}}.
\bjournal{arXiv e-prints}
\bpages{arXiv:1912.04891}.
\end{barticle}
\endbibitem

\bibitem{BenArous}
\begin{barticle}[author]
\bauthor{\bsnm{{Ben-Arous}},~\bfnm{G.}\binits{G.}}
(\byear{2003}).
\btitle{{Aging and spin-glass dynamics}}.
\bjournal{arXiv Mathematics e-prints}
\bpages{math/0304364}.
\end{barticle}
\endbibitem

\bibitem{Bertini1997}
\begin{barticle}[author]
\bauthor{\bsnm{Bertini},~\bfnm{L.}\binits{L.}} \AND
  \bauthor{\bsnm{Giacomin},~\bfnm{G.}\binits{G.}}
(\byear{1997}).
\btitle{Stochastic {B}urgers and {KPZ} Equations from Particle Systems}.
\bjournal{Commun. Math. Phys.}
\bvolume{183}
\bpages{571--607}.
\bdoi{10.1007/s002200050044}
\end{barticle}
\endbibitem

\bibitem{BBW16}
\begin{barticle}[author]
\bauthor{\bsnm{{Borodin}},~\bfnm{A.}\binits{A.}},
  \bauthor{\bsnm{{Bufetov}},~\bfnm{A.}\binits{A.}} \AND
  \bauthor{\bsnm{{Wheeler}},~\bfnm{M.}\binits{M.}}
(\byear{2016}).
\btitle{{Between the stochastic six vertex model and {H}all-{L}ittlewood
  processes}}.
\bjournal{ArXiv e-prints}.
\end{barticle}
\endbibitem

\bibitem{caputo2019}
\begin{barticle}[author]
\bauthor{\bsnm{Caputo},~\bfnm{P.}\binits{P.}},
  \bauthor{\bsnm{Ioffe},~\bfnm{D.}\binits{D.}} \AND
  \bauthor{\bsnm{Wachtel},~\bfnm{V.}\binits{V.}}
(\byear{2019}).
\btitle{{Confinement of Brownian polymers under geometric area tilts}}.
\bjournal{Electron. J. Probab.}
\bvolume{24}.
\end{barticle}
\endbibitem

\bibitem{CIW2019}
\begin{barticle}[author]
\bauthor{\bsnm{Caputo},~\bfnm{P.}\binits{P.}},
  \bauthor{\bsnm{Ioffe},~\bfnm{D.}\binits{D.}} \AND
  \bauthor{\bsnm{Wachtel},~\bfnm{V.}\binits{V.}}
(\byear{2019}).
\btitle{{Tightness and Line Ensembles for Brownian Polymers under Geometric
  Area Tilts}}.
\bjournal{arXiv:1906.06533}.
\end{barticle}
\endbibitem

\bibitem{Corwin12}
\begin{barticle}[author]
\bauthor{\bsnm{Corwin},~\bfnm{I.}\binits{I.}}
(\byear{2012}).
\btitle{The {K}ardar-{P}arisi-{Z}hang equation and universality class}.
\bjournal{Rand. Matr. Theo. Appl.}
\bvolume{1}
\bpages{1130001, 76}.
\bdoi{10.1142/S2010326311300014}
\bmrnumber{2930377}
\end{barticle}
\endbibitem

\bibitem{C18}
\begin{bincollection}[author]
\bauthor{\bsnm{Corwin},~\bfnm{I.}\binits{I.}}
(\byear{2018}).
\btitle{Exactly solving the {KPZ} equation}.
In \bbooktitle{Random growth models}.
\bseries{Proc. Sympos. Appl. Math.}
\bvolume{75}
\bpages{203--254}.
\bpublisher{Amer. Math. Soc., Providence, RI}.
\bmrnumber{3838899}
\end{bincollection}
\endbibitem

\bibitem{Corwin2018}
\begin{barticle}[author]
\bauthor{\bsnm{Corwin},~\bfnm{I.}\binits{I.}} \AND
  \bauthor{\bsnm{Dimitrov},~\bfnm{E.}\binits{E.}}
(\byear{2018}).
\btitle{{Transversal Fluctuations of the ASEP, Stochastic Six Vertex Model, and
  Hall-Littlewood Gibbsian Line Ensembles}}.
\bjournal{Commun. Math. Phys.}
\bvolume{363}
\bpages{435--501}.
\end{barticle}
\endbibitem

\bibitem{corwin2012}
\begin{barticle}[author]
\bauthor{\bsnm{Corwin},~\bfnm{I.}\binits{I.}},
  \bauthor{\bsnm{Ferrari},~\bfnm{P.~L.}\binits{P.~L.}} \AND
  \bauthor{\bsnm{P\'{e}ch\'{e}},~\bfnm{S.}\binits{S.}}
(\byear{2012}).
\btitle{Universality of slow decorrelation in {KPZ} growth}.
\bjournal{Ann. Inst. H. Poincar\'{e} B}
\bvolume{48}
\bpages{134--150}.
\bdoi{10.1214/11-AIHP440}
\end{barticle}
\endbibitem

\bibitem{CG18a}
\begin{barticle}[author]
\bauthor{\bsnm{{Corwin}},~\bfnm{I.}\binits{I.}} \AND
  \bauthor{\bsnm{{Ghosal}},~\bfnm{P.}\binits{P.}}
(\byear{2020}).
\btitle{{ Lower tail of the KPZ equation}}.
\bjournal{Duke Math. J.}
\bvolume{169}
\bpages{1329--1395}.
\bdoi{10.1215/00127094-2019-0079}
\bmrnumber{4094738}
\end{barticle}
\endbibitem

\bibitem{CG18b}
\begin{barticle}[author]
\bauthor{\bsnm{Corwin},~\bfnm{I.}\binits{I.}} \AND
  \bauthor{\bsnm{Ghosal},~\bfnm{P.}\binits{P.}}
(\byear{2020}).
\btitle{{KPZ} equation tails for general initial data}.
\bjournal{Electron. J. Probab.}
\bvolume{25}
\bpages{38 pp.}
\bdoi{10.1214/20-EJP467}
\end{barticle}
\endbibitem

\bibitem{CH14}
\begin{barticle}[author]
\bauthor{\bsnm{Corwin},~\bfnm{I.}\binits{I.}} \AND
  \bauthor{\bsnm{Hammond},~\bfnm{A.}\binits{A.}}
(\byear{2014}).
\btitle{Brownian {G}ibbs property for {A}iry line ensembles}.
\bjournal{Invent. Math.}
\bvolume{195}.
\end{barticle}
\endbibitem

\bibitem{CorHam16}
\begin{barticle}[author]
\bauthor{\bsnm{Corwin},~\bfnm{I.}\binits{I.}} \AND
  \bauthor{\bsnm{Hammond},~\bfnm{A.}\binits{A.}}
(\byear{2016}).
\btitle{K{PZ} line ensemble}.
\bjournal{Probab. Theo. Rel. Fields}
\bvolume{166}
\bpages{67--185}.
\bdoi{10.1007/s00440-015-0651-7}
\bmrnumber{3547737}
\end{barticle}
\endbibitem

\bibitem{COSZ}
\begin{barticle}[author]
\bauthor{\bsnm{Corwin},~\bfnm{I.}\binits{I.}},
  \bauthor{\bsnm{O'Connell},~\bfnm{N.}\binits{N.}},
  \bauthor{\bsnm{Sepp\"al\"ainen},~\bfnm{T.}\binits{T.}} \AND
  \bauthor{\bsnm{Zygouras},~\bfnm{N.}\binits{N.}}
(\byear{2014}).
\btitle{Tropical combinatorics and Whittaker functions}.
\bjournal{Duke Math. J.}
\bvolume{163}
\bpages{513--563}.
\end{barticle}
\endbibitem

\bibitem{CQ11}
\begin{barticle}[author]
\bauthor{\bsnm{Corwin},~\bfnm{I.}\binits{I.}} \AND
  \bauthor{\bsnm{Quastel},~\bfnm{J.}\binits{J.}}
(\byear{2013}).
\btitle{Crossover distributions at the edge of the rarefaction fan}.
\bjournal{Ann. Probab.}
\bvolume{41}.
\end{barticle}
\endbibitem

\bibitem{Corwin2015}
\begin{barticle}[author]
\bauthor{\bsnm{Corwin},~\bfnm{I.}\binits{I.}},
  \bauthor{\bsnm{Quastel},~\bfnm{J.}\binits{J.}} \AND
  \bauthor{\bsnm{Remenik},~\bfnm{D.}\binits{D.}}
(\byear{2015}).
\btitle{{Renormalization Fixed Point of the KPZ Universality Class}}.
\bjournal{J. Stat. Phys.}
\bvolume{160}
\bpages{815--834}.
\end{barticle}
\endbibitem

\bibitem{CorwinShen19}
\begin{barticle}[author]
\bauthor{\bsnm{Corwin},~\bfnm{I.}\binits{I.}} \AND
  \bauthor{\bsnm{Shen},~\bfnm{H.}\binits{H.}}
\btitle{{Some recent progress in singular stochastic PDEs}}.
\bjournal{arXiv:1904.00334, to appear in Bull. Amer. Math. Soc.}
\end{barticle}
\endbibitem

\bibitem{Landscape}
\begin{barticle}[author]
\bauthor{\bsnm{Dauvergne},~\bfnm{D.}\binits{D.}},
  \bauthor{\bsnm{Ortmann},~\bfnm{J}\binits{J.}} \AND
  \bauthor{\bsnm{Vir\'{a}g},~\bfnm{B.}\binits{B.}}
(\byear{2018}).
\btitle{The directed landscape}.
\bjournal{arXiv:1812.00309}.
\end{barticle}
\endbibitem

\bibitem{Basicproperties}
\begin{barticle}[author]
\bauthor{\bsnm{Dauvergne},~\bfnm{D.}\binits{D.}} \AND
  \bauthor{\bsnm{Vir\'{a}g},~\bfnm{B.}\binits{B.}}
(\byear{2018}).
\btitle{Basic properties of the {A}iry line ensemble}.
\bjournal{arXiv:1812.00311}.
\end{barticle}
\endbibitem

\bibitem{de_Nardis_Le_Doussal_2017}
\begin{barticle}[author]
\bauthor{\bparticle{de} \bsnm{Nardis},~\bfnm{J.}\binits{J.}} \AND
  \bauthor{\bsnm{Le~Doussal},~\bfnm{P.}\binits{P.}}
(\byear{2017}).
\btitle{Tail of the two-time height distribution for {KPZ} growth in one
  dimension}.
\bjournal{J. Stat. Mech.}
\bvolume{2017}
\bpages{053212}.
\end{barticle}
\endbibitem

\bibitem{de_Nardis_2018}
\begin{barticle}[author]
\bauthor{\bparticle{de} \bsnm{Nardis},~\bfnm{J.}\binits{J.}} \AND
  \bauthor{\bsnm{Le~Doussal},~\bfnm{P.}\binits{P.}}
(\byear{2018}).
\btitle{Two-time height distribution for 1D {KPZ} growth: the recent exact
  result and its tail via replica}.
\bjournal{J. Stat. Mech.}
\bvolume{2018}
\bpages{093203}.
\end{barticle}
\endbibitem

\bibitem{De_Nardis_Le_Doussal_Takeuchi}
\begin{barticle}[author]
\bauthor{\bparticle{de} \bsnm{Nardis},~\bfnm{J.}\binits{J.}},
  \bauthor{\bsnm{Le~Doussal},~\bfnm{P.}\binits{P.}} \AND
  \bauthor{\bsnm{Takeuchi},~\bfnm{K.}\binits{K.}}
(\byear{2017}).
\btitle{Memory and Universality in Interface Growth}.
\bjournal{Phys. Rev. Lett.}
\bvolume{118}
\bpages{125701}.
\end{barticle}
\endbibitem

\bibitem{Dembo2007}
\begin{barticle}[author]
\bauthor{\bsnm{Dembo},~\bfnm{A.}\binits{A.}} \AND
  \bauthor{\bsnm{Deuschel},~\bfnm{J.~D.}\binits{J.~D.}}
(\byear{2007}).
\btitle{Aging for interacting diffusion processes}.
\bjournal{Ann. Inst. H. Poincar\'{e} B}
\bvolume{43}
\bpages{461-480}.
\end{barticle}
\endbibitem

\bibitem{Dotsenko_2013}
\begin{barticle}[author]
\bauthor{\bsnm{Dotsenko},~\bfnm{V.}\binits{V.}}
(\byear{2013}).
\btitle{Two-time free energy distribution function in $(1 + 1)$ directed
  polymers}.
\bjournal{J. Stat. Mech.}
\bvolume{2013}
\bpages{P06017}.
\bdoi{10.1088/1742-5468/2013/06/p06017}
\end{barticle}
\endbibitem

\bibitem{Dotsenko_2015}
\begin{barticle}[author]
\bauthor{\bsnm{Dotsenko},~\bfnm{V.}\binits{V.}}
(\byear{2015}).
\btitle{{Two-time free energy distribution function in the KPZ problem}}.
\bjournal{arXiv:1507.06135}.
\end{barticle}
\endbibitem

\bibitem{Dotsenko_2016}
\begin{barticle}[author]
\bauthor{\bsnm{Dotsenko},~\bfnm{V.}\binits{V.}}
(\byear{2016}).
\btitle{On two-time distribution functions in $(1 + 1)$ random directed
  polymers}.
\bjournal{J. Phys. A}
\bvolume{49}
\bpages{27LT01}.
\bdoi{10.1088/1751-8113/49/27/27lt01}
\end{barticle}
\endbibitem

\bibitem{Ferrari_2008}
\begin{barticle}[author]
\bauthor{\bsnm{Ferrari},~\bfnm{P.~L.}\binits{P.~L.}}
(\byear{2008}).
\btitle{{Slow decorrelations in Kardar{\textendash}Parisi{\textendash}Zhang
  growth}}.
\bjournal{J. Stat. Mech.}
\bvolume{07}
\bpages{P07022}.
\bdoi{10.1088/1742-5468/2008/07/p07022}
\end{barticle}
\endbibitem

\bibitem{FO18}
\begin{barticle}[author]
\bauthor{\bsnm{{Ferrari}},~\bfnm{P.~L.}\binits{P.~L.}} \AND
  \bauthor{\bsnm{{Occelli}},~\bfnm{A.}\binits{A.}}
(\byear{2019}).
\btitle{{Time-time Covariance for Last Passage Percolation with Generic Initial
  Profile}}.
\bjournal{Math. Phys. Anal. Geom.}
\bvolume{22}.
\end{barticle}
\endbibitem

\bibitem{FS16}
\begin{barticle}[author]
\bauthor{\bsnm{Ferrari},~\bfnm{P.~L.}\binits{P.~L.}} \AND
  \bauthor{\bsnm{Spohn},~\bfnm{H.}\binits{H.}}
(\byear{2016}).
\btitle{On time correlations for {KPZ} growth in one dimension}.
\bjournal{SIGMA: Symm. Integr. Geom. Meth. Appl.}
\bvolume{12}.
\end{barticle}
\endbibitem

\bibitem{GJ2014a}
\begin{barticle}[author]
\bauthor{\bsnm{Gon{\c{c}}alves},~\bfnm{P.}\binits{P.}} \AND
  \bauthor{\bsnm{Jara},~\bfnm{M.}\binits{M.}}
(\byear{2014}).
\btitle{Nonlinear fluctuations of weakly asymmetric interacting particle
  systems}.
\bjournal{Arch. Ration. Mech. Anal.}
\bvolume{212}
\bpages{597--644}.
\bdoi{10.1007/s00205-013-0693-x}
\bmrnumber{3176353}
\end{barticle}
\endbibitem

\bibitem{GIP15}
\begin{barticle}[author]
\bauthor{\bsnm{Gubinelli},~\bfnm{M.}\binits{M.}},
  \bauthor{\bsnm{Imkeller},~\bfnm{P.}\binits{P.}} \AND
  \bauthor{\bsnm{Perkowski},~\bfnm{N.}\binits{N.}}
(\byear{2015}).
\btitle{Paracontrolled distributions and singular {PDE}s}.
\bjournal{Forum Math. Pi}
\bvolume{3}
\bpages{e6, 75}.
\bmrnumber{3406823}
\end{barticle}
\endbibitem

\bibitem{GP17}
\begin{barticle}[author]
\bauthor{\bsnm{Gubinelli},~\bfnm{M.}\binits{M.}} \AND
  \bauthor{\bsnm{Perkowski},~\bfnm{N.}\binits{N.}}
(\byear{2017}).
\btitle{K{PZ} reloaded}.
\bjournal{Commun. Math. Phys.}
\bvolume{349}
\bpages{165--269}.
\bmrnumber{3592748}
\end{barticle}
\endbibitem

\bibitem{GP2015a}
\begin{barticle}[author]
\bauthor{\bsnm{{Gubinelli}},~\bfnm{M.}\binits{M.}} \AND
  \bauthor{\bsnm{{Perkowski}},~\bfnm{N.}\binits{N.}}
(\byear{2018}).
\btitle{{Energy solutions of {KPZ} are unique}}.
\bjournal{J. Amer. Math. Soc.}
\bvolume{31}
\bpages{427--471}.
\end{barticle}
\endbibitem

\bibitem{Hairer13}
\begin{barticle}[author]
\bauthor{\bsnm{Hairer},~\bfnm{M.}\binits{M.}}
(\byear{2013}).
\btitle{Solving the {KPZ} equation}.
\bjournal{Ann. of Math. (2)}
\bvolume{178}
\bpages{559--664}.
\bmrnumber{3071506}
\end{barticle}
\endbibitem

\bibitem{Hammond18c}
\begin{barticle}[author]
\bauthor{\bsnm{{Hammond}},~\bfnm{A.}\binits{A.}}
(\byear{2016}).
\btitle{{Brownian regularity for the Airy line ensemble, and multi-polymer
  watermelons in Brownian last passage percolation.}}
\bjournal{Mem. Amer. Math. Soc., to appear. arXiv:1609.02971}.
\end{barticle}
\endbibitem

\bibitem{Hammond18a}
\begin{barticle}[author]
\bauthor{\bsnm{{Hammond}},~\bfnm{A.}\binits{A.}}
(\byear{2017}).
\btitle{{Exponents governing the rarity of disjoint polymers in Brownian last
  passage percolation.}}
\bjournal{Proc. Lond. Math. Soc. to appear. arXiv:1709.04110}.
\end{barticle}
\endbibitem

\bibitem{Hammond18b}
\begin{barticle}[author]
\bauthor{\bsnm{{Hammond}},~\bfnm{A.}\binits{A.}}
(\byear{2017}).
\btitle{{Modulus of continuity of polymer weight profiles in Brownian last
  passage percolation.}}
\bjournal{Ann. Probab., to appear. arXiv:1709.04115}.
\end{barticle}
\endbibitem

\bibitem{Hammond18}
\begin{barticle}[author]
\bauthor{\bsnm{{Hammond}},~\bfnm{A.}\binits{A.}}
(\byear{2019}).
\btitle{A patchwork quilt sewn from Brownian fabric: regularity of polymer
  weight profiles in Brownian last passage percolation}.
\bjournal{Forum Math. Pi}
\bvolume{7}
\bpages{e2, 69}.
\bmrnumber{3406823}
\end{barticle}
\endbibitem

\bibitem{HammondSarkar}
\begin{barticle}[author]
\bauthor{\bsnm{{Hammond}},~\bfnm{A.}\binits{A.}} \AND
  \bauthor{\bsnm{{Sarkar}},~\bfnm{S.}\binits{S.}}
(\byear{2020}).
\btitle{{Modulus of continuity for polymer fluctuations and weight profiles in
  Poissonian last passage percolation}}.
\bjournal{Electron. J. Prob.}
\bvolume{25}
\bpages{no. 29, 38 pp.}
\end{barticle}
\endbibitem

\bibitem{Johansson2017}
\begin{barticle}[author]
\bauthor{\bsnm{Johansson},~\bfnm{K.}\binits{K.}}
(\byear{2017}).
\btitle{Two Time Distribution in Brownian Directed Percolation}.
\bjournal{Commun. Math. Phys.}
\bvolume{351}
\bpages{441--492}.
\bdoi{10.1007/s00220-016-2660-5}
\end{barticle}
\endbibitem

\bibitem{J18}
\begin{barticle}[author]
\bauthor{\bsnm{{Johansson}},~\bfnm{K.}\binits{K.}}
(\byear{2018}).
\btitle{{The two-time distribution in geometric last-passage percolation}}.
\bjournal{arXiv:1802.00729, to appear in Probab. Theo. Rel. Fields}.
\end{barticle}
\endbibitem

\bibitem{J19}
\begin{barticle}[author]
\bauthor{\bsnm{{Johansson}},~\bfnm{K.}\binits{K.}}
(\byear{2019}).
\btitle{{The Long ans Short Time Asymptotics of the Two-Time Distribution in
  Local Random Growth}}.
\bjournal{arXiv:1904.08195}.
\end{barticle}
\endbibitem

\bibitem{J19B}
\begin{barticle}[author]
\bauthor{\bsnm{{Johansson}},~\bfnm{K.}\binits{K.}} \AND
  \bauthor{\bsnm{{Rahman}},~\bfnm{M.}\binits{M.}}
(\byear{2019}).
\btitle{{Multi-time distribution in discrete polynuclear growth}}.
\bjournal{arXiv:1906.01053}.
\end{barticle}
\endbibitem

\bibitem{KPZ86}
\begin{barticle}[author]
\bauthor{\bsnm{Kardar},~\bfnm{M.}\binits{M.}},
  \bauthor{\bsnm{Parisi},~\bfnm{G.}\binits{G.}} \AND
  \bauthor{\bsnm{Zhang},~\bfnm{Y.~C.}\binits{Y.~C.}}
(\byear{1986}).
\btitle{Dynamic Scaling of Growing Interfaces}.
\bjournal{Phys. Rev. Lett.}
\bvolume{56}
\bpages{889--892}.
\bdoi{10.1103/PhysRevLett.56.889}
\end{barticle}
\endbibitem

\bibitem{Kupiainen16}
\begin{barticle}[author]
\bauthor{\bsnm{Kupiainen},~\bfnm{A.}\binits{A.}}
(\byear{2016}).
\btitle{Renormalization group and stochastic {PDE}s}.
\bjournal{Ann. H. Poincar\'{e}}
\bvolume{17}.
\end{barticle}
\endbibitem

\bibitem{D17}
\begin{barticle}[author]
\bauthor{\bsnm{Le~Doussal},~\bfnm{P.}\binits{P.}}
(\byear{2017}).
\btitle{Maximum of an Airy process plus Brownian motion and memory in
  Kardar-Parisi-Zhang growth}.
\bjournal{Phys. Rev. E}
\bvolume{96}
\bpages{060101}.
\end{barticle}
\endbibitem

\bibitem{Liggett99}
\begin{bbook}[author]
\bauthor{\bsnm{Liggett},~\bfnm{T.~M.}\binits{T.~M.}}
(\byear{1999}).
\btitle{Stochastic interacting systems: contact, voter and exclusion
  processes}.
\bseries{Grundlehren der Mathematischen Wissenschaften}
\bvolume{324}.
\bpublisher{Springer-Verlag, Berlin}.
\bdoi{10.1007/978-3-662-03990-8}
\bmrnumber{1717346}
\end{bbook}
\endbibitem

\bibitem{Liggett85}
\begin{bbook}[author]
\bauthor{\bsnm{Liggett},~\bfnm{T.~M.}\binits{T.~M.}}
(\byear{2005}).
\btitle{Interacting particle systems}.
\bseries{Classics in Mathematics}.
\bpublisher{Springer-Verlag, Berlin}
\bnote{Reprint of the 1985 original}.
\bmrnumber{2108619}
\end{bbook}
\endbibitem

\bibitem{Liu_2019}
\begin{barticle}[author]
\bauthor{\bsnm{Liu},~\bfnm{Z.}\binits{Z.}}
(\byear{2019}).
\btitle{{Multi-time distribution of TASEP}}.
\bjournal{arXiv:1907.09876}.
\end{barticle}
\endbibitem

\bibitem{KPZfixed}
\begin{barticle}[author]
\bauthor{\bsnm{Matetski},~\bfnm{K.}\binits{K.}},
  \bauthor{\bsnm{Quastel},~\bfnm{J.}\binits{J.}} \AND
  \bauthor{\bsnm{Remenik},~\bfnm{D.}\binits{D.}}
(\byear{2017}).
\btitle{The {KPZ} fixed point}.
\bjournal{arXiv:1701.00018}.
\end{barticle}
\endbibitem

\bibitem{ON12}
\begin{barticle}[author]
\bauthor{\bsnm{O'Connell},~\bfnm{N.}\binits{N.}}
(\byear{2012}).
\btitle{Directed polymers and the quantum {T}oda lattice}.
\bjournal{Ann. Probab.}
\bvolume{40}
\bpages{437--458}.
\end{barticle}
\endbibitem

\bibitem{P17}
\begin{barticle}[author]
\bauthor{\bsnm{Pimentel},~\bfnm{L.}\binits{L.}}
(\byear{2018}).
\btitle{{Local Behaviour of Airy Processes}}.
\bjournal{J. Stat. Phys.}
\bvolume{173}
\bpages{1614--1638}.
\end{barticle}
\endbibitem

\bibitem{Prahofer2002}
\begin{barticle}[author]
\bauthor{\bsnm{Pr{\"a}hofer},~\bfnm{M.}\binits{M.}} \AND
  \bauthor{\bsnm{Spohn},~\bfnm{H.}\binits{H.}}
(\byear{2002}).
\btitle{{Scale Invariance of the PNG Droplet and the Airy Process}}.
\bjournal{J. Stat. Phys.}
\bvolume{108}
\bpages{1071--1106}.
\end{barticle}
\endbibitem

\bibitem{Takeuchi2012}
\begin{barticle}[author]
\bauthor{\bsnm{Takeuchi},~\bfnm{K.}\binits{K.}} \AND
  \bauthor{\bsnm{Sano},~\bfnm{M.}\binits{M.}}
(\byear{2012}).
\btitle{{Evidence for Geometry-Dependent Universal Fluctuations of the
  Kardar-Parisi-Zhang Interfaces in Liquid-Crystal Turbulence}}.
\bjournal{J. Stat. Phys.}
\bvolume{147}
\bpages{853--890}.
\bdoi{10.1007/s10955-012-0503-0}
\end{barticle}
\endbibitem

\bibitem{TracyWidom94}
\begin{barticle}[author]
\bauthor{\bsnm{Tracy},~\bfnm{C.~A.}\binits{C.~A.}} \AND
  \bauthor{\bsnm{Widom},~\bfnm{H.}\binits{H.}}
(\byear{1994}).
\btitle{Level-spacing distributions and the {A}iry kernel}.
\bjournal{Commun. Math. Phys.}
\bvolume{159}
\bpages{151--174}.
\bmrnumber{1257246}
\end{barticle}
\endbibitem

\bibitem{Walsh}
\begin{bbook}[author]
\bauthor{\bsnm{Walsh},~\bfnm{J.~B.}\binits{J.~B.}}
(\byear{1986}).
\btitle{An introduction to stochastic partial differential equations}.
\bseries{Ecole d'Et\'{e} de Probabilit\'{e}s de Saint-Flour XIV - 1984}.
\bpublisher{Springer-Verlag, Berlin}.
\end{bbook}
\endbibitem

\end{thebibliography}

\end{document}